\documentclass{article}

\usepackage{amssymb}
\newcommand{\PSbox}[3]{\mbox{\rule{0in}{#3}\includegraphics{#1}\hspace{#2}}} 
\newtheorem{Theorem}{Theorem}[section]
\newtheorem{Definition}{Definition}[section]

\newtheorem{Convention}[Theorem]{Convention}
\newtheorem{Lemma}[Theorem]{Lemma}
\newtheorem{Proposition}[Theorem]{Proposition}

\newtheorem{Corollary}[Theorem]{Corollary}

\newtheorem{Claim}[Theorem]{Claim}
\newtheorem{Fact}{Fact}[section]

\def\sqr#1#2{{\vcenter{\vbox{\hrule  height.#2pt
       \hbox{\vrule width.#2pt height#1pt \kern#1pt \vrule width.#2pt}
        \hrule height.#2pt}}}}
\def\bb{\sqr66}   
\let\epsilon=\varepsilon
\let\Bbb=\mathbb
\let\phi=\varphi
\def\){ \right) }
\def\({ \left( }
\def\[{ \left[ }
\def\]{ \right] }
\def\<{ \langle }
\def\>{ \rangle }
\let\ljunk=\{
\let\rjunk=\}
\def\{{\left\ljunk}
\def\}{\right\rjunk}
\def\p{\partial}

\def\Riem{{\cal R}{\mathrm i}{\mathrm e}{\mathrm m}}
\def\dist{{\mathrm d}{\mathrm i}{\mathrm s}{\mathrm t}}

\def\cyl{{\mathit c}{\mathit y}\ell}

\def\can{{\mathrm c}{\mathrm a}{\mathrm n}}
\def\const{{\mathrm c}{\mathrm o}{\mathrm n}{\mathrm s}{\mathrm t}}
\def\diam{{\mathrm d}{\mathrm i}{\mathrm a}{\mathrm m}}

\newcommand{\supp}{{\mathbf S}{\mathrm u}{\mathrm p}{\mathrm p}}
\newcommand{\fla}{{\mathrm f}{\mathrm l}{\mathrm a}{\mathrm t}}

\def\lr{\mbox{{\Large $\lrcorner$}}}
\def\L{\Bbb L}
\def\Vol{{\mathrm V}{\mathrm o}{\mathrm l}}
\def\Ric{\mbox{Ric}}

\def\ov{\overline}

\def\intt{{\mathrm I}{\mathrm n}{\mathrm t}}

\newcommand{\Cone}{{\cal C}{\mathrm one}}

\newcommand{\R}{{\mathbf R}}

\newcommand{\Z}{{\mathbf Z}}
\newcommand{\inj}{{\mathrm i}{\mathrm n}{\mathrm j}}

\def\P{\mbox{\bf P}}
\newenvironment{Proof}
  {\begin{trivlist}\item[\hskip \labelsep{\bf Proof.}]}
  {\hfill $\Box$\end{trivlist}}
\newenvironment{lproof}[1]{\par\vspace{2mm} \noindent{\bf
#1.} }{\hfill $\bb$\par\vspace{2mm}}
\begin{document} 
\title{Yamabe metrics on cylindrical manifolds} \author{Kazuo
Akutagawa$^*$, Boris Botvinnik\thanks{Partially supported by the
Grants-in-Aid for Scientific Research (C), Japan Society for the
Promotion of Science, No. 14540072.}}  \date{ \ } \maketitle
\markboth{Yamabe metrics on cylindrical manifolds}{K. Akutagawa,
B. Botvinnik, Yamabe metrics on cylindrical manifolds}
\vspace{-10mm}

\begin{abstract}
We study a particular class of open manifolds. In the category of
Riemannian manifolds these are complete manifolds with cylindrical
ends. We give a natural setting for the conformal geometry on such
manifolds including an appropriate notion of the \emph{cylindrical
Yamabe constant/invariant}. This leads to a corresponding version of
the \emph{Yamabe problem on cylindrical manifolds}.  We affirmatively
solve this Yamabe problem: we prove the existence of minimizing
metrics and analyze their singularities near infinity.  These
singularities turn out to be of very particular type: either
\emph{almost conical} or \emph{almost cusp} singularities.  We
describe the \emph{supremum case}, i.e. when the cylindrical Yamabe
constant is equal to the Yamabe invariant of the sphere. We prove that
in this case such a cylindrical manifold coincides conformally with
the standard sphere punctured at a finite number of points. In the
course of studying the supremum case, we establish a \emph{Positive
Mass Theorem for specific asymptotically flat manifolds with two
almost conical singularities}. As a by-product, we revisit known
results on \emph{surgery and the Yamabe invariant}. {\sc Key
words:} \emph{manifolds with cylindrical ends, Yamabe
constant/invariant, Yamabe problem, conical metric singularities, cusp
metric singularities, Positive Mass Theorem, surgery and Yamabe
invariant.}
\end{abstract}

\section{Introduction}\label{int}
The goal of this paper is to formulate and solve a natural version of
the Yamabe problem for complete manifolds with cylindrical ends.
Before describing our main results in detail, we recall the classical
situation.
%

\noindent
{\bf \ref{int}.1. Classical Yamabe problem for compact manifolds.}
Let $M$ be a smooth closed manifold (i.e. a smooth compact manifold
without boundary) of $\dim M = n\geq 3$ and $\Riem(M)$ the space of
all Riemannian metrics on $M$.  We denote by $R_g$ the scalar
curvature and by $d\sigma_{g}$ the volume form for each Riemannian
metric $g\in\Riem(M)$. Then the (normalized) Einstein-Hilbert
functional $I_M : \Riem(M) \to \R$ is defined as
\begin{equation}\label{int:eq1}
I_M : g \mapsto  \frac{\int_M R_{g} d\sigma_{g}}
{\Vol_{g}(M)^{{n-2\over n}}}.
\end{equation}
The classical Yamabe problem is to find a metric $\check{g}$ in a
given conformal class $C$ such that the Einstein-Hilbert functional
attains its minimum on $C$: $\displaystyle I_M(\check{g})=\inf_{g\in
C} I_M(g)=: Y_C(M)$. This minimizing metric $\check{g}$ is called a
Yamabe metric, and the conformal invariant $Y_C(M)$ the Yamabe
constant.

It is a celebrated result in conformal geometry that the Yamabe
problem has an affirmative solution for closed manifolds. This was
proven in a series of papers starting with the work of Yamabe
\cite{Ya}. Although Yamabe's proof overlooked the fundamental analytic
difficulty concerning Sobolev inequalities with a critical exponent,
the general strategy in \cite{Ya} was correct.  This difficulty is a
profound one, and the proof was ultimately only corrected in stages,
first by Trudinger \cite{Tr}, then by Aubin \cite{Au2}, and was
finally completed by Schoen \cite{Sc1}. Text-book style proofs are
now available in \cite{SY4} and \cite{LP}.

The second case when the Yamabe problem has a positive solution is
when a manifold $M$ has a non-empty boundary. In this case, the
Einstein-Hilbert functional has to be restricted to metrics (and
corresponding conformal subclasses) with a minimal boundary condition
(see \cite{AB1}).  The resulting Yamabe problem was solved by Cherrier
\cite{Ch} and Escobar \cite{Es} under some mild
restrictions.
\vspace{2mm}

\noindent
{\bf \ref{int}.2. Yamabe problem for open manifolds.} In the case of
open manifolds, it is important to clarify what is a \emph{suitable
version of the Yamabe problem} which would capture the geometry of
their ends. One reasonable version of the Yamabe problem is to find a
\emph{complete metric of constant scalar curvature in a given
conformal class}. In this case, however, the concept of a minimizing
Yamabe metric does not make sense, and there are serious difficulties
in this area. Indeed, there are simple noncompact Riemannian manifolds
$(N,\bar{g})$ for which there does not exist \emph{any complete metric
(conformal to $\bar{g}$) of constant scalar curvature},
(cf. \cite{McO}). Hence one needs to place some geometric restriction
(e.g. curvature or injectivity radius bounds) near each end of a
noncompact manifold in order to establish the existence of a metric
of constant scalar curvature.

A more specific version of the Yamabe problem for open manifolds is
the \emph{singular Yamabe problem}, i.e.  when an open manifold $N$ is
a complement $N=M\setminus \Sigma$ of a submanifold $\Sigma$ of a
closed Riemannian manifold $(M,g)$. Then the singular Yamabe problem
is to find a complete metric $\tilde g$ that is conformal to $g$ on
$M\setminus\Sigma$ and has constant scalar curvature, see \cite{AM},
\cite{Sc2} for earlier results and also \cite{McO} for a survey
of results in this area. Here again the concept of a minimizing
Yamabe metric is not well-defined.
\vspace{2mm}

\noindent
{\bf \ref{int}.3. Cylindrical manifolds.}  In this paper, we consider
a particular class of open manifolds, which we call \emph{cylindrical
manifolds}. These are open manifolds with \emph{tame ends} equipped
with a \emph{cylindrical metric} on each end. Cylindrical manifolds
are well-known objects in geometry and topology.  First of all, these
manifolds are well-suited for the Dirac operator on complete
manifolds, as was shown first by Gromov and Lawson in
\cite{GL2}. Secondly, these objects are well known in gauge theory,
where cylindrical manifolds have been thoroughly studied, see for
example the books \cite{FU} and \cite{MMR,N,Ta}.  We
remark that $\Z/k$-manifolds (as geometric objects) can be thought of
as cylindrical manifolds with $k$ identical cylindrical ends (see
\cite{Bo} for results on the existence of positive scalar
curvature on $\Z/k$-manifolds).

There are several other constructions in geometry where cylindrical
manifolds show up as natural \emph{limit manifolds}. The first
example comes from the work describing the process of \emph{bubbling
out of Einstein manifolds}, see \cite{An}, \cite{BKN} and
\cite{Ban}. Here it is known that for a sequence of Einstein
manifolds $(X_i,g_i)$ (with some natural restrictions on their
geometry) there exists a subsequence $(X_{i_q},g_{i_q})$ which converges
(in the Gromov-Hausdorff topology) to a compact \emph{Einstein
orbifold} $(X_{\infty},g_{\infty})$ with a finite set of singular
points. Let $U_p$ be an open neighborhood of a singular point $p\in
X_{\infty}$. Then the tangent cone of $U_p\setminus \{p\}$ at $p$ is
conformally cylindrical, although $U_p\setminus \{p\}$ itself need not
be.  Thus a cylindrical manifold (with an appropriate conformal
metric) can be thought of as a \emph{linear approximation} of the
Einstein orbifold $(X_{\infty},g_{\infty})$ (see \cite{AB3} for the
relation of this approximation to the Yamabe invariant of orbifolds).

The second example is related to the Yamabe invariant $Y(M)=\sup_C
Y_C(M)$ of a compact manifold $M$, see \cite{K2} and
\cite{Sc3}. It was shown that $Y(S^{n-1}\times S^1)= Y(S^n)$, i.e.
there exists a sequence of conformal classes $C_i$ and Yamabe metrics
$\check{g}_i\in C_i$ such that the limit $\lim_i
Y_{C_i}(S^{n-1}\times S^1) = Y(S^n)$.  The sequence of Riemannian
manifolds $( S^{n-1}\times S^1, \check{g}_i)$ converges to the
standard sphere $S^n(1)$ identified with two antipodal points, see
\cite{K1}, \cite{Sc3}.  After deleting the singular point,
the punctured sphere is conformally equivalent to the canonical
cylindrical manifold $S^{n-1} \times {\mathbf R}$.  In the general
case of a closed manifold $M$ with positive Yamabe invariant, there
is also a compactness result.  In \cite{Ak1},
\cite{Ak2}, the first author proved that, under some
restrictions, for a sequence of Yamabe metrics $\check{g}_i$ on $M$
satisfying $\lim_i Y_{[\check{g}_i]}(M) = Y(M)$ there exists a
subsequence $\check{g}_{i_{\nu}}$ of $\check{g}_i$ such that
$(M,\check{g}_{i_{\nu}})$ converges (in the Gromov-Hausdorff topology)
to a compact metric space $(M_{\infty},g_{\infty})$, which is a smooth
Riemannian manifold away from a finite number of singular points.
Again, a cylindrical manifold (with an appropriate conformal metric)
serves here as a \emph{linear approximation} of the limit singular
space $(M_{\infty},g_{\infty})$.

For a given closed manifold $M$, it is a challenging problem to find a
nice sequence of Yamabe metrics $\check{g}_i$ satisfying \ $\lim_i
Y_{[\check{g}_i]}(M) = Y(M)$ \ such that \ $(M,\check{g}_{i})$ \
converges to a singular Riemannian space $(M_{\infty}, g_{\infty})$,
and to understand the limit singular space $(M_{\infty},g_{\infty})$
with bubbling out spaces. In our view, cylindrical manifolds may
provide a typical conformal model of such a limit singular space
$(M_{\infty},g_{\infty})$.
\vspace{2mm}

\noindent
{\bf \ref{int}.4.  Yamabe problem  for cylindrical manifolds.}  In the
smooth category, a cylindrical manifold $X$ is an open manifold with a
relatively compact  open submanifold $W\subset  X$ with $\displaystyle
\p  \overline{W} =  Z$  such that  $\displaystyle  X\setminus W  \cong
Z\times [0,\infty)$, where $Z=\bigsqcup_{j=1}^m Z_j$ and each $Z_j$ is
a connected closed $(n-1)$-manifold. Throughout this paper, we always
assume that $X$ is connected.

To formulate an appropriate notion of conformal class on cylindrical
manifolds requires some care. As a \emph{reference metric}, we start
with a \emph{cylindrical} Riemannian metric $\bar{g}$ on $X$, i.e. $
\bar{g}(x,t) = h(x) + dt^2$ on $Z\times [1,\infty)$, for some metric
$h$ on $Z$. We denote $\p_{\infty}\bar{g}=h$. Let $[\bar{g}]$ denote
the conformal class containing $\bar{g}$. Clearly the Einstein-Hilbert
functional $I_X$ (given by (\ref{int:eq1})) restricted to the
conformal class $[\bar{g}]$ is not well-defined since the scalar
curvature may not be integrable and the volume may be
infinite. Consider the Dirichlet form
$$
Q_{(X,\bar{g})}(u) = \frac{\int_X\[\frac{4(n-1)}{n-2}|du|^2 +
R_{\bar{g}}u^2\] d\sigma_{\bar{g}}}{\(\int_X
|u|^{\frac{2n}{n-2}}d\sigma_{\bar{g}}\)^{\frac{n-2}{n}}}
$$
on the space $C^{\infty}_c(X)$ (of smooth functions on $X$ with
compact support) associated with the conformal Laplacian of
$(X,\bar{g})$. The functional $Q_{(X,\bar{g})}$ suggests the following
natural setting for the Yamabe problem on cylindrical manifolds.

First, we define the $L^{k,2}_{\bar{g}}$-conformal class
$[\bar{g}]_{L^{k,2}_{\bar{g}}}$ consisting of all metrics
$u^{\frac{4}{n-2}}\cdot \bar{g}$, where $u$ is a smooth positive
function on $X$ and $u\in L^{k,2}_{\bar{g}}(X)$, $k=1,2$.  Here $
L^{k,2}_{\bar{g}}(X)$ denotes the Sobolev space of square-integrable
functions on $X$ (with respect to $\bar{g}$) up to their $k$-th weak
derivatives. Then the functional $Q_{(X,\bar{g})}$ is well-defined on
the space $C^{\infty}_+(X)\cap L^{1,2}_{\bar{g}}(X)$, where
$C^{\infty}_+(X)= \{u\in C^{\infty}(X) \ | \ u>0 \}$. On the other
hand, the Einstein-Hilbert functional $I_X$ on the conformal class
$[\bar{g}]_{L^{2,2}_{\bar{g}}(X)}$ is also well-defined.  It turns out
that the infima of both functionals coincide:
$$
Y^{\cyl}_{[\bar{g}]}(X):=\inf_{\tilde{g}\in [\bar{g}]_{L^{2,2}_{\bar{g}}(X)}} 
I_X(\tilde{g}) = \inf_{u\in C_+^{\infty}(X)\cap
L^{1,2}_{\bar{g}}(X)} Q_{(X,\bar{g})}(u).
$$
We   call    the   constant   \    $Y^{\cyl}_{[\bar{g}]}(X)$   the
\emph{cylindrical Yamabe  constant of $(X,[\bar{g}])$} and  show that it
does not depend on the choice of a reference cylindrical metric in the
same  conformal class  as $\bar{g}$.  Now we  are ready  to  state the
Yamabe problem on cylindrical manifolds.
\vspace{1mm}

\noindent
{\bf Yamabe Problem.} {\it Given a cylindrical metric $\bar{g}$ on
$X$, does there exist a metric $\check{g}= u^{\frac{4}{n-2}} \cdot \bar{g} \in
[\bar{g}]_{L^{1,2}_{\bar{g}}}$ such that $\displaystyle
Q_{(X,\bar{g})}(u)=Y^{\cyl}_{[\bar{g}]}(X)$?}
\vspace{2mm}

\noindent
We call such a metric \ $\check{g}$ \ (if it exists) a \emph{Yamabe
metric} and such a function $u$ a \emph{Yamabe minimizer}. We
affirmatively solve the Yamabe problem for generic
cylindrical conformal classes in this setting.
\vspace{2mm}

\noindent
{\bf \ref{int}.5. The invariant $\lambda({\cal L}_h)$.} First of all,
we show that the constant $Y^{\cyl}_{[\bar{g}]}(X)$ satisfies the
bounds $ -\infty \leq Y^{\cyl}_{[\bar{g}]}(X) \leq Y(S^n), $ where
$Y(S^n)$ is the Yamabe invariant of the $n$-sphere. In particular, it
is possible that $Y^{\cyl}_{[\bar{g}]}(X)=-\infty$. To determine when
the constant $Y^{\cyl}_{[\bar{g}]}(X)$ is finite, we introduce a new
invariant $\lambda({\cal L}_h)$.  We introduce the operator
$$
{\cal L}_h := -\frac{4(n-1)}{n-2}\Delta_h + R_h \ \ \ \mbox{on} \ \ (Z,h). 
$$
Notice that the operator \ ${\cal L}_h$ \ is different from the conformal
Laplacian $\L_{h}= -\frac{4(n-2)}{n-3}\Delta_{h} +  R_{h}$ of
$(Z,h)$. Let $\lambda({\cal L}_h)$ be the first eigenvalue of the
operator ${\cal L}_h$. The invariant $\lambda({\cal L}_h)$ determines
the finiteness of the cylindrical Yamabe constant as follows:
\begin{enumerate}
\item[$\bullet$] If $\lambda({\cal L}_h)<0$, then
$Y_{\bar{C}}^{\cyl}(X)= -\infty$. In particular, if $R_h< 0$ on $Z$, then
$Y_{\bar{C}}^{\cyl}(X)= -\infty$.
\item[$\bullet$] If $\lambda({\cal L}_h)\geq 0$, then
$Y_{\bar{C}}^{\cyl}(X)>-\infty$. In particular, if $R_h\geq 0$ on $Z$, then
$Y_{\bar{C}}^{\cyl}(X)>-\infty$.
\item[$\bullet$] If $\lambda({\cal L}_h)=  0$, then
$Y_{\bar{C}}^{\cyl}(X)\leq 0$. In particular, if $R_h\equiv 0$ on $Z$, then
$Y_{\bar{C}}^{\cyl}(X)\leq 0$.
\end{enumerate}
Here we remark that if $Z$ is not connected (i.e. $Z=\bigsqcup_{j=1}^m
Z_j$ for $m\geq 2$), then $\lambda({\cal L}_h)=\displaystyle
\min_{1\leq j\leq m} \lambda({\cal L}_{h|_{Z_j}})$. Clearly the Yamabe
problem does not make sense if $\lambda({\cal L}_h)<0$. Hence the case
we study here is when $\lambda({\cal L}_h)\geq 0$.  We also observe
that, in general, if $\lambda({\cal L}_h)=0$ then there is no solution
of the Yamabe problem.
\vspace{1mm}

\noindent
{\bf \ref{int}.6. Solution of the Yamabe problem in a ``generic
case''.}  Our first result repeats, in some sense, the classical
approach to the Yamabe problem. Recall that in the ``generic case''
there, i.e. when the Yamabe constant $Y_C(M)< Y(S^n)$, there is a
well-known technique giving a solution, see \cite{Au3}.  Here,
instead of the standard sphere, we have a \emph{canonical cylindrical
manifold} $(Z\times \R, h+dt^2)$ (for some metric $h$ on $Z$) which
plays a similar role. In fact, for any cylindrical manifold
$(X,\bar{g})$ with the cylindrical end $(Z\times [1,\infty),h+dt^2)$,
we show that the inequality $Y_{[\bar{g}]}^{\cyl}(X)\leq
Y^{\cyl}_{[h+dt^2]}(Z\times
\R)$ always holds. Here we also remark that if $Z$ is not connected
(i.e. $Z=\bigsqcup_{j=1}^m Z_j$ for $m\geq 2$), then
$$
Y^{\cyl}_{[h+dt^2]}(Z\times \R) = \displaystyle \min_{1\leq j\leq m}
Y^{\cyl}_{[h|_{Z_j} +dt^2]}(Z_j\times \R).
$$
First we study the case when $Y_{\bar{C}}^{\cyl}(X)<
Y^{\cyl}_{[h+dt^2]}(Z\times \R)$. We emphasize that the situations
when the invariant $\lambda({\cal L}_h)$ is positive or zero are very
different geometrically. In fact, $Y^{\cyl}_{[h+dt^2]}(Z\times \R)>0$
(resp. $Y^{\cyl}_{[h+dt^2]}(Z\times \R)=0$) if $\lambda({\cal L}_h)>0$
(resp. $\!\lambda({\cal L}_h)=0 $), and see Theorem B below.  However the
existence results are similar.
\vspace{1.5mm}

\noindent
{\bf Theorem A.} {\it Let $X$ be an open manifold of $\dim X\geq 3$
with tame ends $Z\times [0,\infty)$, and $h\in \Riem(Z)$. Let
$\bar{C}$ be a conformal class on $X$ containing a cylindrical metric
$\bar{g}\in \bar{C}$ with $\p_{\infty} \bar{g}=h$.  Assume that {\bf
either}
\begin{enumerate}
\item[{\bf (a)}] $\lambda({\cal
L}_h)>0$ and $ Y_{\bar{C}}^{\cyl}(X)< Y^{\cyl}_{[h+dt^2]}(Z\times
\R)$, {\bf or}
\item[{\bf (b)}] $\lambda({\cal
L}_h)=0$ and $Y^{\cyl}_{\bar{C}}(X)<0$.
\end{enumerate}
Then the Yamabe problem has a solution, i.e. there exists a Yamabe
minimizer $u\in C_+^{\infty}(X)\cap L^{1,2}_{\bar{g}}(X)$ with $\int_X
u^{\frac{2n}{n-2}} d\sigma_{\bar{g}}=1$ such that $Q_{(X,\bar{g})}(u)
= Y_{\bar{C}}^{\cyl}(X)$. In particular, the minimizer $u$ satisfies
the Yamabe equation $\L_{\bar{g}}u =
Y_{\bar{C}}^{\cyl}(X)u^{\frac{n+2}{n-2}}$.}
\vspace{2mm}

Next we study the behavior near infinity of the Yamabe metrics
$\check{g} = u^{\frac{4}{n-2}}\cdot \bar{g}$ given by Theorem A.  As
we mentioned, the cases $\lambda({\cal L}_h)>0$ and $\lambda({\cal
L}_h)=0$ lead to completely different geometric situations; in
particular, the asymptotics of the Yamabe metrics turn out to be
qualitatively distinct.

The first case when $\lambda({\cal L}_h)>0$ leads to \emph{almost
conical metrics}.  A canonical geometric model here is an open cone
over $Z$, i.e.
$$
\Cone(Z)\cong (Z\times (0,\infty), e^{-2t}(h + dt^2)).
$$
Then, a metric $g$ on $X$ is \emph{almost conical on a
connected end} $Z_j\times [1,\infty)$ if on the end $Z_j\times
[1,\infty)$ it is given as $g(x,t)=\phi(x,t)(h(x)+dt^2)$, where
$\phi(x,t)$ is asymptotically bounded by $C_1\cdot e^{-\beta t}\leq
\phi(x,t) \leq C_2\cdot e^{-\alpha t}$ for some constants
$0<\alpha\leq\beta$, $0<C_1\leq C_2$.

The case when $\lambda({\cal L}_h)=0$ leads to \emph{almost cusp
metrics}. A canonical geometric model here is given by the cusp end of
a hyperbolic $n$-manifold
$$
\left(
(R^{n-1}/\Gamma)\times [1,\infty),\frac{1}{t^2}(h_{0} +dt^2)
\right)
$$ 
of curvature $-1$.  Here $(R^{n-1}/\Gamma, h_{0})$ is a closed
Riemannian manifold uniformized by a flat torus \ $T^{n-1}$. \ Then, a
metric $g$ on $X$ is \emph{almost cusp on a connected end}
$Z_j\times [1,\infty)$ if on the end $Z_j\times [1,\infty)$ it is
given as $g(x,t)=\phi(x,t)(h(x)+dt^2)$, \ where the function \
$\phi(x,t)$ is asymptotically bounded by $C_1\cdot t^{-2}\leq
\phi(x,t) \leq C_2\cdot t^{-2}$ for some constants $0<C_1\leq C_2$.
\vspace{1mm}

\noindent
{\bf Theorem B.} {\it Let $X$ be an open manifold of $\dim X\geq 3$
with tame ends $Z\times [0,\infty)$, and $h\in \Riem(Z)$. Let
$\bar{C}$ be a conformal class on $X$ containing a cylindrical metric
$\bar{g}\in \bar{C}$ with $\p_{\infty} \bar{g}=h$.
\begin{enumerate}
\item[{\bf (a)}] If $\lambda({\cal L}_h)>0$ and $
Y_{\bar{C}}^{\cyl}(X)< Y^{\cyl}_{[h+dt^2]}(Z\times \R)$, then the
Yamabe metric $\check{g}=u^{\frac{4}{n-2}}\cdot \bar{g}$ given by
Theorem A is almost conical on each connected end.
\item[{\bf (b)}] If $\lambda({\cal L}_h)=0$ and
$Y^{\cyl}_{\bar{C}}(X)<0$, then the Yamabe metric
$\check{g}=u^{\frac{4}{n-2}}\cdot \bar{g}$ given by Theorem A is
almost cusp on each connected end $Z_j\times [1,\infty)$ with
\ $\lambda({\cal L}_{h|_{Z_j}})=0$ \ and \ almost conical \ on \ each connected
end $Z_k\times [1,\infty)$ with $\lambda({\cal L}_{h|_{Z_j}})>0$.
\end{enumerate}
}
\noindent
{\bf \ref{int}.7. Solution of the Yamabe problem for canonical
cylindrical manifolds.} \ \ As we mentioned, \ the \ \emph{canonical
cylindrical \ manifolds} $(Z\times \R, \bar{h}=h+dt^2)$ play a crucial
role with respect to the Yamabe problem, and these manifolds require a
special treatment. To state the result, we need the dominant canonical
cylindrical manifold $( S^{n-1}\times \R, h_++dt^2)$, where $h_+$ is
the standard metric of constant scalar curvature $1$ on the sphere
$S^{n-1}$. We remark here that the cylindrical Yamabe constant
$Y^{\cyl}_{[h_++dt^2]}( S^{n-1}\times \R)$ coincides with the Yamabe
invariant $Y(S^n)$.
\vspace{2mm}

\noindent
{\bf Theorem C.} {\it Let $(Z\times \R, \bar{h}=h+dt^2)$ be a
connected canonical cylindrical manifold of $\dim (Z\times \R) =n \geq
3$ with $\lambda({\cal L}_h)>0$. Assume that $ Y^{\cyl}_{[\bar{h}]}(
Z\times \R) < Y^{\cyl}_{[h_++dt^2]}( S^{n-1}\times \R)$. Then the
Yamabe problem has a solution, i.e. there exists a Yamabe minimizer $
u\in C_+^{\infty}(Z\times \R)\cap L_{\bar{h}}^{1,2}(Z\times \R)$ with
$ \int_{Z\times \R} u^{\frac{2n}{n-2}}d\sigma_{\bar{h}}= 1 $ such that
$ Q_{(Z\times \R, \bar{h})}(u) = Y^{\cyl}_{[\bar{h}]}( Z\times
\R)$. Furthermore, the Yamabe metric $\check{g}=u^{\frac{4}{n-2}}\cdot
\bar{g}$ is an almost conical metric.}
\vspace{2mm}

\noindent
If $\lambda({\cal L}_h)=0$, then the Yamabe problem does not have a
solution for any canonical cylindrical manifold $(Z\times \R,
\bar{h}=h+dt^2)$, see Proposition \ref{yam:new1-1}.
\vspace{2mm}

\noindent
{\bf \ref{int}.8. Characterization of the supremum case.}  First we
analyze the sumpremum case for a connected canonical cylindrical
manifold, i.e. when the cylindrical Yamabe constant
$Y^{\cyl}_{[\bar{h}]}( Z\times \R) = Y^{\cyl}_{[h_++dt^2]}(
S^{n-1}\times \R)= Y(S^n)$. 
\vspace{2mm}

\noindent
{\bf Theorem D.} {\it Let $(Z\times \R, \bar{h}=h+dt^2)$ be a
connected canonical cylindrical manifold of $\dim (Z\times
\R) =n \geq 3$. Assume that $ Y^{\cyl}_{[\bar{h}]}( Z\times \R) =
Y(S^n).  $ Then $(Z,h)$ is homothetic to the standard sphere
$S^{n-1}(1)=(S^{n-1},h_+)$.}
\vspace{2mm}

\noindent
The proof of Theorem D is somewhat involved. First, as in the compact
manifolds case, we prove Theorem D under the condition that the
manifold $(\!Z\times\R, \bar{h})\!$ is locally conformally flat. The
remaining part of the proof is to show that $
Y^{\cyl}_{[\bar{h}]}(Z\!\times \!\R) < Y(S^n)$ provided that $(Z\times
\R, \bar{h})$ is not locally conformally flat.  This result splits
into two different cases: $n\geq 6$ and $n=3,4,5$. In the case $n\geq
6$, we use the conformal normal coordinates technique and the family
of
\emph{instantons} $u_{\epsilon}(x) =\left(
\frac{\epsilon}{\epsilon^2+|x|^2}
\right)^{\frac{n-2}{2}}$ on $\R^n$ to construct nice test functions and to
prove the desired statement.

The case when $n=3,4,5$ is a truly difficult part.  We first use a
technique similar to the compact manifolds case to construct a
\emph{minimal positive Green function} $\bar{G}_p\in
C^{\infty}((Z\times \R)\setminus \{p\})$ of the conformal Laplacian
$\L_{\bar{h}}$ on $Z\times \R$. Then $\bar{G}_p$ has the expansion
$\bar{G}_p = r^{2-n}+ A + O^{\prime\prime}(r)$ in the conformal normal
coordinates, where $A$ is a constant related to the \emph{mass}.  To
complete the proof of Theorem D, we analyze the following situation.

Let $(Z\times \R, h+dt^2)$ be a connected canonical cylindrical
manifold with $\lambda({\cal L}_h)>0$, and $p\in Z\times \R$. Then for
a point $p\in Z\times \R$, we consider the manifold $\hat{X} =
(Z\times \R)\setminus \{p\}$ with the metric $\hat{h}=
\bar{G}_p^{\frac{4}{n-2}}\cdot \bar{h}$, where $\bar{G}_p$ is the
minimal positive Green function. We obtain that $(\hat{X},\hat{h})$ is
a scalar-flat, asymptotically flat manifold.  We emphasize that
$(\hat{X},\hat{h})$ has two almost conical singularities. Then we note
that in this case the \emph{mass} ${\mathfrak m}(\hat{h})$ is
well-defined in the same way as in the classical case, see
\cite{Bar}. The following Theorem is a generalization of the
classical positive mass theorem \cite{SY4,SY1,SY2} (cf. \cite{Bar}
and \cite{LP}) to the case of specific asymptotically flat manifolds
with two almost conical singularities.
\vspace{2mm}

\noindent
{\bf Theorem E.} (Positive mass theorem) {\it Let $(Z\times \R,
\bar{h}=h+dt^2)$ be a connected canonical cylindrical manifold of $\dim
(Z\times \R)=n$ {\rm (}$n =3,4,5${\rm )} with $\lambda({\cal L}_h)>0$,
and $p\in Z\times \R$. Let $(\hat{X},\hat{h}) = ((Z\times \R)\setminus
\{p\}, \bar{G}_p^{\frac{4}{n-2}}\cdot \bar{h})$ be the scalar-flat,
asymptotically flat manifold with two almost conical singularities as above.
Then the mass ${\mathfrak m}(\hat{h})$ is non-negative.  Furthermore,
if ${\mathfrak m}(\hat{h})=0$, then $(\hat{X},\hat{h})$ is isometric
to $\R^n\setminus\{2 \ {\mathrm p}{\mathrm o}{\mathrm i}{\mathrm
n}{\mathrm t}{\mathrm s}\}$ with the Euclidian metric.}

We use Theorem D to analyze the supremum case for general cylindrical
manifolds, i.e. when $Y^{\cyl}_{[\bar{g}]}(X)= Y(S^n)$.
\vspace{2mm}

\noindent
{\bf Theorem F.} {\it Let $(X,\bar{g})$ be a cylindrical manifold of
$\dim X \geq 3$ with tame ends $Z\times [0,\infty)$ and $\p_{\infty}
\bar{g}=h$. Set $Z =\bigsqcup_{i=1}^m Z_j$, where each $Z_j$ is 
connected. Assume that $ Y^{\cyl}_{[\bar{g}]}(X)= Y(S^n)$.  Then
there exist $k$ points $\{p_1,\ldots, p_m\}$ in $S^n$ such that
\begin{enumerate}
\item[{\bf (i)}] each manifold $(Z_j,h_j)$ is homothetic to
$(S^{n-1},h_+)$,
\item[{\bf (ii)}] the manifold $(X,[\bar{g}])$ is conformally
equivalent to the punctured sphere $(S^n \setminus\!
\{p_1,\ldots,p_m\}, C_{\can})$.
\end{enumerate}
Here $h_j=h|_{Z_j}$ and $C_{\can}$ denotes the canonical conformal class on
$S^n$.
}
\vspace{2mm}

\noindent
{\bf \ref{int}.9. Cylindrical Yamabe invariant and gluing
constructions.} The Yamabe problem on cylindrical manifolds and the
cylindrical Yamabe constant naturally leads us to the Yamabe invariant
for cylindrical manifolds. For each metric $h$ on the slice $Z$, we
define the $h$-\emph{cylindrical Yamabe invariant} $Y^{h\hbox{-}\cyl}(X)$ and
the \emph{cylindrical Yamabe invariant $Y^{\cyl}(X)$} as follows:
$$
Y^{h\hbox{-}\cyl}(X):=\sup_{^{^{\begin{array}{c}
_{\bar{g}\in\Riem^{\cyl}(X)}\\ _{\p_{\infty}\bar{g}=h}
\end{array}}}} Y_{[\bar{g}]}^{\cyl}(X), \ \ \ \ 
Y^{\cyl}(X):=\sup_{h\in \Riem(Z)} Y^{h\hbox{-}\cyl}(X).
$$
We show that the $h$-cylindrical Yamabe invariant
$Y^{h\hbox{-}\cyl}(X)$ is a homothetic invariant,
i.e. $Y^{h\hbox{-}\cyl}(X)=Y^{k\cdot h\hbox{-}\cyl}(X)$ for any
constant $k>0$, but it is not a conformal invariant. This invariant is
a natural generalization of the Yamabe invariant for closed
manifolds. Indeed, let $M$ be a closed manifold of $\dim M=
n\geq 3$ and $p_{\infty}$ a point in $M$. Then the manifold
$X=M\setminus\{p_{\infty}\}$ is an open manifold with the tame end
$Z\times [0,\infty) =S^{n-1}\times [0,\infty)$. Put $h=h_+$. Then we
show that
$$
Y^{\cyl}(M\setminus \{p_{\infty}\})\geq Y^{h_+\hbox{-}\cyl}(M\setminus
\{p_{\infty}\}) = Y(M).
$$
We analyze the dependence of these invariants under gluing in the
following situation.  Let $W_1$, $W_2$ be two compact connected
manifolds with boundaries $\p W_1 =\p W_2 = Z \neq \emptyset$.
Then we have the open manifolds 
$$
X_1 := W_1\cup_Z
Z\times [0,\infty) \ \ \  \mbox{and} \ \ \ 
X_2 := W_2\cup_Z Z\times [0,\infty)
$$ 
with tame ends.  We prove the following Kobayashi-type inequality (see
\cite{K2}).
\vspace{1mm}

\noindent
{\bf Theorem G.} {\it Let $W_1$, $W_2$ be compact manifolds of $\dim
W_i\geq 3$ with $\p W_1 = Z = \p W_2$ and  $X_1$, $X_2$ the
corresponding open manifolds with tame ends. Let $W= W_1\cup_Z W_2$ be
the union of $W_1$ and $W_2$ along the common boundary $Z$. Then, for
any metric $h\in \Riem(Z)$,
$$
Y\!(W_1\!\cup_Z\! W_2\!) \!\geq\!\!
\{\!
\begin{array}{cl}\!\!\!
-\!\!\left(
|Y^{h\!\hbox{-}\!\cyl}\!(X_1\!)|^{\frac{n}{2}}\!\! + \!
|Y^{h\!\hbox{-}\!\cyl}\!(X_2\!)|^{\frac{n}{2}} 
\right)^{\frac{2}{n}} & \!\!\! \!\!\mbox{if} \   Y^{h\!\hbox{-}\!\cyl}(X_i)
\leq 0, \ i\!=\!1,2,
\\
\!\!\!\! \!\!
\min\{Y^{h\hbox{-}\cyl}(X_1), Y^{h\hbox{-}\cyl}(X_2)\} & \!\!\mbox{otherwise}.
\end{array}
\right.\!\!
$$
}

\noindent
Finally, we revisit the surgery construction for the Yamabe invariant
from \cite{PY}. We prove the following result.
\vspace{2mm}

\noindent
{\bf Theorem H.} {\it
Let $M$ be a closed  manifold of $\dim M =n \geq 3$, and $N$ 
an embedded closed submanifold of $M$ with trivial normal bundle. Let
$g_N\in \Riem(N)$ be a given metric on $N$ and $h_+$ the standard
metric on $S^q$. Assume that $q= n-p\geq 3$. 
\begin{enumerate}
\item[{\bf (1)}] If $Y(M)\leq 0$, then for any $\epsilon> 0$ there
exists $\delta= \delta(\epsilon, g_N, |Y(M)|)>0$ such that
$$
Y^{(\kappa^2\cdot h_+ + g_N)\hbox{-}\cyl}(M \setminus N) \geq
Y(M)-\epsilon
$$
for any $0< \kappa\leq \delta$. In particular, $Y^{\cyl}(M \setminus
N)\geq Y(M)$.
\item[{\bf (2)}] If $Y(M)> 0$, then there exists $\delta= \delta(g_N,
Y(M))>0$ such that
$$
Y^{(\kappa^2\cdot h_+ + g_N)\hbox{-}\cyl}(M \setminus N) >0
$$
for any $0< \kappa\leq \delta$. In particular, $Y^{\cyl}(M \setminus N)>0$.
\end{enumerate}
}

\noindent
Theorem H combined with Theorem G gives a refinement of the 
result on surgery and the Yamabe invariant due to Petean and Yun,
see \cite{PY}.
\vspace{2mm}

\noindent
{\bf \ref{int}.10. The plan of the paper.} In Section \ref{s2a}, we
introduce basic terminology and review necessary results.  In Section
\ref{s2b}, we study some properties of the cylindrical Yamabe
constant/invariant under the gluing of manifolds along their
boundary. In particular, we prove Theorem G (Theorem
\ref{cyl:Th7}). In Section \ref{s3}, we revisit the surgery
construction for the Yamabe invariant. In particular, we prove Theorem
H (Theorem \ref{cyl:L29}). Sections \ref{s4} and \ref{s5} form the
core of this paper.  \ In Section \ref{s4}, we first give the setting
for the Yamabe problem.  Then we study the case \ $\lambda({\cal
L}_h)>0$ \ and prove \ Theorem \ref{yam:Th1} \ and \ Theorem
\ref{yam:Th2} \ which imply Theorem A(a) and Theorem B(a). Next we
study the case \ $\lambda({\cal L}_h)=0$ \ and prove \ Theorem
\ref{yam:Th3} \ which implies \ Theorem A(b) \ and Theorem B(b).  In
Section \ref{s5}, we solve the Yamabe problem for canonical
cylindrical manifolds. First we consider the case when
$Y^{\cyl}_{[\bar{h}]}( Z\times \R) < Y(S^n)$ and prove Theorem
\ref{can:Th1} which implies Theorem C. Next we study the supremum case
and prove Theorem \ref{can:Th2} and Corollary \ref{can:Th4} which
imply Theorem D and Theorem F. Theorem E follows from Theorem
\ref{can:Th14} and Theorem \ref{can:Cor17}. In the Appendix, we
discuss the \emph{bottom of the spectrum} of the conformal Laplacian
$\L_{\bar{g}}$ on a cylindrical manifold $(X,\bar{g})$ and its
relationship to the cylindrical Yamabe constant
$Y^{\cyl}_{[\bar{g}]}(X)$ and the sign of scalar curvature.
\vspace{2mm}

\noindent
{\bf \ref{int}.11. Acknowledgments.} The first author would like to
thank the Department of Mathematics at the University of Oregon for
kind hospitality. The second author would like to thank Sergei Novikov
for his support.  Both authors thank Hironori Kumura and Sergei
Preston for useful discussions. Also the authors are grateful to
Richard Schoen for interesting discussions on the Yamabe invariant
during his short visit to the University of Oregon.
\section{Cylindrical manifolds}\label{s2a}
{\bf \ref{s2a}.1. Definition of cylindrical manifolds.}  Here we
review some basic facts and give definitions we need. Let $X$ be an
open smooth manifold of $\dim X = n \geq 3$ without boundary $\p X
=\emptyset$.
\begin{Definition}\label{def-1} {\rm An open smooth
manifold $X$ is called a \emph{manifold with tame ends} if there is a
relatively compact open submanifold $W\subset X$ with $\p \overline{W}
= Z$ such that
\begin{enumerate}
\item[{\bf (1)}] $\displaystyle Z = \bigsqcup_{j=1}^m Z_j$, where each
$Z_j$ is connected,
\item[{\bf (2)}] $\displaystyle X\setminus W \cong Z\times [0,1)= \bigsqcup_{j=1}^m 
\( Z_j \times [0,1)\)$.
\end{enumerate}}
\end{Definition}
We study such manifolds $X$ equipped with cylindrical metrics, which
are defined as follows.

We start by choosing a Riemannian metric $h\in \Riem(Z)$. We identify
the interval $[0,1)$ with the half-line $\R_{\ge 0}=[0,\infty)$, and
hence the product $Z\times [0,1)$ also can be identified with $Z\times
[0,\infty)$.
\begin{figure}[htb]
\hspace*{30mm}\PSbox{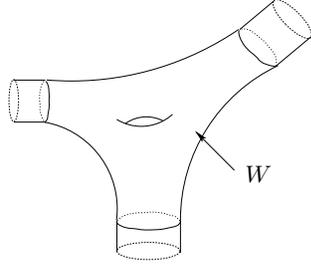}{2cm}{3.2cm}
\begin{picture}(0,0)
\put(30,30){{\small $W$}}
\end{picture}
\caption{A manifold with tame ends}\label{fig2-1}
\end{figure}
\begin{Definition}\label{def-2}
{\rm A complete Riemannian metric $\bar{g}\in \Riem(X)$ is called a
{\sl cylindrical metric modeled by} $(Z,h)$ if there exists a coordinate
system $(x,t)$ on $Z\times [0,\infty)$ such that 
$$
\bar{g}(x,t) = h(x) + dt^2 \ \ \ \mbox{on} \ \ \ Z\times [1,\infty)\subset 
Z\times [0,\infty)\cong X\setminus W.
$$
Let $\Riem^{\cyl}(X)$ ($\subset \Riem(X)$) be the space of cylindrical
metrics modeled by $(Z,h)$ for some Riemannian metric $h\in
\Riem(Z)$.}
\end{Definition}
We also define the map $\p_{\infty}: \Riem^{\cyl}(X)\to \Riem(Z)$
by $\p_{\infty}\bar{g}= h$ if $\bar{g}$ is a cylindrical metric
as above. Clearly the map $\p_{\infty}$ is onto since any metric $h\in
\Riem(Z)$ can be extended to a cylindrical metric on $X$.
\begin{figure}[hb]
\hspace*{10mm}\PSbox{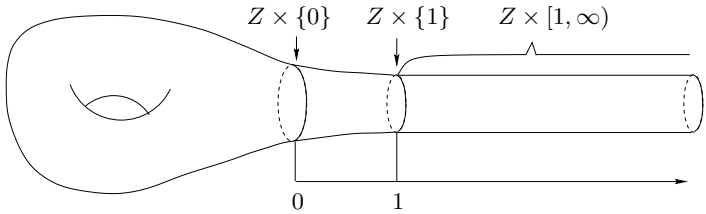}{4cm}{25mm}
\begin{picture}(0,0)
\put(-25,65){{\small $Z\times\{0\}$}}
\put(20,65){{\small $Z\times\{1\}$}}
\put(70,65){{\small $Z\times [1,\infty)$}}
\put(-8,-5){{\small $0$}}
\put(30,-5){{\small $1$}}
\end{picture}
\caption{Cylindrical manifold $X$.}\label{fig2-2}
\end{figure}

\noindent
{\bf Remarks. (1)} For simplicity, we will often assume that the slice
manifold $Z$ is connected. In the general case, i.e.  when
$Z=\bigsqcup_{j=1}^m Z_j$, one can easily obtain the corresponding
results to those we prove here by using that $\lambda({\cal
L}_h)=\displaystyle \min_{1\leq j\leq m} \lambda({\cal L}_h|_{Z_j})$
and $Y^{\cyl}_{[h+dt^2]}(Z\times \R) = \displaystyle \min_{1\leq j\leq
m} Y^{\cyl}_{[h|_{Z_j} +dt^2]}(Z_j\times \R)$.
\vspace{1mm}

\noindent
{\bf (2)} Later on we will also consider the case when a cylindrical
manifold $X$ has a non-empty boundary. In that case the definition is
modified by assuming that $\p X= \p W$ and $\p \overline{W}= \p
X \sqcup Z$ (see Fig. \ref{fig4-1}). 
\vspace{1mm}

\noindent
{\bf (3)} From the viewpoint of conformal geometry, cylindrical
manifolds were considered implicitly in \cite{AB2}.
\vspace{2mm}

\noindent
{\bf \ref{s2a}.2. Cylindrical conformal classes.} Next, we define an
appropriate notion of a conformal class on a cylindrical manifold
$X$. Let $\bar{C}=[\bar{g}]$ denote the regular conformal class for a
metric $\bar{g}\in \Riem(X)$. Let ${\cal C}(X)$ be the space of
conformal classes on $X$.  We define the space of cylindrical
conformal classes
$$
{\cal C}^{\cyl}(X):= \{ [\bar{g}] \ | \ \bar{g}\in
\Riem^{\cyl}(X) \}\subset {\cal C}(X).
$$
For $\bar{C}\in {\cal C}^{\cyl}(X)$, we fix a cylindrical metric
$\bar{g}\in \bar{C}\cap \Riem^{\cyl}(X)$ to define the
$L_{\bar{g}}^{k,2}$-conformal class
$$
\bar{C}_{L_{\bar{g}}^{k,2}}:= \{ u^{\frac{4}{n-2}}\cdot \bar{g} \ |
\ u \in C^{\infty}_+(X)\cap L_{\bar{g}}^{k,2}(X) \} \subset \bar{C},
$$
where $k=0,1,2$. 

\begin{figure}[hpt]
\hspace*{10mm}\PSbox{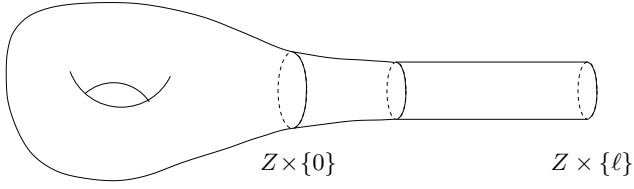}{4cm}{2.5cm}
\begin{picture}(0,0)
\put(-20,5){{\small $Z\!\times\!\{0\}$}}
\put(90,5){{\small $Z\times\{\ell\}$}}
\end{picture}
\caption{Manifold $X(\ell)$.}\label{fig2-3}
\end{figure}

Now we show that the conformal classes $\bar{C}_{L_{\bar{g}}^{k,2}}$
do not depend on the choice of reference cylindrical metrics
$\bar{g}$ for $k=0,1$.
\begin{Proposition}\label{conf:new}
Let $\bar{C}$ be a cylindrical conformal class and
$\bar{g},\tilde{g}\in
\bar{C}\cap \Riem^{\cyl}(X)$ two cylindrical metrics. Then there
exists a constant $K\geq 1$ such that $K^{-1}\cdot \bar{g}\leq
\tilde{g}\leq K\cdot \bar{g}$ on $X$. In particular, 
$$
\begin{array}{l}
\displaystyle
K^{-\frac{n}{2}} \cdot d\sigma_{\bar{g}} \leq d\sigma_{\tilde{g}}\leq
K^{\frac{n}{2}} \cdot d\sigma_{\bar{g}} \ \ \mbox{on} \ \ X, \ \ \mbox{and}
\\
\displaystyle
L^{k,2}_{\bar{g}}(X) = L^{k,2}_{\tilde{g}}(X), \ \ \ \bar{C}_{L_{\bar{g}}^{k,2}} =
\bar{C}_{L_{\tilde{g}}^{k,2}} \ \ \mbox{for} \ \ k=0,1.
\end{array}
$$
\end{Proposition}
\begin{Proof}
We have that $\bar{g}(x,t)= h(x)+dt^2$ on $Z\times [1,\infty)$ and
$\tilde{g}(y,s)= \tilde{h}(y)+ds^2$ on $\tilde{Z}\times [1,\infty)$,
where $y=y(x,t)$ and $s=s(x,t)$ give a diffeomorphism: $Z\times
[1,\infty)\cong \tilde{Z}\times [1,\infty)$.  By the assumption,
$\tilde{g}=\phi\cdot\bar{g}$ for some function $\phi\in
C_+^{\infty}(X)$. We write $\phi=\phi(x,t)$ on $Z\times [1,\infty)$.
It is enough to show
\begin{equation}\label{conf:eq1}
\inf_{X} \phi > 0, \ \ \ \sup_{X} \phi <\infty.
\end{equation}
Indeed, suppose that \ $\inf_{X} \phi = 0$. \ Then there exists a sequence
$\{(x_i, t_i)\}\subset Z\times [1,\infty)$ such that $\phi(x_i,t_i)\to
0$ as $i\to\infty$. In particular, it implies that the injectivity
radius $\inj_X \tilde{g}=0 $.  On the other hand, since $\tilde{g}$ is
a cylindrical metric, the injectivity radius $\inj_X \tilde{g} \geq
\delta$ for some $\delta>0$. This leads to a contradiction. A similar
argument also shows that $\sup_{X} \phi <\infty$. Hence the property
(\ref{conf:eq1}) holds.
\end{Proof}
Now we define the following functional on the space
$L_{\bar{g}}^{1,2}(X)$. Set
$$
Q_{(X,\bar{g})}(u) := \frac{\int_X\[\alpha_n |du|^2 +
R_{\bar{g}}u^2\] d\sigma_{\bar{g}}}{\(\int_X
|u|^{\frac{2n}{n-2}}d\sigma_{\bar{g}}\)^{\frac{n-2}{n}}}, \ \ \mbox{where} \ \ 
\alpha_n= \frac{4(n-1)}{n-2},
$$
for $u\in L_{\bar{g}}^{1,2}(X)$ with $u \not\equiv 0$.
\begin{Lemma}\label{cyl:L1}
The functional $Q_{(X,\bar{g})}(u)$ is well-defined on the space
$L_{\bar{g}}^{1,2}(X)$ with $u \not\equiv 0$.
\end{Lemma}
\begin{Proof}
Let $(N,g)$ be a complete Riemannian manifold. We define
$$
C_c^{\infty}(N) : = \{ f\in C^{\infty}(N) \ | \ \supp(f) \ 
\mbox{is compact} \}.
$$
The following facts are well known (cf. \cite{Au3}).
\begin{Fact}\label{fact} $\mbox{ \ }$
\begin{enumerate}
\item[{\bf (i)}] {\sl The $L_{g}^{k,2}$-completion of the space $C_c^{\infty}(N)$
coincides with the space $L_{g}^{k,2}(N)$ for $k=0,1$.}
\end{enumerate}
Now we assume that the sectional curvature $K_g$ and the injectivity
radius $\iota_g$ are bounded, i.e. there exist constants $C>0$ and
$\delta>0$ such that $|K_g|\leq C$ and $\iota_g\geq \delta$. Under this
assumption, we have that 
\begin{enumerate}
\item[{\bf (ii)}] {\sl The $L_{g}^{2,2}$-completion of the space $C_c^{\infty}(N)$
coincides with the space $L_{g}^{2,2}(N)$.}
\item[{\bf (iii)}] {\sl The Sobolev embedding $L_{g}^{1,2}(N)\subset
L_{g}^{\frac{2n}{n-2}}(N)$ holds.}
\end{enumerate}
\end{Fact}

Now we return to a cylindrical manifold $X$. The Sobolev embedding
(iii) implies that
$$
\int_X |u|^{\frac{2n}{n-2}} d\sigma_{\bar{g}} < \infty \ \ \ \mbox{for} \ \ 
u\in L_{\bar{g}}^{1,2}(X).
$$
We define $X(\ell):= W \cup_Z \( Z\times [0,\ell]\)$ for $\ell\geq 1$,
see Fig. \ref{fig2-3}.

For $u\in L_{\bar{g}}^{1,2}(X)$, we have the following estimates:
$$
\begin{array}{rcl}
\displaystyle 
\left| \int_X R_{\bar{g}}u^2 d\sigma_{\bar{g}}\right|  &\leq &
\displaystyle 
\left| \int_{X(1)} R_{\bar{g}}u^2 d\sigma_{\bar{g}} \right| +
\int_{X\setminus X(1)} | R_{\bar{g}}| u^2 d\sigma_{\bar{g}}
\\
\\
&\leq &\displaystyle 
\left| \int_{X(1)} R_{\bar{g}}u^2 d\sigma_{\bar{g}} \right| +
(\max_{Z} R_h) \cdot \int_{X\setminus X(1)} u^2 d\sigma_{\bar{g}} <\infty.
\end{array}
$$
This completes the proof of Lemma \ref{cyl:L1}.
\end{Proof}
{\bf \ref{s2a}.3. Cylindrical Yamabe constant.}  For $\bar{g}\in
\Riem^{\cyl}(X)$, we define the functional
$$
I_X(\tilde{g}):= \frac{\int_X R_{\tilde{g}}d\sigma_{\tilde{g}}}
{\Vol_{\tilde{g}}(X)^{\frac{n-2}{n}}}
$$
on the space of $L^{2,2}_{\bar{g}}$-conformal metrics,
i.e. $\tilde{g}\in [\bar{g}]_{L^{2,2}_{\bar{g}}}$. Recall that
$\tilde{g} = u^{\frac{4}{n-2}}\cdot\bar{g}$, where $u\in
C_+^{\infty}(X)\cap L^{2,2}_{\bar{g}}(X)$.
\begin{Lemma}\label{cyl:L2}
The functional $I_X(\tilde{g})$ is well-defined for metrics
$\tilde{g}\in [\bar{g}]_{L^{2,2}_{\bar{g}}}$. Furthermore,
$I_X(\tilde{g}) = Q_{(X,\bar{g})}(u)$ if $\tilde{g} =
u^{\frac{4}{n-2}}\cdot\bar{g}$.
\end{Lemma}
\begin{Proof}
For $\tilde{g} = u^{\frac{4}{n-2}}\cdot\bar{g}$ with $u\in
C_+^{\infty}(X)\cap L^{2,2}_{\bar{g}}(X)$, we first notice
$$
\Vol_{\tilde{g}}(X) =
\int_Xd\sigma_{\tilde{g}} = \int_X u^{\frac{2n}{n-2}}
d\sigma_{\bar{g}} < \infty .
$$
Then we have:
$$
\begin{array}{rcl}\displaystyle 
\int_X R_{\tilde{g}}d\sigma_{\tilde{g}}  &= & \displaystyle 
\int_X\! u^{-\frac{n+2}{n-2}}\left(-\alpha_n\Delta_{\bar{g}} u + 
R_{\bar{g}} u \right)
\cdot u^{\frac{2n}{n-2}}\! d\sigma_{\bar{g}}
\\
\\
&= & \displaystyle 
\int_X \!\!\left(-\alpha_n u\cdot \Delta_{\bar{g}} u + R_{\bar{g}} u^2 \right)
\!d\sigma_{\bar{g}}.
\end{array}
$$
On the other hand, we have 
$$
\begin{array}{rcl}\displaystyle 
\int_X R_{\tilde{g}}d\sigma_{\tilde{g}} &=&  \displaystyle 
\lim_{\ell\to\infty}\int_{X(\ell)} R_{\tilde{g}}d\sigma_{\tilde{g}} 
\\
\\
&=&  \displaystyle 
\lim_{\ell\to\infty}\left[\int_{X(\ell)} \!\!
\left( \alpha_n |d u|^2\!  +\!  R_{\bar{g}} u^2 \right)
d\sigma_{\bar{g}}\!  -\!  2\alpha_n \int_{Z\times \{\ell\}} \!\!\! 
u\cdot \p_t u \ d\sigma_h
\right].
\end{array}
$$
Hence to prove that $I_X(\tilde{g}) = Q_{(X,\bar{g})}(u)$, it is enough
to prove the following.
\begin{Claim}\label{cyl:L3}
$ \displaystyle \int_{Z\times \{\ell\}} \!\!\!  u\cdot \p_t u \
d\sigma_h =o(1)$ as $\ell\to\infty$.
\end{Claim}
\begin{Proof}
Let  $f(t)$ be the function defined as
$$
f(t):= \int_{Z\times \{t\}} \!\!\!  u^2 \ d\sigma_h, \ \ \mbox{then} \ \ 
f^{\prime}(t) = 2 \int_{Z\times \{\ell\}} \!\!\!  u\cdot \p_t u \
d\sigma_h.
$$
Since $u\in L^{2,2}_{\bar{g}}(X)$, there exist $\epsilon>0$ and
$\ell_{\epsilon}>0$ such that
$$
\int_{Z\times [\ell_{\epsilon},\infty)} \!\!\! 
(u^2 + | \p_t u|^2 + | \p_t^2 u|^2 )d\sigma_h \leq \epsilon.
$$
Note that the volume form $d\sigma_{\bar{g}}= d\sigma_h \ dt$ on
the cylindrical part $Z\times [1,\infty)$. Now for any $t_2>t_1\geq
\ell_{\epsilon}$ we have
\begin{equation}\label{cyl:eq-2}
\begin{array}{rcl}
\displaystyle 
|f(t_2)-f(t_1)| &=& \left| \int_{t_1}^{t_2} f^{\prime}(t)\right|  =  
\displaystyle 
2  \left| 
\int_{Z\times [t_1,t_2]} \!\!\! 
u\cdot \p_t u \ d\sigma_h \ dt
\right|
\\
\\
&\leq &  \displaystyle 
2 \left( \int_{Z\times [t_1,t_2]} \!\!\! u^2  d\sigma_{\bar{g}}\right)^{\frac{1}{2}} \cdot
\left( \int_{Z\times [t_1,t_2]} \!\!\! |\p_t u|^2  d\sigma_{\bar{g}}\right)^{\frac{1}{2}} \leq
2\epsilon .
\end{array}\!\!\!\!\!\!\!\!\!\!\!\!\!\!
\end{equation}
We claim that $f(\ell_{\epsilon}) \leq 3\epsilon$ and $f(t) \leq
5\epsilon$ for all $t\geq \ell_{\epsilon}$. 

Indeed, suppose that $f(\ell_{\epsilon}) > 3\epsilon$. Then
(\ref{cyl:eq-2}) implies that $f(t) > f(\ell_{\epsilon}) - 2\epsilon
\geq \epsilon $ for $t\geq \ell_{\epsilon}$. Therefore one has
$$
\int_{Z\times [\ell_{\epsilon},\infty)} \!\!\! u^2  d\sigma_{\bar{g}} = 
\int_{\ell_{\epsilon}}^{\infty} f(t) \ dt \geq 
\int_{\ell_{\epsilon}}^{\infty} \epsilon \ dt  = \infty,
$$
which contradicts the condition $u\in L^{2,2}_{\bar{g}}(X)$.  Hence
$f(\ell_{\epsilon}) \leq 3\epsilon$, and $ 0 < f(t) \leq
f(\ell_{\epsilon}) +2\epsilon \leq 5 \epsilon $ for $t\geq
\ell_{\epsilon}$. Now we define
$$
g(t):= \int_{Z\times\{t\}} |\p_t u|^2 \ d\sigma_h \geq 0 \ \ \mbox{with} \ \ 
g^{\prime}(t) = 2 \int_{Z\times\{t\}} \p_t u \cdot \p_t^2 u \ d\sigma_h .
$$
Similarly, one can show that $g(t)\leq 5\epsilon$ for all $t\geq
\ell_{\epsilon}$.  Now we have
$$
\begin{array}{rcl}
\displaystyle 
\int_{Z\times \{t\}} \!\!\!  | u | \cdot |\p_t u | \
d\sigma_h  &\leq&  \displaystyle 
\left(\int_{Z\times \{t\}} \!\!\!  u^2  \
d\sigma_h \right)^{\frac{1}{2}} \cdot
\left(\int_{Z\times \{t\}} \!\!\!  |\p_t u|^2  \
d\sigma_h \right)^{\frac{1}{2}} 
\leq 
\displaystyle 5 \epsilon 
\end{array}
$$
for all  $t\geq \ell_{\epsilon}$. Clearly this implies that
$$
\int_{Z\times \{t\}} \!\!\!  | u | \cdot |\p_t u | \
d\sigma_h = o(1) 
$$
as $t\to\infty$. 
\end{Proof}
This completes the proof of Lemma \ref{cyl:L2}.
\end{Proof}
\begin{Definition}\label{cyl:D1}
For $\bar{g}\in \Riem^{\cyl}(X)$, we define the constant
$$
Y_{\bar{g}}^{\cyl}(X):=\inf_{\tilde{g}\in
[\bar{g}]_{L^{2,2}_{\bar{g}}}} I_X(\tilde{g}).
$$
\end{Definition}
\begin{Lemma}\label{cyl:L4}
The following identities hold:
$$
\begin{array}{rcl}\displaystyle
Y_{\bar{g}}^{\cyl}(X) & = & \displaystyle
\inf_{\tilde{g}\in
[\bar{g}]_{L^{2,2}_{\bar{g}}}} I_X(\tilde{g})
= \inf_{u\in C_+^{\infty}(X)\cap
L^{1,2}_{\bar{g}}(X)} Q_{(X,\bar{g})}(u)
\\
\\
& = & \displaystyle
\inf_{u\in  L^{1,2}_{\bar{g}}(X), u\not\equiv 0} Q_{(X,\bar{g})}(u)
= \inf_{u\in  C_c^{\infty}(X), u\not\equiv 0} Q_{(X,\bar{g})}(u).
\end{array}
$$
\end{Lemma}
\begin{Proof}
The assertion follows from the following inequalities:
$$
\begin{array}{rclcl}\displaystyle
Y_{\bar{g}}^{\cyl}(X) &\geq &
\displaystyle
\inf_{u\in C_+^{\infty}(X)\cap
L^{1,2}_{\bar{g}}(X)} Q_{(X,\bar{g})}(u)
&\geq &\displaystyle
\inf_{u\in 
L^{1,2}_{\bar{g}}(X), u\not\equiv 0} Q_{(X,\bar{g})}(u)
\\
\\
&= &
\displaystyle
\inf_{u\in  C_c^{\infty}(X), u\not\equiv 0} Q_{(X,\bar{g})}(u) 
&= &\displaystyle
\inf_{\begin{array}{c}^{u\in  C_c^{\infty}(X),}\\
^{u\geq 0, u\not\equiv 0}\end{array}} Q_{(X,\bar{g})}(u) 
\\
\\
&\geq &
\displaystyle
\inf_{u\in C_+^{\infty}(X)\cap
L^{2,2}_{\bar{g}}(X)} Q_{(X,\bar{g})}(u) &=&\displaystyle
 Y_{\bar{g}}^{\cyl}(X) .
\end{array}
$$
Here we make use of the fact that the $L_{\bar{g}}^{k,2}$-completion of the
space $C_c^{\infty}(X)$ coincides with the space $L_{\bar{g}}^{k,2}(X)$.
We use essentially that $\bar{g}=h+dt^2$ on $Z\times [1,\infty)$ to obtain the
last inequality.
\end{Proof}
Now we would like to recall the following important observation due to
Schoen and Yau \cite{SY3}.  We state the corresponding result in the
above terms.
\begin{Fact}\label{cyl:Th1}
{\rm \cite[Section 2]{SY3}}
For any $g^{\prime}, \ g^{\prime\prime}\in \bar{C}$, 
$$
\inf_{u\in  C_c^{\infty}(X), u\not\equiv 0} Q_{(X,g^{\prime})}(u) =
\inf_{u\in  C_c^{\infty}(X), u\not\equiv 0} Q_{(X,g^{\prime\prime})}(u). 
$$
In particular, these constants are conformally invariant.
\end{Fact}
We use Fact \ref{cyl:Th1} or Proposition \ref{conf:new} to conclude
the following.
\begin{Corollary}\label{cyl:Th2}
Let $\bar{g}, \ \hat{g}\in \bar{C}\cap \Riem^{\cyl}(X)$ be any two
cylindrical metrics. Then $Y_{\bar{g}}^{\cyl}(X)=Y_{\hat{g}}^{\cyl}(X)$.
\end{Corollary}
Now we define the \emph{cylindrical Yamabe constant} as follows.
\begin{Definition}\label{cyl:D2}
{\rm Let $\bar{C}\in {\cal C}^{\cyl}(X)$ be a cylindrical conformal
class. Then the \emph{cylindrical Yamabe constant
$Y_{\bar{C}}^{\cyl}(X)$ of $(X,\bar{C})$} is defined by
$Y_{\bar{C}}^{\cyl}(X):= Y_{\bar{g}}^{\cyl}(X)$ for any cylindrical
metric $\bar{g}\in \bar{C}\cap \Riem^{\cyl}(X)$.}
\end{Definition}
Corollary \ref{cyl:Th2} implies that the Yamabe constant
$Y_{\bar{C}}(X)$ is well-defined, and Lemma \ref{cyl:L4} gives the formula:
$$
Y_{\bar{C}}^{\cyl}(X) =\inf_{\tilde{g}\in \bar{C}_{L^{2,2}_{\bar{g}}}}
I_X(\tilde{g}) =\inf_{u\in 
L^{1,2}_{\bar{g}}(X), u\not\equiv 0} Q_{(X,\bar{g})}(u) .
$$
{\bf \ref{s2a}.4. Finiteness of cylindrical Yamabe constants.}  Let
$(X,\bar{g})$ be a cylindrical manifold modeled by $(Z,h)$ as above.
We would like to give a complete criterion for the finiteness of
$Y_{[\bar{g}]}^{\cyl}(X)$ in terms of the geometry of the slice
Riemannian manifold $(Z,h)$.  Recall that $\dim X=n$, and hence $\dim
Z= n-1$.  We define the operator
$$
{\cal L}_h := -\frac{4(n-1)}{n-2}\Delta_h + R_h = -\alpha_n \Delta_h + R_h
\ \ \ \mbox{on} \ \ \ (Z,h). 
$$
Recall that the operator $\L_h = -\frac{4(n-2)}{n-3}\Delta_h +
R_h=-\alpha_{n-1} \Delta_h + R_h $ (when $n\geq 4$) is the conformal
Laplacian of $(Z,h)$. Clearly ${\cal L}_h$ resembles the conformal
Laplacian $\L_h$ of $(Z,h)$, but it is not equal to $\L_h$.  The first
eigenvalue $\lambda({\cal L}_h)$ of ${\cal L}_h$ (given by the formula
below) captures a complete information on the finiteness of
$Y_{[\bar{g}]}(X)$. We have:
$$
\lambda({\cal L}_h):= \inf_{u\in  L^{1,2}(Z), u\not\equiv 0} 
\frac{ \int_Z \left( \alpha_n |du|^2 +R_h u^2
\right) d\sigma_h}{\int_Z u^2 \ d\sigma_h}.
$$
\begin{Lemma}\label{cyl:L5}
Assume that $\lambda({\cal L}_h)<0$. Then $Y_{\bar{C}}^{\cyl}(X)=
-\infty$ for any cylindrical conformal class $\bar{C}=[\bar{g}]\in
{\cal C}^{\cyl}(X)$ with $\bar{g}\in \Riem^{\cyl}(X)$ and $\p_{\infty}
\bar{g} =h$.  In particular, $Y_{\bar{C}}^{\cyl}(X)= -\infty$ if the
scalar curvature $R_h<0$ on $Z$.
\end{Lemma}
\begin{Proof}
Without loss of generality, we may assume that $Z$ is connected. We
define a Lipschitz function $f_{\epsilon}(t)\in C_c^{0,1}([0,\infty))$
whose graph is given at Fig. \ref{fig2-4}.
\begin{figure}[h]
\hspace*{15mm}\PSbox{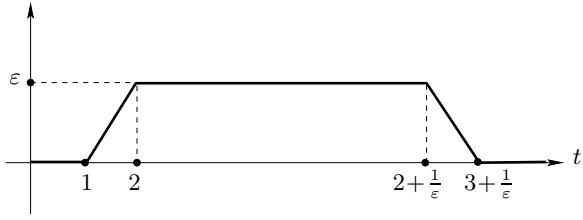}{4cm}{3cm}
\begin{picture}(0,0)
\put(-115,50){{\small $\epsilon$}}
\put(-88,10){{\small $1$}}
\put(-70,10){{\small $2$}}
\put(30,10){{\small $2\!+\!\frac{1}{\epsilon}$}}
\put(57,10){{\small $3\!+\!\frac{1}{\epsilon}$}}
\put(98,20){{\small $t$}}
\end{picture}
\caption{Function $f_{\epsilon}(t)$.}\label{fig2-4}
\end{figure}
The standard elliptic theory implies that there exists a function
$\phi\in C_+^{\infty}(Z)$ such that
$$
\begin{array}{l}
\displaystyle
\int_Z \left( \alpha_n  |d\phi|^2 +R_h \phi^2
\right) d\sigma_h = \lambda({\cal L}_h) \ \ \ \ \mbox{and} \ \ \ 
\displaystyle
\int_Z \phi^2 d\sigma_h = 1.
\end{array}
$$
We define a family of functions depending on $\epsilon$ by
$u_{\epsilon}(x,t):=f_{\epsilon}(t)\cdot \phi(x) \in C_c^{0,1}(X)$. We
estimate $Q_{(X,\bar{g})}(u_{\epsilon})$ from above:
$$
\begin{array}{rcl}
Q_{(X,\bar{g})}(u_{\epsilon}) &=& \displaystyle
\frac{\lambda({\cal L}_h)\cdot\int_{0}^{\infty} f_{\epsilon}^2\ dt + \alpha_n \int_{0}^{\infty}
(f_{\epsilon}^{\prime})^2\ dt}{\left(\int_Z \phi^{\frac{2n}{n-2}}_{_{\mbox{\begin{picture}(0,5)
\end{picture}}}} 
d\sigma_h\right)^{\frac{n-2}{n}}\cdot 
\left( \int_{0}^{\infty} f_{\epsilon}^{\frac{2n}{n-2}}dt\right)^{\frac{n-2}{n}}}
\\
\\
&\leq & \displaystyle
\frac{\lambda({\cal L}_h)\cdot \epsilon^2\cdot\frac{1}{\epsilon} + 
2 \alpha_n\epsilon^2}{C^{-1}
\left(\epsilon^{\frac{2n}{n-2}}\cdot\frac{2}{\epsilon}\right)^{\frac{n-2}{n}}}
\\
\\
&\leq & \displaystyle
\frac{C^{\prime}\left(\lambda({\cal L}_h) + 
2 \alpha_n \epsilon\right)\epsilon}{\epsilon^{\frac{n+2}{n}}} =
\frac{C^{\prime}\left(\lambda({\cal L}_h) + 
2 \alpha_n \epsilon\right)}{\epsilon^{\frac{2}{n}}} .
\end{array}
$$
Here \ $C$ \ and \ $C^{\prime}$ \ are some positive constants. \ Finally 
we have that  $Q_{(X,\bar{g})}(u_{\epsilon})
\to -\infty$ as $\epsilon\to 0$ since $\lambda({\cal L}_h)<0$.
Hence $Y_{\bar{C}}^{\cyl}(X)=-\infty$.
\end{Proof}
To proceed further, we recall the following result, which is essentially
due to Aubin \cite{Au2}. 
\begin{Fact}\label{cyl:Th3}{\rm (Aubin's inequality)} The cylindrical Yamabe
constant $Y_{\bar{C}}^{\cyl}(X)$ is bounded from above:
$Y_{\bar{C}}^{\cyl}(X)\leq Y(S^n)$, where $Y(S^n)$ is the Yamabe
invariant of the $n$-sphere.
\end{Fact}
We conclude that for any cylindrical conformal class $\bar{C}\in {\cal
C}^{\cyl}(X)$
$$
-\infty \leq Y_{\bar{C}}^{\cyl}(X)\leq Y(S^n).
$$
\begin{Definition}\label{cyl:D3}
{\rm Let $X$ be an open manifold of $\dim X =n\geq 3$ with tame
ends. Let $h$ be a metric on $Z$. We define the \emph{$h$-cylindrical
Yamabe invariant of $X$} by
$$
Y^{h\hbox{-}\cyl}(X):=\sup_{^{^{\begin{array}{c}
_{\bar{g}\in\Riem^{\cyl}(X)}\\
_{\p_{\infty}\bar{g}=h}
\end{array}}}} Y_{[\bar{g}]}^{\cyl}(X).
$$
The \emph{cylindrical Yamabe invariant of $X$} is also defined by}
$$
Y^{\cyl}(X):=\sup_{h\in \Riem(Z)} Y^{h\hbox{-}\cyl}(X).
$$
\end{Definition}
Then the following inequalities hold:
$$
-\infty\leq Y^{h\hbox{-}\cyl}(X) \leq Y^{\cyl}(X) \leq Y(S^n).
$$
\begin{Lemma}\label{cyl:L10} The $h$-cylindrical Yamabe invariant
$Y^{h\hbox{-}\cyl}(X)$ is a homothetic invariant:
$$
Y^{h\hbox{-}\cyl}(X) = Y^{kh\hbox{-}\cyl}(X),
$$
where $kh(x):=k\cdot h(x)$ for any constant $k>0$.
\end{Lemma}
\begin{Proof} Choose any $\bar{g}\in \Riem^{\cyl}(X)$ with 
$\p_{\infty} \bar{g} = kh$. Then $\bar{g}=kh+dt^2$ on $Z\times
[1,\infty)$. Set $s=\frac{1}{\sqrt{k}}t$, and then $ds=
\frac{1}{\sqrt{k}}dt$. We have
$$
\bar{g}=kh+kds^2 = k(h+ ds^2) \ \ \ \mbox{on} \ \ Z\times[1/\sqrt{k},\infty).
$$
We obtain: 
$$ 
Y^{\cyl}_{[\bar{g}]}(X) =
Y^{\cyl}_{[\frac{1}{k}\bar{g}]}(X) \leq Y^{h\hbox{-}\cyl}(X).  
$$ 
Hence $Y^{h\hbox{-}\cyl}(X) \geq Y^{kh\hbox{-}\cyl}(X)$. Similarly, $
Y^{h\hbox{-}\cyl}(X)\leq Y^{kh\hbox{-}\cyl}(X)$.
\end{Proof}
\vspace{1mm}

\noindent
{\bf Remark.} We emphasize that, in general, the invariant
$Y^{h\hbox{-}\cyl}(X)$ is \emph{not a conformal invariant} with
respect to a metric in the conformal class $[h]$. 
\vspace{3mm}

\noindent
{\bf Example (1)} Set $Z=T^{n-1}$ for $n\geq 4$. Let $h_0$ be a flat
metric on $T^{n-1}$. Then
$$
{\cal L}_{h_0} = -\alpha_n \Delta_{h_0}, \ \ \ \mbox{and} \ \ \ \ 
\lambda({\cal L}_{h_0})=0.
$$
Hence $0\geq Y^{h_0\hbox{-}\cyl}(X) >-\infty$. On the other hand, let $h\in
[h_0]$ be a non-flat metric on $T^{n-1}$. Clearly there exists a
function $u\in C_+^{\infty}(T^{n-1})$ such that
$$
-\alpha_{n-1}\Delta_{h} u + R_h u = 0.
$$ 
Then we have:
$$
\begin{array}{rcl}
\displaystyle {\cal L}_{h}u & = &\displaystyle -\alpha_n \Delta_{h}u +
R_h u = \displaystyle -\alpha_{n-1}\Delta_{h}u + R_h u +
\frac{4}{(n-2)(n-3)}\Delta_{h} u .
\end{array}
$$
In particular, 
$$
({\cal L}_{h} u, u) = - \frac{4}{(n-2)(n-3)}\int_{Z} |d u|^2 d\sigma_h < 0.
$$
This implies that $\lambda({\cal L}_{h})<0$, and consequently,
$Y^{h\hbox{-}\cyl}(X) =-\infty$. Hence we obtain that
$Y^{h\hbox{-}\cyl}(X)\neq Y^{h_0\hbox{-}\cyl}(X)$.
\vspace{2mm}

\noindent
{\bf Remark.} For any canonical cylindrical manifold, we notice that
$$
Y^{h\hbox{-}\cyl}(Z\times \R) = Y^{\cyl}_{[h+dt^2]}(Z\times \R) ,
$$
see Proposition \ref{cyl:univ-bound}.
\vspace{2mm}

\noindent
{\bf Example (2)} Set $Z= S^{n-1}$ and $X=S^{n-1}\times \R$ for $n\geq
3$.  Then
$$
Y^{h_+\hbox{-}\cyl}(X)=Y^{\cyl}_{[h+dt^2]}(X) = Y(S^n).
$$
On the other hand, let $h\in [h_+]$ be any metric whose sectional
curvature is not identically equal to a positive constant. Then we
will prove in Theorem \ref{can:Th2} in Section \ref{s5} that
$$
Y^{h\hbox{-}\cyl}(X)= Y^{\cyl}_{[h_++dt^2]}(X) < Y(S^n).
$$
Hence $Y^{h\hbox{-}\cyl}(X)\neq Y^{h_+\hbox{-}\cyl}(X)$. 
\begin{Lemma}\label{cyl:L6} Let $(X,\bar{g})$ be a cylindrical manifold with $\p_{\infty}\bar{g}=h$, and set $\bar{C}=[\bar{g}]$.
\begin{enumerate}
\item[{\bf (1)}] If $\lambda({\cal L}_h)\geq 0$, then
$Y_{\bar{C}}^{\cyl}(X)>-\infty$, and hence 
$Y^{h\hbox{-}\cyl}(X) >-\infty$. In particular, if $R_h\geq 0$ on $Z$, then
$Y_{\bar{C}}^{\cyl}(X)>-\infty$ and $Y^{h\hbox{-}\cyl}(X) >-\infty$.
\item[{\bf (2)}] If $\lambda({\cal L}_h)= 0$, then
$Y_{\bar{C}}^{\cyl}(X)\leq 0$, and hence $Y^{h\hbox{-}\cyl}(X) \leq 0
$. In particular, if $R_h\equiv 0$ on $Z$, then
$Y_{\bar{C}}^{\cyl}(X)\leq 0$ and $Y^{h\hbox{-}\cyl}(X) \leq 0$.
\end{enumerate}
\end{Lemma}
\begin{Proof}
{\bf (1)} Let $u\in L_{\bar{g}}^{1,2}(X)$ be a function with 
$
\int_X |u|^{\frac{2n}{n-2}} d\sigma_{\bar{g}} =1.
$
Then we have:
$$
\begin{array}{rcl}
Q_{(X,\bar{g})}(u) &\geq & \displaystyle
\int_{X(1)} \left(\alpha_n |du|^2 + R_{\bar{g}} u^2 \right)d\sigma_{\bar{g}}
\\
\\
&&
\displaystyle \  
+ \ \lambda({\cal L}_h)\cdot\int_{Z\times[1,\infty)} \!\!\!\! u^2 \ d\sigma_h\ dt +
\alpha_n \cdot \int_{Z\times[1,\infty)} \!\!\!\! |\p_tu|^2 \ d\sigma_h\ dt \ .
\end{array}
$$
Since $\lambda({\cal L}_h)\geq 0$ we have:
$$
\begin{array}{rcl}
Q_{(X,\bar{g})}(u) &\geq & \displaystyle
\int_{X(1)} \left(\alpha_n |du|^2 + R_{\bar{g}} u^2 \right)d\sigma_{\bar{g}}
\geq 
-\left(\max_{X(1)} |R_{\bar{g}}|
\right)\int_{X(1)} u^2 \ d\sigma_{\bar{g}}
\\
\\
 &\geq & \displaystyle
-\left(\max_{X(1)} |R_{\bar{g}}|
\right) \left(\int_{X(1)} |u|^{\frac{2n}{n-2}} d\sigma_{\bar{g}}\right)^{\frac{n-2}{n}}
\cdot\Vol_{\bar{g}}(X(1))^{\frac{2}{n}} 
\\
\\
&\geq & \displaystyle
-\left(\max_{X(1)} |R_{\bar{g}}|
\right) 
\cdot\Vol_{\bar{g}}(X(1))^{\frac{2}{n}} .
\end{array}
$$
This gives that $Y^{h\hbox{-}\cyl}(X) >-\infty$.
\vspace{2mm}

\noindent
{\bf (2)} Let $u_{\epsilon}$ be the same function as in the proof of
Lemma \ref{cyl:L5}. Then
$$
Q_{(X,\bar{g})}(u_{\epsilon}) \leq
\frac{C^{\prime\prime}\epsilon}{\epsilon^{\frac{2}{n}}} =
C^{\prime\prime}\epsilon^{\frac{n-2}{n}} \to 0 
$$ 
as $\epsilon\to 0$, where $C^{\prime\prime}>0$ is a constant. Thus
$Y_{\bar{C}}^{\cyl}(X)\leq 0$.
\end{Proof}
We also notice the following fact. The proof is straightforward.
\begin{Lemma}\label{conf:new:L2}
Assume that $\lambda({\cal L}_{h})=0$ and $n\geq 4$. Then either
$R_h\equiv 0$ or there exists a metric $\check{h}\in [h]$ such that
$$
\{
\begin{array}{l}
R_{\check{h}}\equiv \const. > 0,
\\
\lambda({\cal L}_{\check{h}})>0.
\end{array}
\right.
$$ 
\end{Lemma}
Consider a canonical cylindrical manifold $(Z\times \R, h+
dt^2)$. Then the cylindrical Yamabe constant $Y^{\cyl}_{[h+ dt^2]}(Z
\times \R)$ provides a universal upper bound for all $h$-cylindrical
Yamabe invariants. In more detail, we prove the following.
\begin{Proposition}\label{cyl:univ-bound}
{\rm (Generalized Aubin's inequality)} \ Let \ $X$ \ be an open
manifold with tame ends $Z\times [0,\infty)$ and $h\in \Riem(Z)$
any metric. Then $Y^{\cyl}_{[\bar{g}]}(X) \leq Y^{\cyl}_{[h+
dt^2]}(Z\times \R)$ for any metric $\bar{g}\in \Riem^{\cyl}(X)$ with
$\p_{\infty}\bar{g}=h$. Furthermore, $Y^{h\hbox{-}\cyl}(X) \leq
Y^{\cyl}_{[h+ dt^2]}(Z\times \R)$.
\end{Proposition}
\begin{Proof}
If $\lambda({\cal L}_h)<0$, the assertion is obvious from Lemma
\ref{cyl:L5}. Now assume that $\lambda({\cal L}_h)\geq 0$. Then
$Y^{\cyl}_{[h+ dt^2]}(Z\times \R)>-\infty$ by Lemma \ref{cyl:L6}. From
Lemma \ref{cyl:L4}, there exists a sequence $\{u_i\}$ of functions
$u_i\in C^{\infty}_c(Z\times \R)$ with $u_i\not\equiv 0$ such that
$$
\{
\begin{array}{l}
\supp (u_i) \subset Z\times [-i,i],
\\
\displaystyle
\lim_{i\to\infty} Q_{(Z\times\R,h+dt^2)}(u_i)= Y^{\cyl}_{[h+ dt^2]}(Z\times \R).
\end{array}
\right.
$$
We define the functions $\tilde{u}_i\in C^{\infty}_c(X)$ as follows:
$$
\tilde{u}_i(p)= 
\{
\begin{array}{ll}
0 & \mbox{if} \ \ p\in X(1), \\ u_i(x,t-i-1) & \mbox{if} \ \
p=(x,t)\in Z\times [1,\infty).
\end{array}
\right.
$$
Then we have 
$$
Y^{\cyl}_{[\bar{g}]}(X) \leq \lim_{i\to\infty}
Q_{(X,\bar{g})}(\tilde{u}_i) = \lim_{i\to\infty} Q_{(Z\times
\R,h+dt^2)}(u_i) = Y^{\cyl}_{[h+ dt^2]}(Z\times \R).
$$
This completes the proof of Proposition \ref{cyl:univ-bound}.
\end{Proof}
Now we would like to give an upper bound for the constant
$Y^{\cyl}_{[h+ dt^2]}(Z\times \R)$ in terms of the geometry of
$(Z,h)$.
\begin{Proposition}\label{cyl:univ-bound-2} Assume that $Z$ is connected.
The following inequality holds:
$$
Y^{\cyl}_{[h+ dt^2]}(Z\times \R)\leq Y(S^n)\cdot 
\min\{
1,\left(\frac{\lambda({\cal L}_h)}{(n-1)(n-2)}\right)^{\!\!\frac{n-1}{n}}
\!\!\!\!\!\cdot\!\!
\left(\frac{\Vol_h(Z)}{\Vol_{h_+}(S^{n-1})}\right)^{\!\!\frac{2}{n}}
\}.
$$
\end{Proposition}
\begin{Proof}
Lemmas \ref{cyl:L5} and \ref{cyl:L6}(2) imply that it is enough to
consider the case $\lambda({\cal L}_h)>0$. 
\begin{Claim}\label{cyl:univ-bound-3}
Assume that $\lambda({\cal L}_h) = (n-1)(n-2)$. Then
$$
Y^{\cyl}_{[h+ dt^2]}(Z\times \R)\leq Y(S^n)\cdot 
\left(\frac{\Vol_h(Z)}{\Vol_{h_+}(S^{n-1})}\right)^{\frac{2}{n}}.
$$
\end{Claim}
\begin{Proof}
Set $(X,\bar{h}) = (Z\times \R,h+ dt^2)$.  By the assumption, there
exists a function $u\in C_+^{\infty}(Z)$ such that
$$
\{
\begin{array}{l}
{\cal L}_h u = (n-1)(n-2) u,
\\
\int_Z u^2 d\sigma_h =1 .
\end{array}
\right.
$$
We consider the function $\phi(x,t)= f(t)\cdot u(x) \in
C_c^{\infty}(Z\times \R)$ for any $f\in C^{\infty}_c(\R)$ with
$f(t)\not\equiv 0$. Then we have
\begin{equation}\label{lambda1}
\begin{array}{rcl}
\displaystyle
Y^{\cyl}_{[\bar{h}]}(X) &\leq& \displaystyle
\inf_{^{\begin{array}{c}{^{f\in C^{\infty}_c(\R)}}
\\{^{f\not\equiv 0}}\end{array}}}\!\!\! Q_{(X,\bar{h})}(\phi)
\\
\\
 &\leq& \displaystyle
\inf_{\begin{array}{c}{^{f\in C^{\infty}_c(\R)}}
\\{^{f\not\equiv 0}}\end{array}} \!\!\! 
\frac{\int_X\left(\alpha_n |f^{\prime}|^2 u^2 + (n-1)(n-2) f^2 u^2 \right)d\sigma_h dt}
{\left(\int_X |f\cdot u|^{\frac{2n}{n-2}} d\sigma_h dt\right)^{\frac{n-2}{n}}}.
\end{array}\!\!\!\!\!\!
\end{equation}
We notice that
$$
\begin{array}{c}
\displaystyle
1= \int_Z u^2 d\sigma_h \leq 
\left(\int_Z | u|^{\frac{2n}{n-2}} d\sigma_h \right)^{\frac{n-2}{n}}\cdot 
\Vol_h(Z)^{\frac{2}{n}}, \ \ \ \ \mbox{and hence} \ \ \ \ 
\\
\\
\displaystyle
\frac{1}{\left(\int_Z |u |^{\frac{2n}{n-2}} d\sigma_h \right)^{\frac{n-2}{n}}} 
\leq \Vol_h(Z)^{\frac{2}{n}}.
\end{array}
$$
From this inequality combined with (\ref{lambda1}), we obtain
$$
\begin{array}{rcl}
\displaystyle
Y^{\cyl}_{[\bar{h}]}(X) &\leq& \displaystyle
\Vol_h(Z)^{\frac{2}{n}}\cdot \inf_{\begin{array}{c}{^{f\in C^{\infty}_c(\R)}}
\\{^{f\not\equiv 0}}\end{array}} \!\!\! 
\frac{\int_{\R}\left(\alpha_n |f^{\prime}|^2  + (n-1)(n-2) f^2  \right) dt}
{\left(\int_{\R} |f |^{\frac{2n}{n-2}}  dt\right)^{\frac{n-2}{n}}}
\\
\\
 &=& \displaystyle
Y(S^n)\cdot\left( \frac{\Vol_h(Z)}{\Vol_{h_+}(S^{n-1})}\right)^{\frac{2}{n}}.
\end{array}
$$
The last equality follows from Lemma \ref{cyl:L25} in Section
\ref{s3}.  This completes the proof of Claim \ref{cyl:univ-bound-3}.
\end{Proof}
We continue with the proof of Proposition \ref{cyl:univ-bound-2}. \ \ We
have that \ $\lambda({\cal L}_{\kappa^2\cdot h})= \kappa^{-2} \cdot
\lambda({\cal L}_{h}).  $ In particular, if 
$$
\kappa^2 = \frac{ \lambda({\cal L}_{h})}{(n-1)(n-2)}, \ \ \ \mbox{then}
\ \ \ \lambda({\cal L}_{\kappa^2\cdot h})=(n-1)(n-2).
$$
Also we have that $\Vol_{\kappa^2\cdot h}(Z) = \kappa^{n-1}
\Vol_{h}(Z)$. Then Claim \ref{cyl:univ-bound-3} combined with the
proof of Lemma \ref{cyl:L10} gives
$$
\begin{array}{rcl}
Y^{\cyl}_{[\bar{h}]}(Z\times \R) &\leq& \displaystyle Y(S^n) \cdot
\left( \frac{\Vol_{\kappa^2\cdot
h}(Z)}{\Vol_{h_+}(S^{n-1})}\right)^{\frac{2}{n}}
\\
\\
&=& 
Y(S^n) \cdot \left(\frac{ \lambda({\cal L}_{h})}{(n-1)(n-2)}\right)^{\frac{n-1}{n}}\cdot
\left( \frac{\Vol_{h}(Z)}{\Vol_{h_+}(S^{n-1})}\right)^{\frac{2}{n}}.
\end{array}
$$
\end{Proof}
Now let $M$ be a closed manifold of $\dim M= n\geq 3$ and
$p_{\infty}$ a point in $M$. \ Then the manifold \ $M\setminus
\{p_{\infty}\}$ \ is an open manifold with the tame end $Z\times
[0,\infty) =S^{n-1}\times [0,\infty)$.  Similarly, we start with an
open manifold $X$ with tame ends $Z\times [0,\infty)$, and
choose a finite number of points $p_1,\ldots,p_{k}\in X$.  Then \ 
$X^{\prime}=X\setminus\{p_1,\ldots,p_{k}\}$ \ is also an \ open \ manifold\ 
with \ tame ends $Z^{\prime}\times [0,\infty)$, where $Z^{\prime} = Z
\sqcup \left(S_1^{n-1}\sqcup\cdots \sqcup S_k^{n-1}\right)$ and
$S_1^{n-1}\sqcup\cdots \sqcup S_k^{n-1}$ denotes $k$ disjoint
copies of $S^{n-1}$.  We denote by $k\cdot h_+$ the metric
$h_+\sqcup\cdots\sqcup h_+$ on $S_1^{n-1}\sqcup\cdots \sqcup
S_k^{n-1}$. With these understood, we prove the following assertion.
\begin{Lemma}\label{cyl:L7} 
\begin{enumerate}
\item[{\bf (1)}] Let $M$ be a closed manifold, and $p_{\infty}\in M$.
Then
$$
Y^{\cyl}(M\setminus \{p_{\infty}\})\geq 
Y^{h_+\hbox{-}\cyl}(M\setminus \{p_{\infty}\}) = Y(M).
$$
\item[{\bf (2)}] Let $X$ be an open manifold with tame ends $Z\times
[0,\infty)$, and $p_1,\ldots,p_{k}\in X$. Then
$$
\begin{array}{rcl}
Y^{\cyl}(X\setminus\{p_1,\ldots,p_{k}\}) \!\!\!\!\!&\geq& \!\!\displaystyle
\sup_{h\in \Riem(Z)}\!\!
Y^{(h\sqcup (k \cdot h_+))\hbox{-}\cyl}(X\setminus\{p_1,\ldots,p_{k}\})
\\
\\ 
&= &\!\! Y^{\cyl}(X).
\end{array}
$$
\end{enumerate}
\end{Lemma}
\begin{Proof}
{\bf (1)} It is well-known from \cite{K2} that for any
$\epsilon>0$ there exists a conformal class $C_{\epsilon}\in {\cal
C}(M)$ such that
$$
\{
\begin{array}{l}
Y_{C_{\epsilon}}(M) \geq Y(M)-\epsilon, 
\\ 
C_{\epsilon} \ \ \mbox{is
locally conformally flat near $p_{\infty}$.}
\end{array}
\right. 
$$
This implies that there exists a cylindrical metric $\bar{g}_{\epsilon}\in
\Riem^{\cyl}(M\setminus \{p_{\infty}\})$ with $\p_{\infty}
\bar{g}_{\epsilon} = h_+$ such that
$$
Y^{\cyl}_{[\bar{g}_{\epsilon}]}(M\setminus \{p_{\infty}\}) = 
Y_{C_{\epsilon}}(M) \geq Y(M)-\epsilon.
$$
This gives
$$
Y^{\cyl}(M\setminus \{p_{\infty}\}) \geq Y^{h_+\hbox{-}\cyl}(M\setminus
\{p_{\infty}\})\geq Y(M).
$$
On the other hand,
$$
Y^{h_+\hbox{-}\cyl}(M\setminus
\{p_{\infty}\})\ \ = \!\!\!\!\!\sup_{^{^{\begin{array}{c}
_{\bar{g}\in \Riem^{\cyl}(M\setminus \{p_{\infty}\})} 
\\
_{\p_{\infty} \bar{g}= h_+}
\end{array}}}} \!\!\!\!\!
Y^{\cyl}_{[\bar{g}]}(M\setminus\{p_{\infty}\}) \leq Y(M).
$$
Therefore, $Y^{h_+\hbox{-}\cyl}(M\setminus \{p_{\infty}\})=Y(M)$.
The proof of the assertion {\bf (2)} is similar. Hence we omit it.
\end{Proof}
Lemma \ref{cyl:L7} shows that the $h$-cylindrical Yamabe invariant is
a natural generalization of the Yamabe invariant for closed manifolds.
Clearly the dependence of the invariant $Y^{h\hbox{-}\cyl}(X)$ on the
slice metric $h$ is important. There are the following natural
questions.
\begin{enumerate}
\item[{\bf (1)}] Is it true that $Y^{\cyl}(M\setminus
\{p_{\infty}\})=Y(M)$?
\item[{\bf (2)}] Is it true that $Y^{\cyl}(X\setminus\{p_1,\ldots,p_{k}\}) = Y^{\cyl}(X)$?
\end{enumerate}
It is easy to see that if $M=S^{n}$ and
$X=S^{n}\setminus\{\mbox{finite number of points}\}$, then (1) and (2)
hold. We also note the following.
\begin{Claim}
Let $M$ be a closed enlargeable manifold {\rm (}see {\rm
\cite{GL2}}{\rm )} with $Y(M)=0$, and $X:= M\setminus\{\mbox{finite
number of points}\}$. Then
$$
Y^{\cyl}(M\setminus \{p_{\infty}\})=Y(M)\ \ \ \mbox{and} \ \ 
Y^{\cyl}(X\setminus\{p_1,\ldots,p_{k}\}) = Y^{\cyl}(X).
$$
\end{Claim}
\section{Kobayashi-type inequalities}\label{s2b}
\begin{figure}[h]
\hspace*{25mm}\PSbox{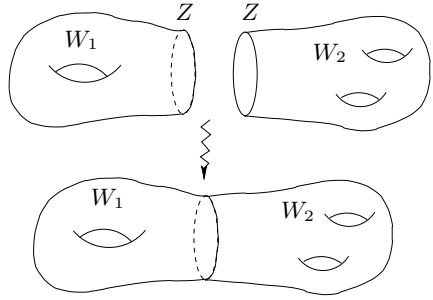}{6cm}{4cm}
\begin{picture}(0,0)
\put(-152,95){{\small $W_1$}}
\put(-58,90){{\small $W_2$}}
\put(-85,105){{\small $Z$}}
\put(-110,105){{\small $Z$}}
\put(-142,35){{\small $W_1$}}
\put(-70,30){{\small $W_2$}}
\end{picture}
\caption{Manifold $W=W_1\cup_Z W_2$.}\label{fig3-1}
\end{figure}

\begin{figure}[h]
\hspace*{25mm}\PSbox{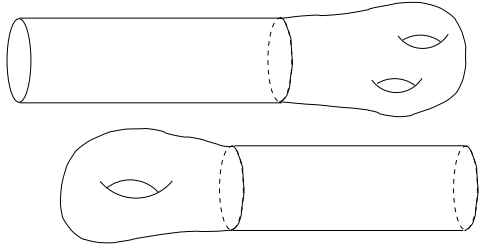}{6cm}{3.3cm}
\vspace{2mm}
\caption{Cylindrical manifolds $X_1$ and $X_2$.}\label{fig3-2}
\end{figure}

{\bf \ref{s2b}.1. Gluing construction.} Let $W_1$, $W_2$ be two
compact connected manifolds with $\p W_1 =\p W_2 = Z \neq
\emptyset$. Then let $W=W_1\cup_Z W_2$ be the union of manifolds $W_1$
and $W_2$ along their common boundary $Z$, see Fig. \ref{fig3-1}.
We define the corresponding open manifolds with tame ends
(see Fig. \ref{fig3-2}.2) by
$$
\{
\begin{array}{rl}
X_1 := & W_1\cup_Z Z\times [0,\infty),
\\
X_2 := & W_2\cup_Z Z\times [0,\infty).
\end{array}
\right.
$$
We denote by $X_1\sqcup X_2$ the disjoint union of $X_1$ and $X_2$.  The
following result is analogous to \cite[Lemma 1.10]{K2}. The
proof is similar to \cite{K2}.
\begin{Lemma}\label{cyl:L8}
For $\bar{C}_i\in {\cal C}^{\cyl}(X_i)$, $i=1,2$, let
$\bar{C}_1\sqcup\bar{C}_2$ denote the disjoint union of the conformal
classes $\bar{C}_1$ and $\bar{C}_2$ on $X_1\sqcup X_2$. Then
$$
Y^{\cyl}_{\bar{C}_1\sqcup\bar{C}_2}\!(X_1\!\sqcup X_2)\!=\!\!
\{\!\!
\begin{array}{cl}
\!\!-\!\left(\!
|Y^{\cyl}_{\bar{C}_1}\!(X_1)|^{\!\frac{n}{2}}\! + \!
|Y^{\cyl}_{\bar{C}_2}\!(X_2)|^{\!\frac{n}{2}}\! 
\!\right)^{\!\frac{2}{n}} & \!\!\!\!\!\mbox{if} \  Y^{\cyl}_{\bar{C}_i}(X_i)
\leq 0, \ i=1,2,
\\
\!\min\{Y^{\cyl}_{\bar{C}_1}(X_1), Y^{\cyl}_{\bar{C}_2}(X_2)\} & 
\!\!\!\mbox{otherwise}.
\end{array}
\right.
$$
\end{Lemma}
Lemma \ref{cyl:L8} implies the following assertions.
\begin{Proposition}\label{cyl:Th4}
Let $X_1$, $X_2$ be the open manifolds with tame ends as above and
$h\in \Riem(Z)$ any metric. Then
$$
\!\!
\begin{array}{lcl}
\displaystyle
Y^{h\hbox{-}\cyl}\!(X_1\sqcup X_2)\!\!\! \!\!\! 
&=&\!\!\!\!\! \displaystyle
\{
\begin{array}{ll}\!\!\!\!
-\!\left(
|Y^{h\hbox{-}\cyl}\!(X_1)|^{\frac{n}{2}}\!+\!
|Y^{h\hbox{-}\cyl}\!(X_2)|^{\frac{n}{2}} 
\!\right)^{\frac{2}{n}}\!\!\! &\!\!\!  \mbox{if} \ \  Y^{h\hbox{-}\cyl}(X_i)
\!\leq \!0, 
\\ &\ \ \ \ \ \ \ \ \ \ \  i = 1,2,
\\
\!\!\!
\min\{Y^{h\hbox{-}\cyl}(X_1), Y^{h\hbox{-}\cyl}(X_2)\} &\!\!\! 
 \mbox{otherwise}.
\end{array}
\right.
\\
\\
\displaystyle
Y^{\cyl}\!(X_1\!\sqcup X_2)\!\! \!\!\! &=&\!\!\!\!\! \displaystyle
\{
\begin{array}{ll}
\!\!\!\!
-\!\left(\!
|Y^{\cyl}\!(X_1)|^{\!\frac{n}{2}} \!+\!
|Y^{\cyl}\!(X_2)|^{\!\frac{n}{2}} 
\right)^{\!\frac{2}{n}} \!\!\! & \!\!\!\!\!\!\!\!\!\!\!\!\mbox{if} \  
Y^{\cyl}(X_i)\!\leq \!0,
\ i\!= \!1,2,
\\
\min\{Y^{\cyl}(X_1), Y^{\cyl}(X_2)\} \ \ \ \ \ \ & \!\!\! \mbox{otherwise}.
\end{array}
\right.
\end{array}
$$
\end{Proposition}
We denote $\Riem^*(Z)= \{ h\in \Riem(Z) \ | \ \lambda({\cal L}_h)>0
\}$.  We notice that the space 
$\Riem^+(Z)$ of positive scalar curvature metrics on $Z$ is contained
in $\Riem^*(Z)$, however  $\Riem^+(Z)\subsetneqq\Riem^*(Z) $.
\begin{Proposition}\label{cyl:Th5}
Let $W_1$, $W_2$ be compact manifolds of $\dim W_i =n\geq 3$ with $\p
W_1 = Z = \p W_2$ and $X_1$, $X_2$ the corresponding open manifolds
with tame ends as above.  Assume that $\Riem^*(Z)\neq
\emptyset$. Let $h\in \Riem^*(Z)$ be any metric. Then
\begin{equation}\label{cyl:eq6}
\!\!\!Y(W_1\cup_Z W_2) \!\geq\!
\{\!
\begin{array}{cl}\!\!\!
-\!\left(
|Y^{h\hbox{-}\!\cyl}\!(X_1)|^{\frac{n}{2}}\!\! + \!
|Y^{h\hbox{-}\!\cyl}\!(X_2)|^{\frac{n}{2}} 
\right)^{\!\frac{2}{n}} \!\!& \!\!\! \mbox{if} \   Y^{h\hbox{-}\cyl}\!(X_i)
\!\leq \!0, 
\\ &\ \ \ \ \ \ \ \ \ i = 1,2,
\\
\!\!\!
\min\{Y^{h\hbox{-}\cyl}(X_1), Y^{h\hbox{-}\cyl}(X_2)\} & \!\!\! 
\mbox{otherwise}.
\end{array}
\right.\!\!\!\!\!\!\!\!\!\!
\end{equation}
\end{Proposition}
{\bf Remark.} We notice that if $\lambda({\cal L}_h)<0$, then
$Y^{h\hbox{-}\cyl}(X_i) =-\infty$, $i=1,2$, and hence the inequality
(\ref{cyl:eq6}) holds. We study the case when $\lambda({\cal L}_h)=0$
later.
\begin{Proof}
Let $h\in \Riem^*(Z)$ be any metric. First we recall that Lemma
\ref{cyl:L6} gives that $Y^{h\hbox{-}\cyl}(X_i)>-\infty$ for $i=1,2$.

For any $\epsilon>0$ there exists a metric $\bar{g}_i\in
\Riem^{\cyl}(X_i)$ for  $i=1,2$ such that
$$
\{
\begin{array}{l}
\p_{\infty} \bar{g}_1 = h = \p_{\infty} \bar{g}_2, \ \ \mbox{and}
\\
Y^{\cyl}_{\bar{C}_1\sqcup\bar{C}_2}(X_1\sqcup X_2) \geq 
Y^{h\hbox{-}\cyl}(X_1\sqcup X_2) -\epsilon.
\end{array}
\right.
$$
Here $\bar{C}_i=[\bar{g}_i]\in {\cal C}^{\cyl}(X_i)$, $i=1,2$. 

We notice that the condition $\lambda({\cal L}_h)>0$ implies that
$\lambda({\cal L}_h^{\delta})>0$ for sufficiently small $\delta>0$,
where ${\cal L}_h^{\delta}$ is the operator defined by
$$
{\cal L}_h^{\delta} := -\left(\alpha_n -\delta\right)\Delta_h + R_h.
$$
We construct a family of metrics $\bar{g}(\ell)$ as below on the manifold
$W=W_1\cup_Z W_2$ by identifying it with the manifold
$$
\begin{array}{c}
\displaystyle
W_1\cup_Z W_2 = X_1(1) \cup_Z \left(Z\times [0,\ell]\right)\cup_Z X_2(1);
\\
\\
\bar{g}(\ell):= \bar{g}_1|_{X_1(1)} \cup (h+dt^2) \cup
\bar{g}_1|_{X_1(1)}; \ \ \ \mbox{(see  Fig. \ref{fig3-3}).}
\end{array}
$$
\begin{figure}[h]
\hspace*{2mm}\PSbox{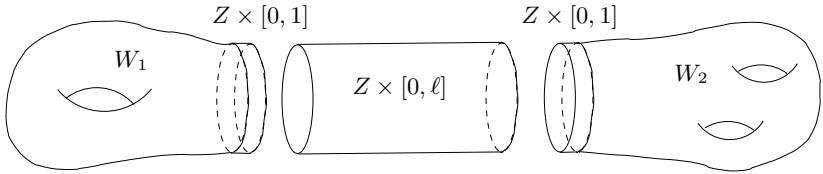}{6cm}{2.6cm}
\begin{picture}(0,0)
\put(-132,40){{\small $W_1$}}
\put(80,35){{\small $W_2$}}
\put(-42,30){{\small $Z\times [0,\ell]$}}
\put(-95,57){{\small $Z\times [0,1]$}}
\put(22,57){{\small $Z\times [0,1]$}}
\end{picture}

\caption{Decomposition of $W=W_1\cup_Z W_2$.}\label{fig3-3}
\end{figure}
By the definition, we have
$$
Y_{[\bar{g}(\ell)]}(W_1\cup_Z W_2) =\inf_{u>0} Q_{(W_1\cup_Z W_2,\bar{g}(\ell))}(u),
$$
and hence that for any $\ell>0$ there exists a function
$u_{\ell}\in C_+^{\infty}(W_1\cup_Z W_2)$ such that
$$
\begin{array}{c}
\displaystyle
\alpha_n \!
\int_{W_1\cup_Z W_2}\!\!\!\!\!|du_{\ell}|^2 d\sigma_{\bar{g}(\ell)}\! +\!
\int_{W_1\cup_Z W_2}\!\!\!\!\! 
R_{\bar{g}(\ell)} u_{\ell}^2 d\sigma_{\bar{g}(\ell)}
\leq Y_{[\bar{g}(\ell)]}(W_1\cup_Z W_2)\! +\!\frac{1}{1+\ell}, 
\\
\\
\displaystyle
\ \ \mbox{and} \ \ \ 
\int_{W_1\cup_Z W_2} |u_{\ell}|^{\frac{2n}{n-2}} d\sigma_{\bar{g}(\ell)} = 1 .
\end{array}
$$
For simplicity, set
$$
\begin{array}{rcl}
\displaystyle
E_{(W_1\cup_Z W_2,\bar{g}(\ell))}(u_{\ell})\!\!\!&:=&\!\!\! \displaystyle
\alpha_n 
\int_{W_1\cup_Z W_2}|du_{\ell}|^2 d\sigma_{\bar{g}(\ell)} +
\int_{W_1\cup_Z W_2} R_{\bar{g}(\ell)} u_{\ell}^2 d\sigma_{\bar{g}(\ell)},
\\
\\
Y_{[\bar{g}(\ell)]}\!\!\! &:=&\!\!\! Y_{[\bar{g}(\ell)]}(W_1\cup_Z W_2).
\end{array}
$$
\begin{Claim}\label{cyl:L9}
There exists $t_{\ell}$, $0\leq t_{\ell} \leq \ell$ such that
$$
\int_{Z\times \{ t_{\ell} \}} (|du_{\ell}|^2 + u_{\ell}^2)d\sigma_h \leq \frac{B}{\ell}
$$
for some positive constant $B$ independent of $\ell$.
\end{Claim}
\begin{Proof}
Using the inequality
$$
E_{(W_1\cup_Z W_2,\bar{g}(\ell))}(u_{\ell}) \leq Y_{[\bar{g}(\ell)]} +
\frac{1}{1+\ell},
$$
we obtain the following
$$
\begin{array}{rcl}
\displaystyle
Y_{[\bar{g}(\ell)]}\! +\!\frac{1}{1+\ell} \!\!\!&\geq& \!\!\!
\displaystyle
E_{(X(1), \bar{g}_1)}(u_{\ell})\! +\! E_{(X(2), \bar{g}_2)}(u_{\ell})\! +\!
\lambda({\cal L}_h^{\delta})\int_{Z\times [0,\ell]}\!\!\!
u_{\ell}^2 d\sigma_h dt 
\\
\\
&&\ \!\!\! \!\!\! 
\displaystyle
+  \ \delta \int_{Z\times [0,\ell]} \!\!\!|du_{\ell}|^2 d\sigma_h \ dt\!
+\! (\alpha_n\! -\!\delta) \int_{Z\times [0,\ell]} \!\!\!
|\p_t u_{\ell}|^2 d\sigma_h  dt .
\end{array}\!\!\!\!\!\!\!\!\!\!\!\!\!\!
$$
We use the fact that the last term on the right-hand side is positive,
and then
\begin{equation}\label{cyl:eq9}
\begin{array}{rcl}
\displaystyle
\!\!Y_{[\bar{g}(\ell)]} +\frac{1}{1+\ell}\!\!\! &\geq& \!\!\!
\displaystyle
(\min_{X_1(1)} R_{\bar{g}_1}^-)\!\cdot \!
\Vol_{\bar{g}_1}(X_1(1))^{\frac{2}{n}}\! +\!
(\min_{X_2(1)} R_{\bar{g}_2}^-)\!\cdot \!
\Vol_{\bar{g}_2}(X_2(1))^{\frac{2}{n}} 
\\
\\
&& 
\displaystyle
\ \ \ \ \ + \ \delta_0\cdot \int_{Z\times [0,\ell]} 
\left(|du_{\ell}|^2+u_{\ell}^2\right)d\sigma_h  \ dt .
\end{array}\!\!\!\!\!\!\!\!\!\!\!
\end{equation}
Here $\delta_0:=\min\{\lambda({\cal L}_h^{\delta}),\delta\}>0$ and
$R^- := \min\{0,R\}$. Hence (\ref{cyl:eq9}) implies that
$$
\begin{array}{l}
\displaystyle
\delta_0\cdot \int_{Z\times [0,\ell]} 
\left(|du_{\ell}|^2+u_{\ell}^2\right)d\sigma_h  \ dt 
\leq Y_{[\bar{g}(\ell)]} +\frac{1}{1+\ell}  + A \ \ \ \ \ \ \ \mbox{with} 
\\
\\
\displaystyle
A := \!-\[\! (\min_{X_1(1)} R_{\bar{g}_1}^{-})\!\cdot \!
\Vol_{\bar{g}_1}(X_1(1))^{\!\frac{2}{n}}\! + \!
(\min_{X_2(1)} R_{\bar{g}_2}^-)\!\cdot \!
\Vol_{\bar{g}_2}(X_2(1))^{\!\frac{2}{n}}\]\! \geq\! 0,
\end{array}
$$
and thus
$$
\int_{Z\times [0,\ell]} \left(|du_{\ell}|^2+u_{\ell}^2\right)d\sigma_h  \ dt
\leq \frac{1}{\delta_0}\left(Y_{[\bar{g}(\ell)]} +\frac{1}{1+\ell}  + A 
\right).
$$
It then follows that there exists $t_{\ell}$, $0\leq t_{\ell} \leq
\ell$ such that
$$
\int_{Z\times \{ t_{\ell} \}} (|du_{\ell}|^2 + 
u_{\ell}^2)d\sigma_h \leq \frac{B}{\ell}
$$
with $B>0$.
\end{Proof}
Now we continue with the proof of Proposition \ref{cyl:Th5}.  We cut
the manifold $W_1\cup_Z W_2$ along the slice $Z\times \{t_{\ell}\}$,
and then attach two copies of the half-cylinder $Z\times [0,\infty)$
to the corresponding pieces to obtain the cylindrical manifolds $X_1$
and $X_2$.  In other words, we regard $X_1\sqcup X_2$ as
\begin{equation}\label{regard}
X_1\sqcup X_2 = (Z\times [0,\infty)) \cup_Z \left(
(W_1\cup_Z W_2) \setminus
(Z\times \{t_{\ell}\}) \right) \cup_Z (Z\times [0,\infty)) .
\end{equation}
We define a non-negative Lipschitz function
$$
U_{\ell}\in C_c^{0}(X_1\sqcup X_2) 
$$
by
\begin{equation}\label{cyl:eq10}
\begin{array}{l}
\displaystyle
U_{\ell} = u_{\ell} \ \ \ \mbox{on} \ \ (W_1\cup_Z W_2) \setminus
(Z\times \{t_{\ell}\}), \ \ \mbox{and}
\\
\\
U_{\ell}(x,t)= \{\begin{array}{cl}
(1-t)u_{\ell}(x,t_{\ell}) & \ \mbox{on} \  Z\times [0,1],
\\
0 & \ \mbox{on} \  Z\times [1,\infty).
\end{array}
\right.
\end{array}
\end{equation}
The conditions (\ref{cyl:eq10}) imply that
$$
\begin{array}{l}
\displaystyle
E_{(X_1\sqcup X_2, \bar{g}_1\sqcup \bar{g}_2)}(U_{\ell}) \leq
Y_{[\bar{g}(\ell)]} +\frac{C}{\ell},
\ \ \ \ \ 
\int_{X_1\sqcup X_2} U_{\ell}^{\frac{2n}{n-2}} 
d\sigma_{\bar{g}_1\sqcup \bar{g}_2} > 1 .
\end{array}
$$
This gives that
\begin{equation}\label{cyl:eq11}
\begin{array}{rcl}
\displaystyle
Y(W_1\cup_Z W_2) + \frac{C}{\ell} &\geq& \displaystyle
Y_{[\bar{g}(\ell)]} +\frac{C}{\ell} \geq
Q_{(X_1\sqcup X_2, \bar{g}_1\sqcup \bar{g}_2)}(U_{\ell})
\\
\\
\displaystyle
&\geq&\displaystyle
 \inf_{u\in C_c^{\infty}(X_1\sqcup X_2)}
Q_{(X_1\sqcup X_2, \bar{g}_1\sqcup \bar{g}_2)}(u) 
\\
\\
&=& 
Y^{\cyl}_{\bar{C}_1\sqcup\bar{C}_2}(X_1\sqcup X_2).
\end{array}\!\!\!\!\!\!\!\!\!\!\!\!\!\!\!\!\!\!
\end{equation}
Now we take $\ell\to\infty$ in (\ref{cyl:eq11}) to obtain
$$
Y(W_1\cup_Z W_2)\geq Y^{\cyl}_{\bar{C}_1\sqcup\bar{C}_2}(X_1\sqcup X_2)
\geq Y^{h\hbox{-}\cyl}(X_1\sqcup X_2) -\epsilon.
$$
Taking $\epsilon\to 0$, we conclude that
$$
\begin{array}{rcl}
\displaystyle
\!\!
Y(W_1\cup_Z W_2)\!\!\!\!\!&\geq& \!\!\!\!\! Y^{h\hbox{-}\cyl}(X_1\sqcup X_2) 
\\
\\
\displaystyle
&\geq& \!\!\!\!\!\!\
\{\!
\begin{array}{cl}\!\!\!
-\!\left(
|Y^{h\hbox{-}\cyl}\!(X_1)|^{\frac{n}{2}}\! +\! 
|Y^{h\hbox{-}\cyl}\!(X_2)|^{\frac{n}{2}}\! 
\right)^{\!\frac{2}{n}} \!\!\!&  \!\!\mbox{if} \   Y^{h\hbox{-}\cyl}\!(X_i)
\!\leq \!0, 
\\ & \ \ \ \ \ \ \ \ \ \ \ i\!=\!1,2,
\\
\!\!\!
\min\{Y^{h\hbox{-}\cyl}(X_1), Y^{h\hbox{-}\cyl}(X_2)\} & \mbox{otherwise}.
\end{array}
\right.\!\!\!\!\!  \end{array} $$ 
This completes the proof of Proposition \ref{cyl:Th5}.
\end{Proof}
Now we establish the Kobayashi-type inequality in the case when
$\lambda({\cal L}_{h})=0$. 
\begin{Proposition}\label{cyl:Th6}
Let $h\in \Riem(Z)$ be a metric with $\lambda({\cal L}_{h})=0$. Then
$$
Y(W_1\cup_Z W_2) \geq -\left(|Y^{h\hbox{-}\cyl}(X_1)|^{\frac{n}{2}}+ 
|Y^{h\hbox{-}\cyl}(X_2)|^{\frac{n}{2}}\right)^{\frac{2}{n}}.
$$
\end{Proposition}
\begin{Proof}
Let $h\in \Riem(Z)$ be any metric with $\lambda({\cal L}_{h})=0$.
Recall that this condition implies that $0\geq Y^{h\hbox{-}\cyl}(X_i)
> -\infty$ for our cylindrical manifold $X_i$, $i=1,2$. By the
definitions, for each $\epsilon>0$ there exists a metric $\bar{g}_i\in
\Riem^{\cyl}(X_i)$, $i=1,2$, such that
$$
\{
\begin{array}{l}
\p \bar{g}_1 = h = \p \bar{g}_2, \\ Y^{\cyl}_{\bar{C}_1\sqcup
\bar{C}_2}(X_1\sqcup X_2) \geq Y^{h\hbox{-}\cyl}(X_1\sqcup X_2) -\epsilon.
\end{array}
\right.
$$
Here $\bar{C}_i := [\bar{g}_i]\in {\cal C}^{\cyl}(X_i)$, $i=1,2$.  As
before, we use the decomposition
$$
W_1\cup_Z W_2 \cong X_1(1) \cup_Z (Z\times [0,\ell^2]) \cup_Z X_2(1)
$$
with the metric 
$$
\bar{g}(\ell^2) = \bar{g}_1|_{X_1(1)} \cup (h+dt^2) \cup \bar{g}_2|_{X_2(1)}
$$
respectively on $X_1(1)$, $Z\times [0,\ell^2]$ and $X_2(1)$. We set
the segments
$$
\{
\begin{array}{rcll}
I_j &:=& [(j-1)\ell, j\ell], & j=1,\ldots, \ell,
\\
I_{j,k} &:=& [(j-1)\ell+(k-1), (j-1)\ell+k ], & j,k=1,\ldots, \ell.
\end{array}
\right.
$$
We have
$$
Y_{[\bar{g}(\ell^2)]}:= Y_{[\bar{g}(\ell^2)]}(W_1\cup_Z W_2) =
\inf_{u>0} Q_{(W_1\cup_Z W_2,\bar{g}(\ell^2))}(u).
$$
Then for any $\ell>>1$ there exists $u_{\ell}\in C^{\infty}_+(W_1\cup_Z
W_2)$ such that
$$
\{
\begin{array}{l}
\displaystyle
E_{(W_1\cup_Z W_2,\bar{g}(\ell^2))}(u_{\ell}) \leq 
Y_{[\bar{g}(\ell^2)]} +\frac{1}{1+\ell}\ , 
\\
\\
\displaystyle
\int_{W_1\cup_Z W_2}u_{\ell}^{\frac{2n}{n-2}} d\sigma_{\bar{g}(\ell^2)} = 1.
\end{array}
\right.
$$
\begin{Claim}\label{cyl:L11}
There exists $t_{\ell}$, $0\leq t_{\ell}\leq \ell^2$ such that
$$
\int_{Z\times\{t_{\ell}\}}\left(|du_{\ell}|^2 + u_{\ell}^2\right) d\sigma_h 
\leq \frac{B}{\ell^{\frac{n-2}{n}}},
$$
where $B$ is independent of $\ell$.
\end{Claim}
\begin{Proof}
FRom the above, we have
$$
Y(W_1\cup_Z W_2)+\frac{1}{1+\ell} \geq  Y_{[\bar{g}(\ell^2)]} +\frac{1}{1+\ell} 
\geq E_{(W_1\cup_Z W_2,\bar{g}(\ell^2))}(u_{\ell})  .
$$
We use the same constant $A$ as in the proof of Claim \ref{cyl:L9} to
give the following estimate
$$
\begin{array}{rcl}
\displaystyle 
E_{(W_1\cup_Z W_2,\bar{g}(\ell^2))}(u_{\ell}) &\geq &
\displaystyle -A +
\int_{Z\times[0,\ell^2]}\left(\alpha_n |d u_{\ell}|^2 +
R_h u_{\ell}^2\right) d\sigma_h \ dt
\\
\\
&&\displaystyle 
\ \ \ \ \ \ \ \ 
\ + \ \alpha_n\int_{Z\times[0,\ell^2]} |\p_t u_{\ell}|^2 d\sigma_h \ dt
\\
\\
&\geq &
\displaystyle
-A +\alpha_n \int_{Z\times[0,\ell^2]} |\p_t u_{\ell}|^2 d\sigma_h \ dt
\end{array}\!\!\!\!\!\!\!\!\!\!\!\!\!\!\!\!\!\!
$$
since $\lambda({\mathcal L}_h)=0$.  Clearly there exists an integer $j$
($1\leq j\leq \ell$) such that
$$
\int_{Z\times I_j} u^{\frac{2n}{n-2}} d\sigma_h \ dt \leq \frac{1}{\ell}.
$$
Then there exists also an integer $k$ ($1\leq k\leq \ell$)  such that
$$
\begin{array}{l}
\displaystyle
\int_{Z\times I_{j,k}} \left(\alpha_n |du_{\ell}|^2 +
R_h u_{\ell}^2 \right) d\sigma_h \ dt 
\leq \frac{A^{\prime}}{\ell}, \ \ \ \ 
\mbox{where} \ \ \ \ 
\\
\\
\displaystyle
A^{\prime}:= A + Y(W_1\cup_Z W_2) +\frac{1}{1+\ell}.
\end{array}
$$
We have 
\begin{equation}\label{cyl:eq13}
\int_{Z\times I_{j,k}} u_{\ell}^2  d\sigma_h \ dt \leq \Vol_h(Z)^{\frac{2}{n}}
\cdot \left(\int_{Z\times I_{j,k}}  u_{\ell}^{\frac{2n}{n-2}} 
d\sigma_h \ dt\right)^{\frac{n-2}{n}} \leq \frac{C_1}{\ell^{\frac{n-2}{n}}} 
\end{equation}
for some positive constant $C_1$. We define the function $f(t)$ by
$$
f(t):= \int_{Z\times \{t\}} u_{\ell}^2  d\sigma_h, \ \ \ \mbox{where} \ \ \ 
f^{\prime}(t)= 2 \int_{Z\times \{t\}} u_{\ell}\cdot \p_{t}u_{\ell}  d\sigma_h.
$$
It follows from (\ref{cyl:eq13}) that there exists $t_0\in I_{j,k}$ such that
\begin{equation}\label{cyl:eq13a}
f(t_0) \leq \frac{C_1}{\ell^{\frac{n-2}{n}}} .
\end{equation}
We have the following estimate
$$
\begin{array}{rcl}
\displaystyle
|f(t_2)-f(t_1)| &\leq &
\displaystyle
2\left| \int_{Z\times [t_1,t_2]} u_{\ell}\cdot \p_{t}u_{\ell}  
d\sigma_h \ dt\right|
\\
\\
&\leq & 
\displaystyle
2\left( \int_{Z\times [t_1,t_2]}\!\! u_{\ell}^2  
d\sigma_h \ dt\right)^{\frac{1}{2}}
\!\cdot\!
\left( \int_{Z\times [t_1,t_2]} \!\!|\p_t u_{\ell}|^2  
d\sigma_h \ dt\right)^{\frac{1}{2}}
\\
\\
&\leq &
\displaystyle
2\cdot\frac{C_1^{\frac{1}{2}}}{\ell^{\frac{n-2}{2n}}} \cdot 
\frac{A^{\prime\prime}}{\ell^{\frac{1}{2}}} 
\ = \ 
\displaystyle
\frac{2A^{\prime\prime}C_1^{\frac{1}{2}}}{\ell^{\frac{n-1}{n}}}
\end{array}
$$
for any $t_1,t_2\in I_{j,k}$. Here $C_1>0$ is the same constant as
above.  This combined with (\ref{cyl:eq13a}) gives that
\begin{equation}\label{cyl:eq14}
f(t)\leq
\frac{2A^{\prime\prime}C_1^{\frac{1}{2}}}{\ell^{\frac{n-1}{n}}} +
\frac{C_1}{\ell^{\frac{n-2}{n}}} \leq
\frac{C_2}{\ell^{\frac{n-2}{n}}}.
\end{equation}
for any $t\in I_{j,k}$. Now  we have the estimate
$$
\begin{array}{rcl}
\displaystyle
\alpha_n \int_{Z\times I_{j,k}} |d u_{\ell}|^2  d\sigma_h \ dt &\leq &
\displaystyle
\frac{A^{\prime}}{\ell} + \left(\max_Z |R_h| \right) \cdot 
\int_{Z\times I_{j,k}} u_{\ell}^2 \ d\sigma_h \ dt
\\
\\
&\leq & \displaystyle 
\frac{A^{\prime}}{\ell} + \frac{C_3}{\ell^{\frac{n-2}{n}}}.
\end{array}
$$
Here we used the estimate (\ref{cyl:eq13}). Then this gives
$$
\int_{Z\times I_{j,k}} | d u_{\ell}|^2  d\sigma_h \ dt 
\leq \frac{C_4}{\ell^{\frac{n-2}{n}}}
$$
for some positive constant $C_4$. We conclude that
there exists $t_{\ell}\in I_{j,k}$ such that
\begin{equation}\label{cyl:eq15}
\int_{Z\times \{t_{\ell}\}} | d u_{\ell}|^2 d\sigma_h \leq
\frac{C_4}{\ell^{\frac{n-2}{n}}}.
\end{equation}
Now we use the estimates (\ref{cyl:eq14}) and (\ref{cyl:eq15}) to obtain 
the estimate
$$
\int_{Z\times \{t_{\ell}\}} \left( | d u_{\ell}|^2 + u_{\ell}^2 \right) d\sigma_h 
\leq \frac{C_2+C_4}{\ell^{\frac{n-2}{n}}}.
$$
This completes the proof of Claim \ref{cyl:L11}.
\end{Proof}
We continue now with the proof of \ Proposition \ref{cyl:Th6}.  \ As
above, \ we regard $X_1\sqcup X_2$ as in (\ref{regard}).  \ We define a
non-negative \ Lipschitz function \ $ U_{\ell}\in C_c^0(X_1\sqcup X_2) $
by
$$
\begin{array}{l}
U_{\ell} = u_{\ell} \ \ \mbox{on} \ \ 
(W_1\cup_Z W_2) \setminus (Z\times \{t_{\ell}\}), \ \ 
\mbox{and}
\\
\\
U_{\ell}(x,t) =
\{
\begin{array}{cl}
(1-t)u_{\ell}(x,t_{\ell}) & \mbox{on} \ \ Z\times [0,1],
\\
0 & \mbox{on} \ \ Z\times [1,\infty).
\end{array}
\right.
\end{array}
$$
Then we obtain the following estimate
$$
\{
\begin{array}{l}
\displaystyle 
E_{(X_1\sqcup X_2, \bar{g}_1\sqcup\bar{g}_2)}(U_{\ell}) \leq Y(W_1\cup_Z W_2) + 
\frac{C}{\ell^{\frac{n-2}{n}}}, \ \ \ \mbox{with}  
\\
\\
\displaystyle 
\int_{X_1\sqcup X_2} U_{\ell}^{\frac{2n}{n-2}} d\sigma_{\bar{g}_1\sqcup\bar{g}_2} > 1.
\end{array}
\right.
$$
This implies
$$
\begin{array}{rcl}
\displaystyle 
Y(W_1\cup_Z W_2) +  \frac{C}{\ell^{\frac{n-2}{n}}}
& \geq & Q_{(X_1\sqcup X_2, \bar{g}_1\sqcup\bar{g}_2)}(U_{\ell}) 
\\
\\
& \geq &\displaystyle 
\!\!\!\!\!\!\!
\inf_{u\in C_c^{\infty}(X_1\sqcup X_2)} Q_{(X_1\sqcup X_2, \bar{g}_1\sqcup\bar{g}_2)}(u)
\\
\\
&= & Y^{\cyl}_{\bar{C}_1\sqcup \bar{C}_2}(X_1\sqcup X_2).
\end{array}
$$ 
We take $\ell\to\infty$ to obtain
$$
Y(W_1\cup_Z W_2) \geq Y^{\cyl}_{\bar{C}_1\sqcup \bar{C}_2}(X_1\sqcup X_2) \geq
Y^{h\hbox{-}\cyl}_{\bar{C}_1\sqcup \bar{C}_2}(X_1\sqcup X_2) -\epsilon.
$$
Finally, we take $\epsilon\to 0$ and conclude that
$$
Y(W_1\cup_Z W_2) \geq Y^{h\hbox{-}\cyl}_{\bar{C}_1\sqcup
\bar{C}_2}(X_1\sqcup X_2).
$$
This combined with Proposition \ref{cyl:Th4} completes the
proof of Proposition \ref{cyl:Th6}.
\end{Proof}
Now we combine Proposition \ref{cyl:Th6} with Proposition \ref{cyl:Th5} 
and Lemma \ref{cyl:L5}:
\begin{Theorem}\label{cyl:Th7}
Let $W_1$, $W_2$ be compact manifolds of $\dim W_i\geq 3$ with $\p W_1
= Z = \p W_2$ and $X_1$, $X_2$ the corresponding
open manifolds with tame ends $Z\times [0,\infty)$ as above. Let $h\in
\Riem(Z)$ be any metric. Then
$$
\!\!\!
Y(W_1\cup_Z W_2) \!\geq\!
\{
\begin{array}{cl}\!\!\!
-\!\left(
|Y^{h\!\hbox{-}\!\cyl}\!(X_1)|^{\!\frac{n}{2}} + 
|Y^{h\!\hbox{-}\!\cyl}\!(X_2)|^{\!\frac{n}{2}} 
\right)^{\!\frac{2}{n}} \!\!\! & \!\! \mbox{if} \   
Y^{h\!\hbox{-}\!\cyl}\!(X_i)
\!\leq \!0, 
\\
& \ \ \ \ \ \ \ \ \ \ i\!=\!1,2,
\\
\!\!\!
\min\{Y^{h\hbox{-}\cyl}(X_1), Y^{h\hbox{-}\cyl}(X_2)\} & \mbox{otherwise}.
\end{array}
\right.\!\!\!\!\!
$$
\end{Theorem}
{\bf Remark.}  We proved a similar formula in \cite{AB1} in terms of
the \emph{relative Yamabe invariant}. 
\vspace{2mm}

\noindent
We notice that Theorem \ref{cyl:Th7} and Lemma \ref{cyl:L7} recover
the original Kobayashi inequality (see \cite[Theorem 2(a)]{K2}).
\begin{Corollary}\label{cyl:Th8} 
Let $M_1$, $M_2$ be closed manifolds of $\dim M_i= n\geq
3$, $i=1,2$. Then
$$
Y(M_1\# M_2) \geq
\{
\begin{array}{cl}\!\!\!
-\left(
|Y(M_1)|^{\frac{n}{2}} + 
|Y(M_2)|^{\frac{n}{2}} 
\right)^{\frac{2}{n}} &  \mbox{if} \   Y(M_i)
\leq 0, \ i=1,2,
\\
\!\!\!
\min\{Y(M_1), Y(M_2)\} & \mbox{otherwise}.
\end{array}
\right.\!\!\!\!\!
$$
\end{Corollary}
\begin{Proof}
Theorem \ref{cyl:Th7} combined with Lemma \ref{cyl:L7} gives directly
the following 
$$
\begin{array}{rcl}
\displaystyle
Y(M_1\# M_2) \!\!\!
& = & \!\!\! Y((M_1\setminus D^n)\cup_{S^{n-1}} (M_2\setminus D^n)) 
\\
\\
&\geq& \!\!\!\!\!\! \displaystyle
\{
\begin{array}{l}
\!-\left(
|Y^{h_+\!\hbox{-}\!\cyl}(M_1\setminus \{pt\})|^{\frac{n}{2}}\! + \!
|Y^{h_+\!\hbox{-}\!\cyl}(M_2\setminus \{pt\})|^{\frac{n}{2}} 
\right)^{\frac{2}{n}} 
\\
\ \ \ \ \ \  \ \ \ \ \ \  \mbox{if} \   Y^{h_+\hbox{-}\cyl}(M_i\setminus \{pt\})
\leq 0, \ i=1,2,
\\
\\
\!\!\min\{Y^{h_+\!\hbox{-}\!\cyl}(M_1\setminus \{pt\}), 
Y^{h_+\!\hbox{-}\!\cyl}(M_2\setminus \{pt\})\} \  \mbox{otherwise}.
\end{array}
\right.\!\!\!\!\!
\\
\\
&=& \!\!\!\!\!\!\displaystyle
\{
\begin{array}{cl}\!\!\!
-\left(
|Y(M_1)|^{\frac{n}{2}} + 
|Y(M_2)|^{\frac{n}{2}} 
\right)^{\frac{2}{n}} &  \mbox{if} \   Y(M_i)
\leq 0, \ i=1,2,
\\
\!\!\!
\min\{Y(M_1), Y(M_2)\} & \mbox{otherwise}.
\end{array}
\right.\!\!\!\!\!
\end{array}
$$
\end{Proof}

\section{Surgery and cylindrical Yamabe invariant}\label{s3}
{\bf \ref{s3}.1.  Relative Yamabe constant for a cylindrical manifold
with boundary.}  Now we concentrate our attention on the case when a
cylindrical manifold has a nonempty boundary.  Let $X$ be a noncompact
manifold with $\p X= M\neq \emptyset$. As before, we assume that $X$
is a noncompact manifold with tame ends. More precisely, there exists
a relatively compact open submanifold $W\subset X$ with $\p W= \p X =
M$ and $\p \ov{W}= M\sqcup Z$ such that
$$
\{
\begin{array}{l}
\displaystyle
Z = \bigsqcup_{j=1}^m Z_j, \ \ \mbox{where each $Z_j$
is connected,}
\\
\displaystyle 
X\setminus W \cong Z\times [0,1)= \bigsqcup_{j=1}^m 
\( Z_j \times [0,1)\).
\end{array}
\right.
$$
We should refine the definition of suitable cylindrical metrics to the
case of a non-empty boundary.

Let $\bar{g}\in \Riem(X)$. We denote by $H_{\bar{g}}$ the mean
curvature of $\bar{g}$ on $M$. Then we define
$$
\Riem^{\cyl,0}(X):= \{ \bar{g}\in \Riem^{\cyl}(X) \ | \ H_{\bar{g}}
\equiv 0 \ \ \mbox{on} \ \ M \}.
$$
Let $\bar{C}\in {\cal C}^{\cyl}(X)$ be a cylindrical conformal class,
and $\bar{g}\in \bar{C}\cap \Riem^{\cyl,0}(X)$. We define the {\it
normalized $L^{k,2}_{\bar{g}}$-conformal class} for $k=1,2$ by
$$
[\bar{g}]_{L^{k,2}_{\bar{g}}}^0 := 
\bar{C}_{L^{k,2}_{\bar{g}}}^0 := \{ u^{\frac{4}{n-2}}\bar{g} \ | \ 
u \in C_+^{\infty}(X)\cap L^{k,2}_{\bar{g}}(X) , \ \left.\frac{\p u}{\p \nu}\right|_M = 0 \} .
$$
Here $\nu$ is the outward unit vector field normal to the boundary $
M$.  Recall that the functional
$$
Q_{(X,\bar{g})}(u) = \frac{\int_X\[\alpha_n|du|^2 +
R_{\bar{g}}u^2\] d\sigma_{\bar{g}}}{\(\int_X
|u|^{\frac{2n}{n-2}}d\sigma_{\bar{g}}\)^{\frac{n-2}{n}}}
$$
is well-defined on the space
$L_{\bar{g}}^{1,2}(X)$ with $u \not\equiv 0$.
\begin{figure}[h]
\hspace*{15mm}
\PSbox{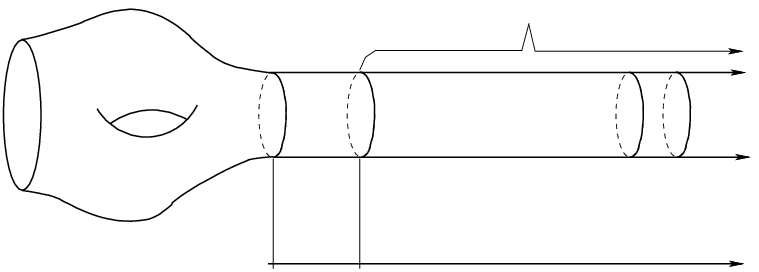}{10mm}{32mm}  
\begin{picture}(0,1)
\put(33,70){{\small $Z\!\times\!\!\{0\}$}}
\put(65,70){{\small $Z\!\times\!\!\{1\}$}}
\put(110,78){{\small $Z\!\times\!\! [1,\infty)$}}
\put(45,-9){{\small $0$}}
\put(70,-9){{\small $1$}}
\put(185,0){{\small $\R$}}
\put(-30,73){{\small $M$}}
\put(0,60){{\small $W$}}
\end{picture}
\caption{A cylindrical manifold $X$ with boundary.}\label{fig4-1}
\end{figure}

\noindent
Now we consider the functional 
$$
I_X(\tilde{g}):= \frac{\int_X R_{\tilde{g}}d\sigma_{\tilde{g}}}
{\Vol_{\tilde{g}}(X)^{\frac{n-2}{n}}}
$$
on the space of normalized $L^{2,2}_{\bar{g}}$-conformal metrics,
i.e. $\tilde{g}\in [\bar{g}]^0_{L^{2,2}_{\bar{g}}}$ with $\bar{g}\in
\Riem^{\cyl,0}(X)$. The following lemma is similar to Lemma
\ref{cyl:L2}.
\begin{Lemma}\label{cyl:L20}
The functional $I_X(\tilde{g})$ is well-defined for metrics
$\tilde{g}\in [\bar{g}]^0_{L^{2,2}_{\bar{g}}}$. Furthermore,
$I_X(\tilde{g}) = Q_{(X,\bar{g})}(u)$ if $\tilde{g} =
u^{\frac{4}{n-2}}\cdot\bar{g}$.
\end{Lemma}
Let $\bar{g}\in \Riem^{\cyl,0}(X)$. Then we define the
constant
$$
Y^{\cyl}_{\bar{g}}(X,M; [\bar{g}|_M]):= \inf_{\tilde{g}\in
[\bar{g}]^0_{L^{2,2}_{\bar{g}}}} I_X(\tilde{g}).
$$
The following result follows from \cite{Ch}.
\begin{Fact}\label{cyl:L21}
$ Y^{\cyl}_{\bar{g}}(X,M; [\bar{g}|_M]) \leq Y(S^n_+, S^{n-1})$, \ 
where \ $Y(S^n_+, S^{n-1})$ \ is the relative Yamabe invariant of
the round hemisphere $S^n_+$ with the equator $S^{n-1} = \p S^n_+$.
\end{Fact}
The following lemma is is an analogue of Lemma \ref{cyl:L4}. The proof
is essentially the same.
\begin{Lemma}\label{cyl:L22}{\rm (cf. \cite[Lemma 2.1]{Ak3})}
The following identities hold:
$$
\begin{array}{rcl}\displaystyle
Y^{\cyl}_{\bar{g}}(X,M; [\bar{g}|_M])& = & \displaystyle
\inf_{\tilde{g}\in
[\bar{g}]^0_{L^{2,2}_{\bar{g}}}} I_X(\tilde{g}) \ \ \ 
= 
\inf_{u\in C_+^{\infty}(X)\cap
L^{1,2}_{\bar{g}}(X)} Q_{(X,\bar{g})}(u)
\\
\\
& = & \displaystyle\!\!\!\!
\inf_{u\in  L^{1,2}_{\bar{g}}(X), u\not\equiv 0} Q_{(X,\bar{g})}(u)\ \  
= \!\!
\inf_{u\in  C_c^{\infty}(X), u\not\equiv 0} Q_{(X,\bar{g})}(u).
\end{array}
$$
\end{Lemma}
We conclude the following:
\begin{Claim}\label{cyl:L23}
Let $\bar{g}, \hat{g}\in \bar{C}\cap \Riem^{\cyl,0}(X)$ be any two metrics. Then
$$
Y^{\cyl}_{\bar{g}}(X,M; [\bar{g}|_M]) = Y^{\cyl}_{\hat{g}}(X,M; [\hat{g}|_M]).
$$
\end{Claim}
This allows us to define the relative cylindrical Yamabe constant. 
\begin{Definition}
{\rm Let \ $\bar{C}\in {\cal C}^{\cyl}(X)$ \ be \ a \ cylindrical \ conformal
class \ and \ $\bar{g}\in \bar{C}\cap \Riem^{\cyl,0}(X)$ any
cylindrical metric. Then the \emph{relative cylindrical Yamabe constant}
$Y^{\cyl}_{\bar{C}}(X,M; [\bar{g}|_M])$ is defined as}
$$
\begin{array}{rcl}
Y^{\cyl}_{\bar{C}}(X,M; \p \bar{C}) & := &
Y^{\cyl}_{\bar{g}}(X,M; [\bar{g}|_M]) =
\inf_{\tilde{g}\in
\bar{C}^0_{L^{2,2}_{\bar{g}}}} I_X(\tilde{g})
\\
\\
&=& 
\inf_{u\in  L^{1,2}_{\bar{g}}(X), u\not\equiv 0} Q_{(X,\bar{g})}(u).
\end{array}
$$
\end{Definition}
It then follows that $Y^{\cyl}_{\bar{C}}(X,M; \p \bar{C}) \leq
Y(S^n_+, S^{n-1})$. We remark here that in the case of compact
manifolds with non-empty boundaries, this definition gives the
\emph{relative Yamabe constant} defined in \cite{AB1}.
\vspace{2mm}

\noindent
{\bf \ref{s3}.2. Finiteness of relative cylindrical Yamabe constants.}
Let $h\in \Riem(Z)$ be a metric on the cylindrical end $Z$, and $C\in
{\cal C}(M)$ be a conformal class. First, we define the \emph{relative
$h$-cylindrical Yamabe invariant of the triple $(X,M;C)$} as
$$
Y^{h\hbox{-}\cyl}(X,M; C)\ \ \  :=\!\!\!\!\sup_{^{_{\begin{array}{c}
	\bar{g}\in \Riem^{\cyl,0}(X),
\\
\p_{\infty} \bar{g}=h, \ \p[\bar{g}] = C.
\end{array}}}} \!\!\!\!\!\!\!\! Y^{\cyl}_{[\bar{g}]}(X,M; C) .
$$
Second, we define the \emph{relative $h$-cylindrical Yamabe invariant
of the pair $(X,M)$} as
$$
Y^{h\hbox{-}\cyl}(X,M) : = \sup_{C\in {\cal C}(M)} Y^{h\hbox{-}\cyl}(X,M; C).
$$
Finally, we define the \emph{relative cylindrical Yamabe invariant of
$X$} as
$$
Y^{\cyl}(X,M) : = \sup_{h\in \Riem(Z)} Y^{h\hbox{-}\cyl}(X,M) .
$$
The definition yields the following inequalities:
$$
\begin{array}{rcl}
-\infty \leq Y^{h\hbox{-}\cyl}(X,M; C) & \leq & Y^{h\hbox{-}\cyl}(X,M) 
\\
\\
&\leq &
 Y^{\cyl}(X,M) \leq Y(S^n_+,S^{n-1}).
\end{array}
$$
The following lemma is an analogue of Lemmas \ref{cyl:L5},
\ref{cyl:L6}.
\begin{Lemma}\label{cyl:L24} $\mbox{\ }$ 
\begin{enumerate}
\item[{\bf (0)}] If $\lambda({\cal L}_h)<0$, then
$Y^{h\hbox{-}\cyl}(X,M; C)=-\infty$.
\item[{\bf (1)}] If $\lambda({\cal L}_h)\geq 0$, then
$Y^{h\hbox{-}\cyl}(X,M; C)>-\infty$.
\item[{\bf (2)}] If $\lambda({\cal L}_h)= 0$, then $0\geq
Y^{h\hbox{-}\cyl}(X,M; C)>-\infty$.
\end{enumerate}
\end{Lemma}
Now we define the following two constants which play a technical role.
$$
\begin{array}{rlc}
A &:=&  \displaystyle
\inf_{^{\begin{array}{c} _{f\in L^{1,2}(\R)} \\  ^{f\not\equiv 0}\end{array}}}
\frac{\int_{\R} \left(\alpha_n  (f^{\prime})^2 + (n-1)(n-2)f^2\right)\ dt}
{\left(\int_{\R} |f|^{\frac{2n}{n-2}} dt\right)^{\frac{n-2}{n}}},
\\
\\
A_0 &:=&  \displaystyle
\inf_{^{\begin{array}{c} _{f\in L^{1,2}(\R_{\geq 0})} 
\\  ^{f\not\equiv 0}\end{array}}}   
\frac{\int_{\R_{\geq 0}} \left(\alpha_n (f^{\prime})^2 + (n-1)(n-2)f^2\right)\ dt}
{\left(\int_{\R_{\geq 0}} |f|^{\frac{2n}{n-2}} dt\right)^{\frac{n-2}{n}}}.
\end{array}
$$
\begin{Lemma}\label{cyl:L25}
Let $A$, $A_0$ be the above constants. Then
\begin{enumerate}
\item[{\bf (1)}] $\displaystyle A = 
n(n-1)\left(\frac{\Vol(S^n(1))}{\Vol(S^{n-1}(1))}\right)^{\frac{2}{n}} =
\frac{Y(S^n)}{\Vol(S^{n-1}(1))^{\frac{2}{n}}}$.
\item[{\bf (2)}] $\displaystyle A_0 =
n(n-1)\left(\frac{\Vol(S^n_+(1))}{\Vol(S^{n-1}(1))}\right)^{\frac{2}{n}}
= \frac{Y(S^n_+, S^{n-1})}{\Vol(S^{n-1}(1))^{\frac{2}{n}}}$.
\item[{\bf (3)}] $\displaystyle A = 2^{\frac{2}{n}} A_0$.
\end{enumerate}
\end{Lemma}
{\bf Proof.}  {\bf (1)} 
From \cite{K1} and \cite{Sc3}, the Yamabe constant \ \ $
Y^{\cyl}_{[h_+ + dt^2]}(S^{n-1} \times \R)$ ($= Y(S^n) $) is attained
by the metric $f^{\frac{4}{n-2}}(h_+ + dt^2)$ with a function $f =
f(t)$ depending only on $t \in \R$.  This gives that
$\Vol(S^{n-1}(1))^{\frac{2}{n}}\cdot A = Y(S^n) 
 = n(n-1)\Vol(S^n(1))^{\frac{2}{n}}$,
and hence the equality {\bf (1)}.
\vspace{2mm}

{\bf (2), (3)}. Similarly, we have
$$
\begin{array}{rl}
\Vol(S^{n-1}(1))^{\frac{2}{n}} \cdot A_0 &=  Y(S^n_+, S^{n-1}) =
n(n-1)\Vol(S^n_+(1))^{\frac{2}{n}} 
\\
\\
&= 2^{-\frac{2}{n}}\cdot n(n-1)\Vol(S^n(1))^{\frac{2}{n}} 
\\
\\
&= 2^{-\frac{2}{n}}\cdot \Vol(S^{n-1}(1))^{\frac{2}{n}}\cdot A \ .
\ \ \ \ \ \ \ \ \Box
\end{array}
$$
\begin{Proposition}\label{cyl:L26}
Let $h\in \Riem (Z)$ be any metric and $\bar{h}=h+dt^2$ a cylindrical
metric on $Z\times \R$.
\begin{enumerate}
\item[{\bf (1)}] If $\lambda({\cal L}_h)\!<\!0$, then $ Y^{\cyl}_{[\bar{h}]}
\!(Z\!\times \!\R)\! =\! -\infty$, 
$Y^{\cyl}_{[\bar{h}]}\!(Z\!\times\!\R_{\geq 0}, Z\!\times\!\{0\};
[h])\! = \!-\infty.  $
\item[{\bf (2)}] If $\lambda({\cal L}_h)\geq 0$, then $ Y^{\cyl}_{[\bar{h}]}
(Z\times \R)\geq 0$, $Y^{\cyl}_{[\bar{h}]} (Z\times \R_{\geq 0}, Z\times\{0\};
[h]) \geq 0$.
\item[{\bf (3)}] If $\lambda({\cal L}_h)= 0$, then $ Y^{\cyl}_{[\bar{h}]}
(Z\times \R)=0$, $Y^{\cyl}_{[\bar{h}]} (Z\times \R_{\geq 0},
Z\times\{0\}; [h]) = 0$.
\item[{\bf (4)}] If $\lambda({\cal L}_h)> 0$, then $
Y^{\cyl}_{[\bar{h}]} (Z\times \R)> 0$, $Y^{\cyl}_{[\bar{h}]}
(Z\times \R_{\geq 0}, Z\times\{0\}; [h]) > 0$.
\end{enumerate}
\end{Proposition}
\begin{Proof}
The assertions (1), \ (2) \ and \ (3) \ follow from \ arguments similar to
the proof of \ Lemmas \ref{cyl:L5}, \ref{cyl:L6}. \ \ Concerning (4), \
we\ postpone \ the proof \ that $Y^{\cyl}_{[\bar{h}]} (Z\times \R)> 0$\ \ 
to \ Proposition \ref{yam:P1-new}. \ Here we only show that\ \ \ 
$Y^{\cyl}_{[\bar{h}]} (Z\times \R_{\geq 0}, Z\times\{0\}; [h]) > 0$.

Assume that $\hat{Y}= Y^{\cyl}_{[\bar{h}]} (Z\times \R_{\geq 0},
Z\times\{0\}; [h]) \leq 0$. Then there exists a sequence of
nonnegative functions ${u_i}$ with $u_i\in C^{\infty}(Z\times \R_{\geq
0})\cap L^{1,2}_{\bar{h}}(Z\times \R_{\geq 0})$ and $u_i \not\equiv 0$
such that
$$
\{
\begin{array}{l}
Q_{(Z\times \R_{\geq 0},\bar{h})}(u_i)\to \hat{Y}\leq 0, \ \ \mbox{as}
\ \ i\to\infty, \\ \displaystyle \frac{\p u_i}{\p t} = 0 \ \ \
\mbox{on} \ \ Z\times \{0\}.
\end{array}
\right.
$$
We set $\bar{u}_i\in C^{1}(Z\times \R)\cap
L^{1,2}_{\bar{h}}(Z\times \R)$ by
$$
\bar{u}_i(x,t) =
\{
\begin{array}{ll}
u_i(x,t) &\mbox{for} \ \ \ (x,t)\in Z\times \R_{\geq 0},
\\
u_i(x,-t) &\mbox{for} \ \ \ (x,t)\in Z\times \R_{\leq 0}.
\end{array}
\right.
$$
Then $\displaystyle \lim\sup_{i\to\infty} Q_{(Z\times
\R,\bar{h})}(\bar{u}_i) \leq 0$, and hence $Y^{\cyl}_{[\bar{h}]}
(Z\times \R)\leq 0$. This contradicts that $Y^{\cyl}_{[\bar{h}]}(Z \times
\R) > 0$.
\end{Proof}
Let $X$ be an open manifold with tame ends $Z \times [0, \infty)$ (and
without boundary). We decompose $X$ as
$$
X = X(1)\cup_Z (Z\times [1, \infty)).
$$
Here the manifold $Z\times [1, \infty)$ endowed with a cylindrical
metric $h + dt^2$ is considered as a cylindrical manifold with the
boundary
$$
\p(Z\times [1, \infty)) = Z\times \{1\}.
$$ 
For $\bar{g} \in \Riem^{\cyl}(X)$ with $\p_{\infty}\bar{g} = h$, set
$\bar{C} = [\bar{g}] \in {\mathcal C}^{\cyl}(X)$. Then we denote
$$
\begin{array}{rcl}
Y_1 &:=& \displaystyle
Y_{\bar{C}|_{X(1)}}(X(1), Z; [h]),
\\
\\
Y_2 &:=& \displaystyle
Y^{\cyl}_{\bar{C}|_{(Z\times [1, \infty))}}( Z\times [1, \infty), Z\times \{1\}; [h]).
\end{array}
$$
Here $Y_{\bar{C}|_{X(1)}}(X(1), Z; [h])$ is the relative Yamabe
constant defined in \cite{AB1}.
\begin{Theorem}\label{cyl:L27}
Under the above assumptions, we have
$$
Y_{\bar{C}}^{\cyl}(X)\geq 
\{
\begin{array}{cl}
-\left(|Y_1|^{\frac{n}{2}} + |Y_2|^{\frac{n}{2}}\right)^{\frac{2}{n}}&
\mbox{if} \ \ Y_1, Y_2 \leq 0,
\\
\min\{Y_1, Y_2\} & \mbox{otherwise}.
\end{array}
\right.
$$
\end{Theorem}
\begin{Proof}
First, we notice that if $Y_2=-\infty$, then there is nothing to
prove. Hence we assume that $Y_2>-\infty$.

For any function $u \in C^{\infty}_c(X)$ with
$u \not\equiv 0$,
$$
Q_{(X,\bar{g})}(u) = \frac{\int_X\[\alpha_n|du|^2 +
R_{\bar{g}}u^2\] d\sigma_{\bar{g}}}{\(\int_{X(1)}
|u|^{\frac{2n}{n-2}}d\sigma_{\bar{g}} +\int_{Z\times [1, \infty)}
|u|^{\frac{2n}{n-2}}d\sigma_{\bar{g}} \)^{\frac{n-2}{n}}} .
$$
We denote:
$$
\alpha = 
\int_{X(1)} |u|^{\frac{2n}{n-2}}d\sigma_{\bar{g}}, \ \ \ \ \ 
\beta = 
\int_{Z\times [1, \infty)} |u|^{\frac{2n}{n-2}}d\sigma_{\bar{g}}.
$$
It is enough to consider the case $\alpha, \ \beta >0$ and
$\alpha+\beta=1$. Then we have
$$
\begin{array}{rcl}
\!\!
Q_{(X,\bar{g})}(u)\!\! \!\! &=& \!\! \!\! \displaystyle 
\frac{1}{\left(1+\frac{\beta}{\alpha}\right)^{\frac{n-2}{n}}}
\frac{\int_{X(1)}\[\alpha_n |du|^2 +
R_{\bar{g}}u^2\] d\sigma_{\bar{g}}}{\alpha^{\frac{n-2}{n}}}
\\
\\
& & \displaystyle \ \ \ \ \  \  \ \ \ \ \ \  
+ \ \ 
\frac{1}{\left(1+\frac{\alpha}{\beta}\right)^{\frac{n-2}{n}}}
\frac{\int_{Z\times[1, \infty)}\[\alpha_n |du|^2 +
R_{\bar{g}}u^2\] d\sigma_{\bar{g}}}{\beta^{\frac{n-2}{n}}}
\\
\\
\!\!&= &\!\! \!\! \displaystyle 
\alpha^{\frac{n-2}{n}}\frac{\int_{X(1)}\!\[\alpha_n |du|^2\! + \!
R_{\bar{g}}u^2\] d\sigma_{\bar{g}}}{\alpha^{\frac{n-2}{n}}}
\\
\\
& & \displaystyle \ \ \ \ \  \  \ \ \ \ \ \ 
+ \ \ 
\beta^{\frac{n-2}{n}}\frac{\int_{Z\times[1, \infty)}\!\[\alpha_n |du|^2\! +\!
R_{\bar{g}}u^2\] d\sigma_{\bar{g}}}{\beta^{\frac{n-2}{n}}}
\\
\\
\!\!&\geq &\!\! \displaystyle 
\alpha^{\frac{n-2}{n}} Y_1 + (1-\alpha)^{\frac{n-2}{n}} Y_2
\end{array}
$$
for any $\alpha\in (0,1)$. Clearly one has
$$
\begin{array}{rcl}
Y_{\bar{C}}^{\cyl}(X) &\geq & \displaystyle 
\inf_{\alpha\in[0,1]}\{ \alpha^{\frac{n-2}{n}} Y_1 + (1-\alpha)^{\frac{n-2}{n}} Y_2 \}
\\
\\
&= & \displaystyle 
\{
\begin{array}{cl}
-\left(|Y_1|^{\frac{n}{2}} + |Y_2|^{\frac{n}{2}}\right)^{\frac{2}{n}}&
\mbox{if} \ \ Y_1, Y_2 \leq 0,
\\
\min\{Y_1, Y_2\} & \mbox{otherwise}.
\end{array}
\right.
\end{array}
$$
This proves Theorem \ref{cyl:L27}.
\end{Proof}
Theorem \ref{cyl:L27} and Proposition \ref{cyl:L26}(4) imply the
following result.
\begin{Corollary}\label{cyl:L28} 
Let \ $X$ \ be an open manifold with tame ends \ $Z \times [0, \infty)$ and
$\bar{g}\in \Riem^{\cyl}(X)$ any cylindrical metric with $\p_{\infty}
\bar{g} =h\in \Riem(Z)$. Assume that $\lambda({\cal L}_h)>0$ on
$Z$. Then $Y^{\cyl}_{[\bar{g}]}(X) \geq
\min\{Y_1, Y_2\} $.
\end{Corollary}
{\bf \ref{s3}.3. Surgery and the cylindrical Yamabe invariant.}  
Let $M^n$ be a closed manifold of $\dim M=n\geq 3$ and $N
\subset M^n$ be an embedded closed submanifold of $\dim N = p\leq n-1$
with trivial normal bundle $\nu_N$. We observe that the open manifold
$M\setminus N$ is a manifold with a tame end $Z\times [0,\infty)$, $Z
= S^{q-1}\times N$ since a tubular neighborhood $D^q\times N$ of $N$
without $N$ itself is diffeomorphic to
$$
S^{q-1}\times N \times (0,1) \cong S^{q-1}\times N \times \R_{> 0} 
\subset M\setminus N, \ \ \ q=n-p.
$$  
With these understood, we prove the following result.
\begin{Theorem}\label{cyl:L29} 
Let $M$ be a closed compact manifold of $\dim M =n \geq 3$ and $N$ an
embedded closed submanifold of $M$ with trivial normal bundle. Let
$g_N\in \Riem(N)$ be a given metric on $N$ and $h_+$ the standard
metric on $S^{q-1}$. Assume that $q= n-p\geq 3$.
\begin{enumerate}
\item[{\bf (1)}] If $Y(M)\leq 0$, then for any $\epsilon> 0$ there
exists $\delta= \delta(\epsilon, g_N, |Y(M)|)>0$ such that
$$
Y^{(\kappa^2\cdot h_+ + g_N)\hbox{-}\cyl}(M \setminus N) \geq Y(M)-\epsilon
$$
for any $0< \kappa\leq \delta$. In particular, $Y^{\cyl}(M \setminus
N)\geq Y(M)$.
\item[{\bf (2)}] If  $Y(M)> 0$, then there exists $\delta= \delta(g_N, Y(M))>0$ such that 
$$
Y^{(\kappa^2\cdot h_+ + g_N)\hbox{-}\cyl}(M \setminus N) >0
$$
for any $0< \kappa\leq \delta$. In particular, $Y^{\cyl}(M \setminus N)>0$.
\end{enumerate}
\end{Theorem}
\begin{Proof}
We choose a "reference" metric $g$ on $M$ and $\varepsilon > 0$. Let
$U_{\epsilon}(N)$ be an open tubular $\varepsilon$-neighborhood of
$N$. Then we define a manifold $X$ with a tame end $Z\times
[0,\infty)$ as follows. Let
$$
\begin{array}{lcl}
W &:=& M\setminus U_{\epsilon}(N),
\\
Z &:=& \p \ov{W} \cong S^{q-1}\times N,
\\
X &:=& W \cup_Z (Z\times \R_{\geq 0}) \cong M \setminus N. 
\end{array}
$$
From \cite{PY} (cf. \cite{GL1}), we recall the following:
\begin{enumerate}
\item[$\bullet$] {\sl For any Riemannian metric $g\in \Riem(M)$ and
any $\epsilon>0$, there exist $\delta= \delta(\epsilon, g_N, |Y(M)|)$,
$L=L(\epsilon, g_N, |Y(M)|)>0$ and a metric $\hat{g}\in
\Riem^{\cyl}(X)$ such that}
$$
\!\!\!
\{
\begin{array}{ll}
\mbox{{\bf (a)}} & \hat{g}= g \ \ \mbox{on} \ \ W,
\\
\mbox{{\bf (b)}} & \displaystyle
R_{\hat{g}} > R_g -\frac{\epsilon}{4} \ \ \mbox{on} \ \ X\cong M\setminus N,
\\
\mbox{{\bf (c)}} & \displaystyle
\Vol_{\hat{g}}(X(L)) \leq \Vol_{g}(M) +\frac{\epsilon}{4(1+|Y(M)|)},
\\
\mbox{{\bf (d)}} & \displaystyle
\hat{g}=(\kappa^2\cdot h_+ + g_N) + dt^2 \ \ \mbox{on} \ \ Z\times [L,\infty)
\ \ \mbox{for any} \ \ 0<\kappa\leq \delta.
\end{array} \ \ \ \ \ \  \ \ \  \ \ \ \ \ \ \  \ \ \  \
\right.
$$
{\sl In particular, $R_{\hat{g}} > 0$ on the cylinder $Z\times [L,\infty)$.}
\end{enumerate}
{\bf (1)} From the assumption $Y(M)\leq 0$, for any $\epsilon>0$ there
exists a conformal class $C\in {\cal C}(M)$ such that
$$
0\geq Y_C(M) \geq Y(M)-\frac{\epsilon}{4}.
$$
Then there exists a Yamabe metric $g\in C$ with $\Vol_g(M)=1$ and and
$R_g \equiv Y_C(M) \leq 0$. The above assertion implies that there
exists a metric $\hat{g}\in
\Riem^{\cyl}(X)$ satisfying (a)--(d). Now Proposition \ref{cyl:L26}(4)
implies that
\begin{equation}\label{cyl:eq23}
Y^{\cyl}_{[\hat{g}|_{Z\times[L,\infty)}]}(Z\times [L,\infty), Z;
[\kappa^2\cdot h_++g_N]) >0 .
\end{equation}
Theorem \ref{cyl:L27} now gives
$$
Y^{\cyl}_{[\hat{g}]}(X) \geq 
\min \{
\begin{array}{l}
Y_{[\hat{g}|_{X(L)}]}(X(L), Z; [\kappa^2\cdot h_++g_N]),
\\
Y^{\cyl}_{[\hat{g}|_{Z\times [L,\infty)}]}(Z\times [L,\infty), Z; [\kappa^2\cdot h_++g_N])
\end{array}
\} .
$$
We see that if $Y^{\cyl}_{[\hat{g}]}(X)=
Y^{\cyl}_{[\hat{g}]}(M\setminus N)\geq 0$, then
$$
Y^{(\kappa^2\cdot h_++g_N)\hbox{-}\cyl} (M\setminus N) \geq 
Y^{\cyl}_{[\hat{g}]}(M\setminus N) \geq 0\geq Y(M)\geq Y(M)-\epsilon.
$$
Furthermore, (\ref{cyl:eq23}) implies that if $Y^{\cyl}_{[\hat{g}]}(X)<0$, then
$$
\begin{array}{rcl}
0> Y^{\cyl}_{[\hat{g}]}(X) &\geq& \displaystyle
Y^{\cyl}_{[\hat{g}|_{X(L)}]}(X(L), Z;
[\kappa^2\cdot h_++g_N])
\\
\\
&\geq&
\Vol_{\hat{g}}(X(L))^{\frac{2}{n}}\cdot
\min_{X}\{R_g-\frac{\epsilon}{4}\}
\\
\\
 &\geq& \displaystyle
\left(\Vol_{g}(M) + \frac{\epsilon}{4(1+|Y(M)|)}\right)^{\frac{2}{n}}\cdot
\left(Y_C(M) -\frac{\epsilon}{4}\right)
\\
\\
&\geq& \displaystyle
Y_C(M) -\frac{3 \epsilon}{4}\geq Y_C(M) - \epsilon.
\end{array}
$$
These imply that 
$
Y^{(\kappa^2\cdot h_++g_N)\hbox{-}\cyl}(M\setminus N) \geq 
Y^{\cyl}_{[\hat{g}]}(M\setminus N) \geq Y(M)-\epsilon.$

{\bf (2)} From the assumption $Y(M) > 0$, there exists $C\in {\cal
C}(M)$ such that $Y_{C}(M) \geq \frac{1}{2}Y(M)>0$. Then there exists
a Yamabe metric $g\in C$ with $\Vol_g(M)=1$ and $R_g \equiv Y_C(M)
>0$. From the above (a), (b) and (d), there exist $\delta=\delta(g_N,
Y(M))>0$ and a metric $\hat{g}\in \Riem^{\cyl}(X)$ such that
$$
\{
\begin{array}{ll}
R_{\hat{g}}\geq \frac{1}{2}R_g > 0 &\mbox{on} \ \ X, \\ \hat{g} =
(\kappa^2\cdot h_+ +g_N)+dt^2 &\mbox{on} \ \ Z\times [L,\infty) \ \
\mbox{for} \ \ 0<\kappa\leq \delta, \\ R_{\hat{g}} =R_{\kappa^2\cdot
h_+ +g_N} \geq \delta_0>0 \ \ &\mbox{on} \ \ Z\times [L,\infty) \ \
\mbox{for some} \ \ \delta_0>0.
\end{array}
\right.
$$
Therefore Corollary \ref{cyl:L28} and Proposition \ref{cyl:L26}(4)
imply that there exists a constant $\delta_1 > 0$ such that $
Y^{(\kappa^2\cdot h_+ +g_N)\hbox{-}\cyl}(X)\geq
Y^{\cyl}_{[\hat{g}]}(X) \geq \delta_1 $ for some $\delta_1>0$.
\end{Proof}
\begin{Corollary}\label{cyl:L30}{\rm (cf. \cite[Theorem 1.1]{PY})}
Let $M_1, M_2$ be closed manifolds of $\dim M_i = n \geq 3$ and $N
 \subset M_i$ an embedded closed submanifold of $\dim N = p$ with
 trivial normal bundle {\rm ($i = 1, 2$)}.  Assume that $q = n - p
 \geq 3$. Let $M_{1, 2}$ be the manifold obtained by gluing $M_1$ and
 $M_2$ along $N$. Then
\begin{enumerate}
\item[{\bf (1)}] If $Y(M_1),Y(M_2)\leq 0$, then $Y(M_{1,2}) \geq
-\left(|Y(M_1)|^{\frac{n}{2}} +
|Y(M_2)|^{\frac{n}{2}}\right)^{\frac{2}{n}}$.
\item[{\bf (2)}] If $Y(M_1)\leq 0$ and $Y(M_2)>0$ , then $Y(M_{1,2})
\geq Y(M_1)$.
\end{enumerate} 
\end{Corollary}
\begin{Proof} Let $g_N\in \Riem(N)$ be a metric and $\varepsilon > 0$ a small constant. 

{\bf (1)} Theorem \ref{cyl:L29}(1) gives that there exists
$\kappa_{\epsilon}>0$ ($\kappa_{\epsilon}\to 0$ as $\epsilon\to 0$)
such that
$$
Y^{(\kappa^2_{\epsilon}\cdot h_+ +g_N)\hbox{-}\cyl}(M_i\setminus N)
\geq Y(M_i) -\epsilon, \ \ i=1,2.
$$
Set $U = N \times D^q$.  Then $M_{1, 2} = (M_1 \setminus U) \cup_{\p
U}(M_2 \setminus U)$.  From Theorem \ref{cyl:Th7}, we have
$$
\begin{array}{rcl}
Y(M_{1,2}) \!\!\!&= & \!\!\!\displaystyle Y((M_1 \setminus U) \cup_{\p U} (M_2
\setminus U)) 
\\ 
\\ 
&\geq& \!\!\!\displaystyle -\left(
|Y^{(\kappa^2_{\epsilon}\cdot h+g_N)\hbox{-}\cyl}(M_1 \setminus U
)|^{\frac{n}{2}} + |Y^{(\kappa^2_{\epsilon}\cdot
h+g_N)\hbox{-}\cyl}(M_2 \setminus U )|^{\!\frac{n}{2}}
\right)^{\frac{2}{n}} 
\\ 
\\ 
&\geq& \!\!\!\displaystyle -\left( (|Y(M_1)
+\epsilon |^{\frac{n}{2}} + (|Y(M_2) +\epsilon
|^{\frac{n}{2}}\right)^{\frac{2}{n}} 
\\ 
\\ 
&\geq& \!\!\displaystyle
-\left( |Y(M_1)|^{\frac{n}{2}} +
|Y(M_2)|^{\frac{n}{2}}\right)^{\frac{2}{n}} -C\epsilon.
\end{array}
$$
Hence $\displaystyle Y(M_{1,2}) \geq -\left( |Y(M_1)|^{\frac{n}{2}} +
|Y(M_2)|^{\frac{n}{2}}\right)^{\frac{2}{n}}$.
\vspace{2mm}

{\bf (2)} From Theorem \ref{cyl:L29}, there exists
$\kappa_{\varepsilon} > 0$ such that
$$
Y^{(\kappa_{\varepsilon}^2\cdot h + g_N)\hbox{-}\cyl}(M_1 \setminus N)
\geq Y(M_1) - \varepsilon, \quad Y^{(\kappa_{\varepsilon}^2\cdot h +
g_N)\hbox{-}\cyl}(M_2 \setminus N) > 0.
$$
Then Theorem \ref{cyl:Th7} implies that
$$
Y(M_{1,2}) \geq Y^{(\kappa^2_{\epsilon}\cdot h+g_N)\hbox{-}\cyl}(M_1
\setminus N ) \geq Y(M_1)-\epsilon.
$$
Hence $Y(M_{1,2})\geq Y(M_1)$.
\end{Proof}
\section{The Yamabe problem for cylindrical \\ manifolds}\label{s4}
{\bf \ref{s4}.1. Yamabe problem for cylindrical manifolds.}  Let $X$
be an open manifold of $\dim X = n \geq 3$ with tame ends $Z \times
[0, \infty)$ (and without boundary). Let $\bar{g}$ be a cylindrical
metric on $X$ with $\p_{\infty}\bar{g} = h$ for $h \in \Riem(Z)$, and
$\bar{C} = [\bar{g}] \in {\mathcal C}^{\cyl}(X)$. For simplicity,
throughout this section we assume that $Z$ is connected.  However, the
corresponding results in this section hold even when $Z$ is not
connected.

Recall the following on the Yamabe constant
$Y_{\bar{C}}^{\cyl}(X)$:
$$
\begin{array}{rcl}
Y_{\bar{C}}^{\cyl}(X) \!\!\!&=&\!\!\! \displaystyle
\inf_{u\in L^{1,2}_{\bar{g}}(X), \ u\not\equiv 0} Q_{(X,\bar{g})}(u), \ \ \ \mbox{where}
\\
\\
Q_{(X,\bar{g})}(u)\!\!\! &=&\!\!\! \displaystyle
\frac{E_{(X,\bar{g})}}{\(\int_X
|u|^{\frac{2n}{n-2}}d\sigma_{\bar{g}}\)^{\frac{n-2}{n}}} \ ,
\ \  
\displaystyle E_{(X,\bar{g})}(u) = \int_X\!\[\alpha_n |du|^2\! +\!
R_{\bar{g}}u^2\] d\sigma_{\bar{g}}.
\end{array}
$$
The following is the Yamabe problem on cylindrical manifolds.
\vspace{2mm}

\noindent
{\bf Yamabe Problem.} {\sl Given a cylindrical metric $\bar{g}$ on
$X$, does there exist a metric $\check{g}= u^{\frac{4}{n-2}} \cdot \bar{g} \in
[\bar{g}]_{L^{1,2}_{\bar{g}}}$ such that $\displaystyle
Q_{(X,\bar{g})}(u)=Y^{\cyl}_{[\bar{g}]}(X)$?}
\vspace{2mm}

\noindent
If such a function $u$ exists, we shall call $u$ a \emph{Yamabe
minimizer} with respect to the metric $\bar{g}$ and the metric
$\check{g}= u^{\frac{4}{n-2}}\cdot \bar{g}$ \ a \emph{Yamabe metric} in
the cylindrical conformal class $\bar{C}=[\bar{g}]$.
First we prove the following. 
\begin{Proposition}\label{yam:new1-1}
Let $(X,\bar{h})=(Z\times \R,h+dt^2)$ be a canonical cylindrical
manifold with $\lambda({\cal L}_h)=0 $. Then there does not exist a
Yamabe minimizer $u\in C_+^{\infty}(X)\cap L^{1,2}_{\bar{h}}(X)$ with
respect to the metric $\bar{h}$.
\end{Proposition}
\begin{Proof}
First \ recall \ that from \ Proposition \ref{cyl:L26}, if \ $\lambda({\cal
L}_h)=0 $, \ then $Y^{\cyl}_{\bar{h}}(Z\times \R)=0$.  The condition
$\lambda({\cal L}_h) = 0$ also implies that there exists $\phi\in
C_+^{\infty}(Z)$ such that ${\cal L}_{h}\phi =0$. We set
$\tilde{h}:=\phi^{\frac{4}{n-2}}\cdot \bar{h}$. Clearly we have $
R_{\tilde{h}} \equiv 0$, and $d\sigma_{\tilde{h}}
=\phi^{\frac{2n}{n-2}}d\sigma_{\bar{h}}.  $ Suppose that there exists
a Yamabe minimizer $u\in C_+^{\infty}(X)\cap L^{1,2}_{\bar{h}}(X)$
with $\int_X u^{\frac{2n}{n-2}} d\sigma_{\bar{h}} =1$, and hence 
$E_{(X,\bar{h})}(u)=0$. Then from Fact \ref{fact}(i) and Fact
\ref{app:L1} (in Appendix), $E_{(X,\tilde{h})}(\phi^{-1}u)=0$.  Set 
$v=\phi^{-1}u\in C_+^{\infty}(X)\cap L_{\tilde{h}}^{1,2}(X)$. Then
\begin{equation}\label{yam:new:eq1}
\int_X v^{\frac{2n}{n-2}} d\sigma_{\tilde{h}} < \infty \ \ \mbox{and} \ \
E_{(X,\tilde{h})}(v) =0,
\end{equation}
and hence 
$$
0 = E_{(X,\tilde{h})}(v) =\int_X\left(\alpha_n |dv|^2 + R_{\tilde{h}}
v^2\right)d\sigma_{\tilde{h}} =\alpha_n \int_X |dv|^2d
\sigma_{\tilde{h}}.
$$
This combined with (\ref{yam:new:eq1}) gives that $v \equiv 0$. 
This contradicts that $v > 0$.
\end{Proof}
\vspace{2mm}

\noindent
{\bf \ref{s4}.2. Solution of the Yamabe problem.}  \ \ 
Recall that \ $Y^{\cyl}_{\bar{C}}(X) \leq
 Y^{\cyl}_{[h + dt^2]}(Z \times \R)$.  First, we consider the case when
 $\lambda(L_h) > 0$ and $Y^{\cyl}_{\bar{C}}(X) < Y^{\cyl}_{[h + dt^2]}(Z
 \times \R)$.
\begin{Theorem}\label{yam:Th1}
Let \ $X$ \ be an open manifold with a connected
tame end $Z\times [0,\infty)$ and $h\in \Riem(Z)$ a metric with
$\lambda({\cal L}_h)>0$. Let $\bar{g}$ be a cylindrical metric on $X$
with $\p_{\infty}\bar{g} = h$, and $\bar{C} = [\bar{g}]$.  Assume that
$$
Y_{\bar{C}}^{\cyl}(X)< Y^{\cyl}_{[h+dt^2]}(Z\times \R).
$$
Then there exists a \ Yamabe minimizer \ $u\in C_+^{\infty}(X)\cap
L^{1,2}_{\bar{g}}(X)$ \ with \ \ $\int_X u^{\frac{2n}{n-2}}
d\sigma_{\bar{g}}=1$ such that $Q_{(X,\bar{g})}(u) =
Y_{\bar{C}}^{\cyl}(X)$. In particular, the minimizer $u$ satisfies the
Yamabe equation:
$$
\L_{\bar{g}} u = -\alpha_n \Delta_{\bar{g}} u + R_{\bar{g}} u = 
Y_{\bar{C}}^{\cyl}(X)u^{\frac{n+2}{n-2}}.
$$
\end{Theorem}
For the proof of Theorem \ref{yam:Th1}, we first prove 
several lemmas and propositions.

Recall that $X(L)= X\setminus (Z\times (L,\infty))$ for each
$L>0$. From the condition $Y^{\cyl}_{\bar{C}}(X) <
Y^{\cyl}_{[h+dt^2]}(Z \times \R)$, there exists $L_0 >> 1$ such that
$$
\begin{array}{l} 
Q_L : = {\displaystyle
\inf_{^{\begin{array}{c} _{f \in L^{1, 2}_{\bar{g}}(X), \ \ 
f \not\equiv 0} \\
^{f \equiv 0 \ {\small \textrm{on}} \ Z \times [L, \infty)}\end{array}}}}
 Q_{(X, \bar{g})}(f) < Y^{\cyl}_{[h+dt^2]}(Z \times \R) 
 \quad \textrm{for any} \quad  L \geq L_0. 
\end{array}
$$
Then the standard argument combined with the inequality
$Y^{\cyl}_{[h+dt^2]}(Z \times \R) \leq Y(S^n)$ implies the following
fact (cf. \cite[Chapter 5]{SY4}).
\begin{Lemma}\label{yam:L1}
Let $(X, \bar{g})$ be a cylindrical manifold as above. Then there
exists a nonnegative function $u_L\in C^0(X)\cap
C^{\infty}(\intt(X(L)))$ for each $L \geq L_0$ such that
$$
\{
\begin{array}{l}
\displaystyle
Q_{(X,\bar{g})}(u_L) = Q_L \ ,  
\\
\\
\displaystyle
u_L > 0 \ \ \mbox{on} \ \ \intt(X(L)), \ \ \ u_L \equiv 0 \ 
 \ \mbox{on} \ \ Z\times [L,\infty),
\\
\\
\displaystyle
\int_X u^{\frac{2n}{n-2}}_L d\sigma_{\bar{g}}=1.
\end{array}
\right.
$$
\end{Lemma}
In particular, the function $u_L$ satisfies the equation
$$
-\alpha_n \Delta_{\bar{g}}u_L+ R_{\bar{g}}u_L = Q_L u_L^{\frac{n+2}{n-2}}
\ \ \ \mbox{on} \ \ \intt(X(L)). 
$$
Moreover, the following properties of the constant $Q_L$ hold:
\begin{enumerate}
\item[$\bullet$] $Q_{L_1}\geq Q_{L_2}$ if $L_1\leq L_2$. 
\item[$\bullet$] $Q_L \to Y_{\bar{C}}^{\cyl}(X)$ as $L\to \infty$. 
\end{enumerate}
\begin{Lemma}\label{yam:L2}
There exists a positive constant $C_0$ such that $u_L \leq C_0$ on $X(2)$ for
any $L\geq L_0$.
\end{Lemma}
\begin{Proof}
 Suppose that there exist $L_i>0$ and a point $p_i\in X(2)$ for
$i=1,2,\ldots$ such that
$$
L_i\to \infty, \ \ \mbox{and} \ \ u_{L_i}(p_i)=\max_{X(2)} u_{L_i} =:
m_i\to \infty.
$$
Since $X(2)$ is compact, we may assume that $p_i\to p_0\in X(2)$. 

Let $\{U, x=(x^1,\ldots,x^n)\}$ be a normal coordinate system centered
at $p_0$. We may also assume that $\{|x|<1 \} \subset U$. We define
the functions
$$
v_i(x):= m_i^{-1}\cdot u_{L_i}\left( m_i^{-\frac{2}{n-2}}\!\!\cdot x \!+\!
x(p_i) \right) \ \ \mbox{for \  $x\in \{ |x|< m_i^{\frac{2}{n-2}}(1-|x(p_i)|) \}$.}
$$
Similarly to the proof of Theorem 2.1 \cite[Chapter 5]{SY4}, there
exists a positive function $v\in C^{\infty}_+(\R^n)$ such that $v_i$
converges to $v$ in the $C^2$-topology on each relatively compact
domain in $\R^n$. Let $\Delta_0$ be the Laplacian on $\R^n$ with
respect to the Euclidean metric.  Then the function  $v$
satisfies  the following:
$$
\{
\begin{array}{l}
-\alpha_n \Delta_0 v = Y^{\cyl}_{\bar{C}}(X)\cdot v^{\frac{n+2}{n-2}}
\ \ \mbox{on} \ \ \R^n, \\ \displaystyle
\int_{\R^n} v^{\frac{2n}{n-2}} dx \leq
\lim\inf_{i\to\infty} \int_{X(2)} u_{L_i}^{\frac{2n}{n-2}} \leq 1.
\end{array}
\right.
$$
Hence $Y_{\bar{C}}^{\cyl}(X)\!\geq\!Y(S^n)$. This contradicts
that $Y^{\cyl}_{\bar{C}}(X) < Y^{\cyl}_{[h+dt^2]}(Z \times \R) \leq
Y(S^n)$.
\end{Proof}
We consider $u_L$ on the end $Z\times [1,\infty)$ with the metric
$\bar{g}=h+dt^2$. Note that the metric $\bar{g}$ is invariant under
parallel translations along the $t$-coordinate. By using this fact, one
can obtain (similar to the proof of Lemma \ref{yam:L2}) that
\begin{equation}\label{yam:new1}
u_L \leq K \ \ \ \mbox{on} \ \ \ Z\times [2,\infty)
\end{equation}
for any $L\geq L_0$, where $K>0$ is a constant independent of $L$. Set
$K_0:=\max\{C_0,K\}>0$. Then Lemma \ref{yam:L2} and (\ref{yam:new1})
imply the following.
\begin{Lemma}\label{yam:L3} 
There exists a constant $K_0>0$ such that $u_L \leq K_0$ on $X$ for
any $L\geq L_0$.
\end{Lemma}
\begin{Convention}\label{convention}
{\rm 
Let \ $(X,\bar{g})$ \ be a cylindrical manifold with \ $\bar{g}= h+
dt^2$ \ on $Z\times [1,\infty)$, where $h$ is a metric on $Z$ with
$\lambda({\mathcal L}_h) > 0$ (resp. $\lambda({\mathcal L}_h) = 0$).
Throughout this and the next sections, we will use the following
convention.  

By the reason below (replacing the metric $\bar{g}$ by a
suitable pointwise conformal metric if necessary), we may assume that
the cylindrical metric $\bar{g}$ satisfies the following
$$
R_{\bar{g}} = R_h \geq R_{\min} :=\min_Z R_h >0 \ \
\mbox{(resp. $R_{\bar{g}} = R_h =0$) \ on} \ \ Z\times [1,\infty).
$$
With this understood, we use some constants depending on $R_{\min}$ in
this and the next sections.  One can easily replace those constants by
other ones which depend on the invariant $\lambda({\cal L}_h)$ and
$\delta_h$ defined below. 

The condition $\lambda({\cal L}_h)>0$ (resp. $\lambda({\cal L}_h)=0$)
gives that there exists a function $\phi_h\in C_+^{\infty}(Z)$ such
that
$$
\{
\begin{array}{l}
{\cal L}_h \phi_h = \lambda({\cal L}_h)\cdot \phi_h \ \ \mbox{on} \ \ Z,
\\
\displaystyle
\max_{Z}\phi_h =1.
\end{array}
\right.
$$
Let $\varphi \in C^{\infty}_+(X)$ be a positive smooth function with
$\varphi \equiv 1$ on $X(0)$ and $\varphi \equiv \varphi_h$ on $Z
\times [1, \infty)$.  Set $\tilde{g} = \varphi^{\frac{4}{n-2}}\cdot
\bar{g}$ and $\delta_h : = \min_Z \varphi > 0$. Then we have
$$
\{
\begin{array}{l}
\delta_h^{\frac{4}{n-2}}\cdot \bar{g} \leq \tilde{g}\leq\bar{g}
\\
\delta_h^{\frac{2n}{n-2}}\cdot d\sigma_{\bar{g}} \leq  d\sigma_{\tilde{g}} \leq 
d\sigma_{\bar{g}} 
\end{array}
\right. \ \ \mbox{on} \ \ Z\times [1,\infty).
$$
Clearly $\Delta_{\tilde{g}} f = \phi^{-\frac{4}{n-2}}_h\cdot
f^{\prime\prime}$ for any function $f = f(t)\in C^{\infty}(Z\times
[1,\infty)) $ which depends only on $t\in [1,\infty)$. We can use this
property to construct comparison functions on $Z\times
[1,\infty)$. Now we have
$$
R_{\tilde{g}} = \lambda({\cal L}_h)\cdot \phi^{-\frac{4}{n-2}}\geq
\lambda({\cal L}_h)>0 \ \ \mbox{(resp. $R_{\tilde{g}}= \lambda({\cal
L}_h )\cdot \phi^{-\frac{4}{n-2}}\equiv 0$)}
$$
on \ $Z\times [1,\infty)$. \ We emphasize that the metric \ \
$\tilde{g}=\phi^{\frac{4}{n-2}}\cdot\bar{g}$ \ \ is no longer a product
metric on \ $Z\times [1,\infty)$, \ but $\phi=\phi_h$ is a positive smooth
function depending only on $x\in Z$. One can use this property of the
metric $\tilde{g}=\phi^{\frac{4}{n-2}}\cdot\bar{g}$ to obtain higher
derivative estimates on Yamabe minimizers similarly to the cylindrical
metric case. \hfill $\Box$}
\end{Convention}
Consider the canonical cylindrical manifold $(X, \bar{h}) = (Z \times
 R, h + dt^2)$ associated to $(X, \bar{g})$.
\begin{Proposition}\label{yam:P1-new}
Let $(X, \bar{h}) = (Z \times \R, h + dt^2)$ be a canonical
cylindrical manifold with $\lambda({\cal L}_h)>0$. Then
$Y^{\cyl}_{[\bar{h}]}(X)>0$.
\end{Proposition}
\begin{Proof}
As above, we may assume that $\displaystyle R_{\min} =\min_Z R_h >0$.
Then $\Ric_h \geq (n-1)\kappa$ on $Z$ since $Z$ is compact, where
$\kappa$ is a constant (not necessarily positive).  Denote by $B_r(x)$
a geodesic ball of radius $r$ centered at $x\in X$.  Let
$\Delta_{\bar{h}} -\p/\p \tau$ be the heat operator on $X$ and
$p=p(x,y,\tau)$ the heat kernel of $\Delta_{\bar{h}} -\p/\p
\tau$. Then the following estimate on $p(x, x, \tau)$ was proved by
Li-Yau (see \cite{LY}):
\begin{equation}\label{yam:new2}
p(x,x,\tau) \leq
\frac{c(n,\delta)}
{\Vol_{\bar{h}}(B_{\sqrt{\tau}}(x))}\exp(-c(n)\delta\kappa\tau)
\end{equation}
for all $\delta>0$, $\tau>0$ and $x\in X$. Here $c(n,\delta)$
(resp. $c(n)$) is a positive constant depending only on $n$ and
$\delta>0$ (resp. $n$). We also notice that
\begin{equation}\label{yam:new3}
\Vol_{\bar{h}}(B_{\sqrt{\tau}}(x)) \geq C^{-1}\cdot \tau^{\frac{n}{2}}
\end{equation}
for $0< \tau \leq \sqrt{\frac{1}{2}\diam(Z,h)}$, where $C>0$ is a
constant.  Set
$$
\tau_0:=\max\{\frac{\alpha_n}{R_{\min}}, \frac{1}{2}\diam(Z,h)\}
>0
$$ 
and $\delta=1$. Using (\ref{yam:new3}) in (\ref{yam:new2}), we then
obtain the following estimate
\begin{equation}\label{yam:new4}
p(x,x,\tau) \leq \frac{C\cdot
c(n,\delta)}{\tau^{\frac{n}{2}}}\exp(c(n)|\kappa|\tau_0)
\end{equation} 
for $0<\tau\leq\tau_0$. Set $ C_4 := C\cdot
c(n,\delta)\exp(c(n)|\kappa|\tau_0)$. Now we use \cite[Theorem
2.2]{Sa} and (\ref{yam:new4}) to obtain
$$
\begin{array}{rcl}
\displaystyle
\left(\int_X |f|^{\frac{2n}{n-2}} d\sigma_{\bar{h}}\right)^{\frac{n-2}{n}}
&\leq & \displaystyle
C^{\prime}\cdot C_4^{\frac{2}{n}}\alpha_n
\int_X \left(
|df|^2 +\frac{1}{\tau_0} f^2
\right)d\sigma_{\bar{h}}
\\
\\
&\leq & \displaystyle
C^{\prime}\cdot C_4^{\frac{2}{n}}\alpha_n
\int_X \left(
|df|^2 +\frac{R_h}{\alpha_n} f^2
\right)d\sigma_{\bar{h}}
\end{array}
$$
for any $f\in C_c^{\infty}(X)$, where $C^{\prime}>0$ is a constant
independent of $f$. This implies that
$$
Y^{\cyl}_{[\bar{h}]}(X) = Y^{\cyl}_{[h+dt^2]}(Z\times \R) \geq 
\frac{1}{C^{\prime}\cdot C_4^{\frac{2}{n}}} > 0.
$$
This completes the proof of Proposition \ref{yam:P1-new}.
\end{Proof}
Now we return to consider the cylindrical manifold $(X, \bar{g})$ in
Theorem \ref{yam:Th1}.  If $Y^{\cyl}_{\bar{C}}(X) > 0$, we set
$$
Y:= Y_{[h+dt^2]}^{\cyl}(Z\times \R)>0, \ \ \ \delta:=
\frac{Y^{\cyl}_{\bar{C}}(X)}{Y}>0 .
$$
Clearly $\delta <1$ by the assumption in Theorem \ref{yam:Th1}. In this
case, we may assume that 
$$
0< \frac{Q_L}{Y} \leq \frac{3\delta+1}{4} < 1
$$
for any $L\geq L_0$. If $Y^{\cyl}_{\bar{C}}(X) \leq 0$, we set $\delta
= 0$.  Using the Moser iteration technique, we then show the following
decay estimate of $u_L$ on the cylindrical end $Z \times [1,
\infty)$.
\begin{Proposition}\label{yam:P1}
Let \ $u_L$ be the minimizer obtained in \ Lemma \ref{yam:L1} \ for each $L
 \geq L_0$. Then there exists
$$
L_1=L_1\left(\frac{1}{(1-\delta)^2 R_{\min}}, \ n \right)\geq 2
$$ 
such that
\begin{equation}\label{yam:P1-eq1}
\sup_{Z\times [t_0-\frac{1}{2}r,t_0+\frac{1}{2}r]} u_L \leq \frac{C_n}{Y^{\frac{n-2}{4}}}
\cdot\frac{1}{r^{\frac{n-2}{2}}}
\end{equation}
for any $t_0$, $r>0$ satisfying $t_0-r\geq L_1$.
Here $C_n>0$ is a constant depending only on $n$. 
\end{Proposition}
\begin{Proof} We give the proof  only for the case 
\ $Y^{\cyl}_{\bar{C}}(X)>0$. \ The case $Y^{\cyl}_{\bar{C}}(X)\leq 0$ is
much easier, hence we omit it. For any $f\in C^{\infty}_{c}(X)$ with
$\supp(f)\subset Z\times [1,\infty)$, we have
\begin{equation}\label{eq-new1}
\int_X |f|^{\frac{2n}{n-2}} d\sigma_{\bar{g}} 
\leq \frac{\alpha_n}{Y}\int |df|^2 d\sigma_{\bar{g}} + 
\frac{1}{Y} \int_X R_{h} f^2 d\sigma_{\bar{g}}.
\end{equation}
The minimizer $u_L$ satisfies the equation
\begin{equation}\label{yam:P1-eq2}
-\Delta_{\bar{g}} u_L + \frac{R_{\bar{g}}}{\alpha_n} u_L = 
\frac{Q_L}{\alpha_n} u_L^{\frac{n+2}{n-2}}
\end{equation}
with $u_L\equiv 0$ on $X\setminus X(L)$ and $u_L>0$ on
$\intt(X(L))$. Then $u_L$ satisfies the following differential
inequality
\begin{equation}\label{eq-new2}
-\Delta_{\bar{g}} u_L + \frac{R_{\bar{g}}}{\alpha_n} u_L \leq 
\frac{Q_L}{\alpha_n} u_L^{\frac{n+2}{n-2}} \ \ \ \mbox{on} \ \ \ X
\end{equation}
in the distributional sense.

Let $\eta $ be a cut-off function with $\supp(\eta)\subset Z\times
[1,\infty)$ and $0\leq \eta\leq 1$.  Recall that $\bar{g}=h+dt^2$ and
$R_{\bar{g}}= R_h$ on the cylindrical end $Z\times [1,\infty)$. We use
(\ref{eq-new1}) and (\ref{eq-new2}) to prove the following assertion.
\begin{Claim}\label{Claim1}
For any $\alpha\geq 1$ and any small $\epsilon>0$, the
following estimate holds:
\begin{equation}\label{new-eq3}
\!\!\!\!\!\!
\begin{array}{rcl}
\displaystyle
\left(
\int_X \left(
\eta^2 u_L^{\alpha+1}
\right)^{\frac{n}{n-2}}
d\sigma_{\bar{g}}
\right)^{\frac{n-2}{n}} \!\!\!& \leq	&\!\!\!\displaystyle
\frac{(\alpha+1)(\alpha+1+\epsilon)}{4(\alpha-\frac{\epsilon}{2})}\cdot
\frac{Q_L}{Y} 
\int_X  \eta^2 u_L^{\alpha+\frac{n+2}{n-2}} d\sigma_{\bar{g}}
\\
\\
&& \displaystyle \ \ \ \ \ \ \  + \ \frac{\bar{C}_n\alpha\epsilon^{-1}}{Y}
\int_X |d\eta|^2 u_L^{\alpha+1} d\sigma_{\bar{g}},
\end{array}\!\!\!\!\!\!\!\!\!\!\!\!
\end{equation}
where $\bar{C}_n>0$ is a constant depending only on $n$.
\end{Claim}
\begin{lproof}{Proof of Claim \ref{Claim1}}
We multiply both sides of (\ref{eq-new2}) by $\eta^2 u^{\alpha}$ for
any $\alpha>0$. Then integrating by parts, we obtain
\begin{equation}\label{new-eq4}
\begin{array}{rcl}
& & \displaystyle
\frac{4\alpha}{(\alpha+1)^2}
\int_X \eta^2 \left| d u_L^{\frac{\alpha+1}{2}}\right|^2 d\sigma_{\bar{g}}
-2 \int_X \eta |d\eta| u_L^{\alpha} |du_L| d\sigma_{\bar{g}}
\\
\\
&\leq & \displaystyle
-\frac{1}{\alpha_n} \int_X R_h\eta^2 u_L^{\alpha+1} d\sigma_{\bar{g}} +
\frac{Q_L}{\alpha_n} \int_X \eta^2 u_L^{\alpha+\frac{n+2}{n-2}} d\sigma_{\bar{g}}.
\end{array}
\end{equation}
The Young inequality implies
\begin{equation}\label{new:eq7}
\eta |d\eta| u_L^{\alpha} |du_L| \leq \epsilon^{-1} \cdot |d\eta|^2 u_L^{\alpha+1} +
\epsilon\cdot \frac{\eta^2}{(\alpha+1)^2} \left| du_L^{\frac{\alpha+1}{2}} \right|^2
\end{equation}
for any $\epsilon>0$. We use (\ref{new-eq4}) and (\ref{new:eq7}) to
obtain the estimate
\begin{equation}\label{new:eq8}
\!\!\!\!
\begin{array}{l}
\displaystyle
\int_X \eta^2 \left| d u_L^{\frac{\alpha+1}{2}}\right|^2 d\sigma_{\bar{g}}
\\
\\
\displaystyle
\leq \frac{(\alpha+1)^2}{4(\alpha-\frac{\epsilon}{2})}\left[
\int_X \!\left(\!
2\epsilon^{-1} \!|d\eta|^2 \!- \!\frac{1}{\alpha_n}R_h\eta^2
\!\right) u_L^{\alpha+1} d\sigma_{\bar{g}} +
\frac{Q_L}{\alpha_n}\int_X \!\eta^2 \!u_L^{\alpha+\frac{n+2}{n-2}} d\sigma_{\bar{g}}
\right]
\end{array}\!\!\!\!\!\!\!\!\!\!\!\!\!\!\!\!
\end{equation}
for $0<\epsilon <2\alpha$. It then follows from (\ref{eq-new1}),
(\ref{new:eq7}) and (\ref{new:eq8}) that 
\begin{equation}\label{new:eq9}
\!\!\!\!
\begin{array}{l}
\displaystyle \left( \int_X \left( \eta u_L^{\frac{\alpha+1}{2}}
\right)^{\frac{2n}{n-2}} d\sigma_{\bar{g}} \right)^{\frac{n-2}{n}} 
\\
\\ 
\ \ \ \ \leq \displaystyle 
\frac{\alpha_n}{Y} \int_X
\left|d\left(\eta u_L^{\frac{\alpha+1}{2}}\right) \right|^2
d\sigma_{\bar{g}} + 
\frac{1}{Y}\int_X R_h\left(\eta u_L^{\frac{\alpha+1}{2}}\right)^2 d\sigma_{\bar{g}}
\\
\\ 
\ \ \ \ \leq  \displaystyle 
\frac{(\alpha+1)(\alpha+1+\epsilon)}{4(\alpha-\frac{\epsilon}{2})}
\cdot\frac{Q_L}{Y} \int_X \eta^2 u_L^{\alpha+\frac{n+2}{n-2}}  d\sigma_{\bar{g}}
\\
\\
\ \ \ \ \ \  
\displaystyle
+ \ \frac{\alpha_n}{Y} 
\{ 
\frac{(\alpha+1)(\alpha+1+\epsilon)\epsilon^{-1}}{2(\alpha-\frac{\epsilon}{2})}
+ 1 + (\alpha+1)\epsilon^{-1}
\}
\int_X |d\eta|^2 u_L^{\alpha+1}  d\sigma_{\bar{g}}
\\
\\
\ \ \ \ \ \ 
\displaystyle
- \ \frac{1}{Y} 
\{
\frac{(\alpha+1)(\alpha+1+\epsilon)}{4(\alpha-\frac{\epsilon}{2})} -1 
\}
\int_X R_h \eta^2 u_L^{\alpha+1} d\sigma_{\bar{g}}. 
\end{array}\!\!\!\!\!\!\!\!\!\!\!\!\!\!\!\!\!\!
\end{equation}
We notice that
$$
\frac{(\alpha+1)(\alpha+1+\epsilon)}{4(\alpha-\frac{\epsilon}{2})} -1 > 0
$$
for $0<\epsilon < 2\alpha$. From $R_h > 0$, 
$$
\int_X R_h\eta^2 u_L^{\alpha+1} d\sigma_{\bar{g}} > 0 .
$$
Combining these observations with (\ref{new:eq9}), we obtain the
inequality (\ref{new-eq3}) for $\alpha\geq 1$.
\end{lproof}
Now we continue with the proof of Proposition \ref{yam:P1}. Set
$\alpha = 1$ in (\ref{new-eq3}), then
\begin{equation}\label{new:eq10}
\begin{array}{rl}
\displaystyle
\!\!\!\!\!\!\!
\left(
\int_X\! \eta^{\frac{2n}{n-2}} u_L^{\frac{2n}{n-2}}  d\sigma_{\bar{g}}
\right)^{\!\frac{n-2}{n}}\!\!\!\!
&\leq \displaystyle 
\frac{2+\epsilon}{2-\epsilon} \cdot \frac{Q_L}{Y} 
\int_X \eta^2 u_L^{\frac{2n}{n-2}} d\sigma_{\bar{g}}
\\
\\
&\displaystyle
\ \ \ \ \ \ \ \ \ \ +\ \frac{\bar{C}_n\epsilon^{-1}}{Y}
\int_X |d\eta^2|u_L^2 d\sigma_{\bar{g}}.
\end{array}
\end{equation}
Then we use
$$
\int_X \eta^{\frac{2n}{n-2}} u_L^{\frac{2n}{n-2}}  d\sigma_{\bar{g}} \leq
\int_X u_L^{\frac{2n}{n-2}}  d\sigma_{\bar{g}} = 1
$$
to obtain the estimate
\begin{equation}\label{new:eq11}
\int_X \eta^{\frac{2n}{n-2}} u_L^{\frac{2n}{n-2}}  d\sigma_{\bar{g}} \leq
\left(
\int_X \eta^{\frac{2n}{n-2}} u_L^{\frac{2n}{n-2}}  d\sigma_{\bar{g}}
\right)^{\frac{n-2}{n}}.
\end{equation}
From  (\ref{new:eq10}) and  (\ref{new:eq11}), we also obtain
\begin{equation}\label{new:eq12}
\left(
1- \frac{2+\epsilon}{2-\epsilon} \cdot \frac{Q_L}{Y} 
\right)
\int_X \eta^{\frac{2n}{n-2}} u_L^{\frac{2n}{n-2}}  d\sigma_{\bar{g}}
\leq 
\frac{\bar{C}_n\epsilon^{-1}}{Y} 
\int_X |d\eta^2|u_L^2 d\sigma_{\bar{g}}.
\end{equation}
Take $\epsilon>0$ in (\ref{new:eq12}) as
$\epsilon=\frac{2(1-\delta)}{3+5\delta}$. From
$$
0 < \frac{Q_L}{Y} \leq \frac{3\delta +1}{4} < 1,
$$
we notice that
\begin{equation}\label{new:eq13}
1- \frac{2+\epsilon}{2-\epsilon} \cdot \frac{Q_L}{Y} \geq \frac{1-\delta}{2}.
\end{equation}
It then follows from (\ref{new:eq12}),  (\ref{new:eq13}) that
\begin{equation}\label{new:eq14}
\int_X \eta^{\frac{2n}{n-2}} u_L^{\frac{2n}{n-2}} d\sigma_{\bar{g}}
\leq \frac{C_n^{\prime}}{(1-\delta)^2Y} \int_X |d\eta|^2 u_L^2
d\sigma_{\bar{g}},
\end{equation}
where $C_n^{\prime}>0$ is a constant depending only on $n$.

With these understood, we also show the following.
\begin{Claim}\label{Claim2}
\begin{equation}\label{new:eq15}
\int_{Z\times [r+1,\infty)} u^{\frac{2n}{n-2}}_L d\sigma_{\bar{g}} \leq
\frac{4\cdot C_n^{\prime}}{(1-\delta)^2 R_{\min}}\cdot \frac{1}{r^2}
\end{equation}
for any $r\geq 1$.
\end{Claim}
\begin{lproof}{Proof of Claim \ref{Claim2}}
For any $T>0$, choose a cut-off function $\eta$ in (\ref{new:eq14})
satisfying 
$$
\{
\begin{array}{ll}
\eta=1 & \mbox{on $ Z\times [r+1,T+r+1]$},
\\
\eta=0 & \mbox{on $X\setminus (Z\times [1,T+2r+1])$},
\\
|d\eta|\leq \frac{2}{r}& \mbox{on} \ \ Z\times\left( [1,1+r]\sqcup [T+r+1,T+2r+1]\right).
\end{array}
\right.
$$
From (\ref{eq-new2}) and $\int_X u_L^{\frac{2n}{n-2}}
d\sigma_{\bar{g}}=1$, we notice that
$$
\begin{array}{rcl}
\displaystyle \int_X u_L^2 d\sigma_{\bar{g}} &\leq&\displaystyle
\frac{1}{R_{\min}}\int_X\left( \alpha_n|du_L|^2 + R_h u_L^2 \right)
d\sigma_{\bar{g}} \\ \\ &\leq&\displaystyle \frac{Q_L}{R_{\min}}\int_X
u_L^{\frac{2n}{n-2}} d\sigma_{\bar{g}} \leq \frac{Q_L}{R_{\min}}.
\end{array}
$$
This combined with (\ref{new:eq14}) implies
$$
\int_{Z\times[r+1,T+r+1]} u_L^{\frac{2n}{n-2}} d\sigma_{\bar{g}} \leq
\frac{4C^{\prime}_n}{(1-\delta)^2R_{\min}}\cdot\frac{Q_L}{Y}\cdot\frac{1}{r^2},
$$
and hence from $\frac{Q_L}{Y}\leq 1$,
$$
\int_{Z\times[r+1,T+r+1]} u_L^{\frac{2n}{n-2}} d\sigma_{\bar{g}} 
\leq \frac{4C^{\prime}_n}{(1-\delta)^2 R_{\min}}\cdot\frac{1}{r^2}.
$$
Letting $T \rightarrow \infty$, we then obtain the estimate
 (\ref{new:eq15}).
\end{lproof}
We return to the proof of Proposition \ref{yam:P1}. From H\"{o}lder's
inequality, we notice
$$
\begin{array}{rcl}
\displaystyle
\int_X \eta^2 u_L^{\frac{2(n+2)}{n-2}} d\sigma_{\bar{g}} \!\!\!&=&\!\!\!
\displaystyle
\int_X u_L^{\frac{4}{n-2}}\cdot\eta^2 u_L^{\frac{2n}{n-2}} d\sigma_{\bar{g}}
\\
\\
&\leq &\!\!\!
\displaystyle
\left(
\int_{\supp(\eta)} u_L^{\frac{2n}{n-2}} d\sigma_{\bar{g}}
\right)^{\!\frac{2}{n}}\!\!\cdot\!
\left(
\int_X \left(
\eta^2 u_L^{\frac{2n}{n-2}} 
\right)^{\!\frac{n}{n-2}} d\sigma_{\bar{g}}
\right)^{\!\frac{n-2}{n}}.
\end{array}
$$
Set $\alpha=\frac{n+2}{n-2}$ and $\epsilon=1$ in (\ref{new-eq3}). Then
this combined with the above inequality implies
\begin{equation}\label{new:eq16}
\!\!\!\!\!\!
\begin{array}{l}
\displaystyle
\{
1-\frac{\frac{2n}{n-2}\left(\frac{2n}{n-2}+1\right)}{4\left(\
\frac{n+2}{n-2} -\frac{1}{2}\right)}
\left(
\int_{\supp(\eta)}\!\!\!u_L^{\frac{2n}{n-2}} d\sigma_{\bar{g}}
\right)^{\!\frac{2}{n}}
\}\cdot
\left(
\int_X\!\!
\left(\eta^2 u_L^{\frac{2n}{n-2}}  \right)^{\frac{n}{n-2}}d\sigma_{\bar{g}}
\right)^{\!\frac{n-2}{n}}
\\
\\
\ \ \ \ \ \ \ \ \ \ \  \ \ \  \ \ \ \ \ \ \ \ \ \ \  \ \ \ \ \ \ \ \ \ \ \ \ \ \ \ \  \ \ \ \ \ \ \ 
\displaystyle
\leq \ \frac{\bar{C_n}}{Y}\cdot\frac{n+2}{n-2}
\int_X|d\eta|^2 u^{\frac{2n}{n-2}} d\sigma_{\bar{g}}.
\end{array}\!\!\!\!\!\!\!\!\!\!\!\!\!\!\!\!\!\!
\end{equation}
From (\ref{new:eq15}), there exists $ L_1 =L_1
\left(\frac{1}{(1-\delta)^2 R_{\min}}, \ n \right) \geq 2$ such that
\begin{equation}\label{new:eq17}
\frac{\frac{2n}{n-2}\left(\frac{2n}{n-2}+1\right)}{4\left(\
\frac{n+2}{n-2} -\frac{1}{2}\right)}\left(
\int_{\supp(\eta)}u_L^{\frac{2n}{n-2}} d\sigma_{\bar{g}}
\right)^{\frac{2}{n}} \leq \frac{1}{2}
\end{equation}
for any $\eta$ satisfying $\supp(\eta)\subset Z\times [L_1,\infty)$. It then follows from
(\ref{new:eq16}) and (\ref{new:eq17}) that
$$
\left( \int_X\left(\eta^2 u_L^{\frac{2n}{n-2}}
\right)^{\frac{n}{n-2}}d\sigma_{\bar{g}} \right)^{\frac{n-2}{n}} \leq
\frac{10\cdot \bar{C}_n}{Y}\int_X |d\eta|^2 u_L^{\frac{2n}{n-2}}
d\sigma_{\bar{g}}.
$$
$t_0-r\geq L_1$. Then we choose a cut-off function $\eta$ satisfying
the following conditions:
$$
\{
\begin{array}{cl}
\eta=1 & \mbox{on} \ \ Z\times [t_0-\frac{3}{4}r, t_0+ \frac{3}{4}r],
\\
\eta=0 & \mbox{on} \ \ X\setminus \left(Z\times [t_0-r, t_0+ r]\right),
\\
|d\eta|\leq \frac{5}{r}& \mbox{on} \ \ 
Z\times \left([t_0-r, t_0-\frac{3}{4}r]\sqcup [t_0+ \frac{3}{4}r,t_0+ r]\right).
\end{array}
\right.
$$
Let $\chi:=\chi_{ Z\times [t_0-\frac{3}{4}r, t_0+ \frac{3}{4}r]}$ denote the
characteristic function on the cylinder $Z\times [t_0-\frac{3}{4}r, t_0+
\frac{3}{4}r]$. Then we have
\begin{equation}\label{new:eq19}
\|\chi\cdot u_L \|_{L_{\bar{g}}^{\frac{2n^2}{(n-2)^2}}} \leq 
\left(
\frac{250\cdot \bar{C}_n}{Y}\cdot\frac{1}{r^2}
\right)^{\frac{n-2}{2n}}.
\end{equation}
Now we specify the cut-off function as follows. Let $\lambda^{-}$,
$\lambda^+$ be positive constants satisfying $\frac{1}{2} \leq
\lambda^- < \lambda^+\leq \frac{3}{4}$. We set $\eta$ as
$$
\{
\begin{array}{cl}
\eta=1 & \mbox{on} \ \ Z\times [t_0-\lambda^-r, t_0+ \lambda^-r],
\\
\eta=0 & \mbox{on} \ \ X\setminus \left(Z\times [t_0-\lambda^+r, t_0+ \lambda^+r]\right),
\\
|d\eta|\leq \frac{2}{r(\lambda^+- \lambda^-)} & \mbox{on} \ \ 
Z\times \left([t_0-\lambda^+r, t_0-\lambda^-r]\sqcup [t_0+ \lambda^-r,t_0+ \lambda^+r]\right).
\end{array}
\right.
$$
Let $\Phi(q,r) := \displaystyle \left(\int_{Z\times[t_0-r,t_0+r]}u_L^q  
d\sigma_{\bar{g}}\right)^{\frac{1}{q}}$ for $q\geq\frac{2n}{n-2}$.
Set $\gamma =\frac{n}{n-2}>1$. We prove the following assertion.
\begin{Claim}{\label{Claim3}}
There exists a constant $\hat{C}_n > 0$ depending only on $n$ such
 that the following inequality holds for any $q \geq \frac{2n}{n-2}$:
\begin{equation}\label{new:eq20}
\Phi(\gamma q,\lambda^{-}r)\leq
\{
\frac{\hat{C}_n \ q}{Y \ r^2}
\left(
q^{\frac{n-2}{2}} +1
\right)
\}^{\frac{1}{q}}\cdot \Phi(q,\lambda^+r).
\end{equation}
\end{Claim}
\begin{lproof}{Proof of Claim \ref{Claim3}}
By H\"older's inequality, (\ref{new:eq19}) and the Young inequality,
we obtain
\begin{equation}\label{new:eq21}
\begin{array}{l}
\displaystyle
\!\!\!\!
\int_X \eta^2 u_L^{\alpha+\frac{n+2}{n-2}} d\sigma_{\bar{g}} = 
\displaystyle
\int_X u_L^{\frac{4}{n-2}}\cdot \eta^2 u_L^{\alpha+1} d\sigma_{\bar{g}}
\\
\\
\ \ \ \
\displaystyle
\leq \|\chi\cdot u_L\|^{\frac{4}{n-2}}_{L_{\bar{g}}^{\frac{2n^2}{(n-2)^2}}(X)}
\cdot \|\eta \cdot u_L^{\frac{\alpha+1}{2}}
\|^2_{L_{\bar{g}}^{\frac{2n^2}{n^2-2n+4}}(X)}
\\
\\
\ \ \ \ 
\displaystyle
\leq \left(
\frac{250\cdot \bar{C}_n}{Y}\cdot \frac{1}{r^2}
\right)^{\frac{2}{n}}\!\! \!\cdot
\{
\epsilon^{\prime}\cdot \|\eta u_L^{\frac{\alpha+1}{2}}
\|_{L_{\bar{g}}^{\frac{2n}{n-2}}(X)}\!\!\!\!
+ \left(\epsilon^{\prime}\right)^{-\frac{n-2}{2}}\!\! \!\cdot 
\|\eta u_L^{\frac{\alpha+1}{2}}\|_{L^2_{\bar{g}}(X)}
\}^2
\end{array}\!\!\!\!\!\!\!\!\!\!\!\!\!\!\!\!\!
\end{equation}
for any $\epsilon^{\prime}>0$. Set $\epsilon=1$ in
(\ref{new-eq3}). It then follows from (\ref{new-eq3}) and (\ref{new:eq21}) that
for any $\alpha>\frac{n+2}{n-2}$
\begin{equation}\label{new:eq22}
\!\!\!
\begin{array}{l}
\displaystyle
\left(
\int_X\left(
\eta u_L^{\frac{\alpha+1}{2}}
\right)^{\frac{2n}{n-2}}\!\! d\sigma_{\bar{g}}
\right)^{\frac{n-2}{n}}\!\!\!\!\!
\\
\\
\leq \ 
\displaystyle
\! \frac{C_n^{\prime\prime}\alpha}{Y}\!
\left[
\frac{Y^{\frac{n-2}{n}}}{r^{\frac{4}{n}}}\!
\{
(\epsilon^{\prime})^2\! \cdot \!
\|\eta u_L^{\frac{\alpha+1}{2}}\|^2_{L_{\bar{g}}^{\frac{2n}{n-2}}(X)}\!\!\!
+\!(\epsilon^{\prime})^{-(n-2)}\!\cdot \!
\|\eta u_L^{\frac{\alpha+1}{2}}\|_{L_{\bar{g}}^2(X)}^2
\}\right.
\\
\\
\displaystyle
\ \ \ \ \ \ \ \ \ \ \ \ \ \ \ \ \ \ \ \ \ \ \ \ \ \ \ \ \ \ 
\ \ \ \ \ \ \ \ \ \ \ \ \ \ \ \ \ \ \ \ \ \ \ \ \ \ \ \ \ \ \ \ \ \ \ 
\left.
+\int_X |d\eta|^2 u_L^{\alpha+1} d\sigma_{\bar{g}}
\!\right].\!\!\!\!\!\!\!\!\!\!\!\!\!\!\!\!\!\!\!\!\!\!\!\!\!\!
\end{array}
\end{equation}
Set $\displaystyle (\epsilon^{\prime})^2 =
\frac{Y^{\frac{2}{n}}r^{\frac{4}{n}}}{2 \ C_n^{\prime\prime}\alpha}$
in (\ref{new:eq22}). Then 
$\displaystyle (\epsilon^{\prime})^{-(n-2)} =
\frac{\left(2 \ C_n^{\prime\prime}\right)^{\frac{n-2}{2}}\alpha^{\frac{n-2}{2}}}
{Y^{\frac{n-2}{n}}r^{\frac{2(n-2)}{n}}}$,  and hence
$$
\|\eta u_L^{\frac{\alpha+1}{2}}\|^2_{L_{\bar{g}}^{2\gamma}(X)} \leq
\frac{C_n^{\prime\prime\prime}\alpha^{\frac{n}{2}}}{Y r^2} 
\|\eta u_L^{\frac{\alpha+1}{2}}\|^2_{L_{\bar{g}}^{2}(X)} +
\frac{C_n^{\prime\prime}\alpha}{Y} \| \ |d\eta|\cdot u_L^{\frac{\alpha+1}{2}}
\|^2_{L_{\bar{g}}^{2}(X)}.
$$
Finally we set $q=\alpha+1$, and then this implies the estimate
(\ref{new:eq20}).
\end{lproof}
Now we complete the proof of Proposition \ref{yam:P1}. The estimate
(\ref{new:eq20}) can be iterated to yield the desired result. Indeed, we
set
$$
\{
\begin{array}{rcl}
q_m &=& \gamma^m\frac{2n}{n-2},
\\
\lambda_m^-&=& \frac{1}{2}+ 2^{-(m+3)},
\\
\lambda_m^+&=& \frac{1}{2}+ 2^{-(m+2)} = \lambda_{m-1}^-
\end{array}
\right.
$$
for $m=0,1,2,\ldots$. Then we use (\ref{new:eq20}) to obtain
\begin{equation}\label{new:eq23}
\!\!\!\!\!\!\!\!
\begin{array}{rcl}
\displaystyle
\Phi\left(\gamma^m\frac{2n}{n-2},\frac{1}{2}r\right) \!\!&\leq&\!\!
\displaystyle
\left(
\frac{2\bar{C}_n}{Yr^2}
\right)^{\frac{n-2}{2n}\sum_{i=0}^{m-1}\gamma^{-i}}\!\!\!\!\!\!\cdot
\gamma^{\frac{n}{2} \sum_{i=0}^{m-1}i\gamma^{-i}} \cdot 
\Phi\left(\frac{2n}{n-2},r\right)
\\
\\
&\leq&\!\!
\displaystyle
\frac{C_n}{Y^{\frac{n-2}{4}}}\cdot \frac{1}{r^{\frac{n-2}{2}}}.
\end{array}
\end{equation}
Letting $m\to\infty$ in (\ref{new:eq23}), we obtain the estimate
(\ref{yam:P1-eq1}). This completes the proof of Proposition
\ref{yam:P1}.
\end{Proof}
Now we show a more precise decay estimate of the minimizers $u_L$ on
the cylindrical end $Z\times [1,\infty)$. 
\begin{Proposition}\label{prop3}
Let $u_L$ be the minimizer obtained in Lemma \ref{yam:L1} for each $L
\geq L_0$.  For any $a>0$ satisfying
$a<\sqrt{\frac{R_{\min}}{\alpha_n}}$, there exist constants
$\bar{K}=\bar{K}(a)>0$ and $\ell=\ell(a)>0$ such that
\begin{equation}\label{new:eq24}
u_L\leq \bar{K} \cdot e^{-at} \ \ \ \mbox{on} \ \ \ Z\times [\ell,\infty)
\end{equation}
for any $L>\ell$.
\end{Proposition}
\begin{Proof}
From (\ref{yam:P1-eq1}), there exists $\ell=\ell(a)>0$ such that
\begin{equation}\label{new:eq25}
\sup_{Z\times [\ell,\infty)} u_L \leq
\left(
\frac{R_{\min}-\alpha_n a^2}{Y(S^n)}
\right)^{\frac{n-2}{4}}
\end{equation}
for any $L>\ell$. We set
$$
\bar{K}= e^{a\ell} \left(
\frac{R_{\min}-\alpha_n a^2}{Y(S^n)}
\right)^{\frac{n-2}{4}}>0.
$$
Consider the function $w_L=u_L - \bar{K}\cdot e^{-at}$ on $Z\times
[\ell,L]$. It then follows from (\ref{yam:P1-eq2}) and
(\ref{new:eq25}) that
$$
\begin{array}{rcl}
\Delta_{\bar{g}}w_L &=&
\displaystyle
\frac{R_h}{\alpha_n} u_L -\frac{Q_L}{\alpha_n} u_L^{\frac{n+2}{n-2}} - a^2 \bar{K}\cdot 
e^{-at}
\\
\\
&\geq&
\displaystyle
\{
\left(
\frac{R_{\min}}{\alpha_n} -a^2
\right)
- \frac{Y(S^n)}{\alpha_n} u^{\frac{4}{n-2}}_L
\} u_L + a^2\cdot w_L 
\\
\\
&\geq&
\displaystyle
a^2\cdot w_L \ \ \ \ \mbox{on} \ \ \ Z\times(\ell,L),
\end{array}
$$
and $w_L|_{Z\times\{\ell\}} \leq 0$, $w_L|_{Z\times\{L\}} <
0$. By the maximum principle, we obtain that
$$
w_L = u_L - \bar{K}\cdot e^{-at} \leq 0 \ \ \ \mbox{on} \ \ \ Z\times[\ell,L].
$$
This completes the proof of Proposition \ref{prop3}.
\end{Proof}
\begin{lproof}{Proof of Theorem \ref{yam:Th1}}
The $L^p$ and Schauder interior estimates combined with 
Lemma \ref{yam:L3} imply that
$$
\| u_L\|_{C^{2,\alpha}(X(i))} \leq C(i)
$$
for any $L\geq L_0$ and $i=1,2,\ldots$, where each constant $C(i) > 0$
is independent of $L$. Then by the argument of diagonal subsequence,
there exist a subsequence ${u_{L_j}}$ and a nonnegative function $u
\in C^2(X) \cap L^{1, 2}_{\bar{g}}(X)$ such that $u_{L_j} \to u$ with
respect to the $C^2$-topology on each $X(i)$. Furthermore, the
function $u$ satisfies
\begin{equation}\label{new:eq26}
-\alpha_n\Delta_{\bar{g}}u +R_{\bar{g}} u = Y^{\cyl}_{\bar{C}}(X)\cdot u^{\frac{n+2}{n-2}}
\ \ \ \ \mbox{on} \ \ \ X.
\end{equation}
Combining Lemma \ref{yam:L3} with \cite[Proposition 3.75]{Au3}, we
notice that $u>0$ on $X$ or $u\equiv 0$ on $X$. Hence from
(\ref{new:eq26}), $u\in C^{\infty}(X)$.

On the other hand, it follows from (\ref{new:eq24}) that for any small
$\epsilon>0$ there exists $L(\epsilon)>0$ such that
$$
\int_{X(L(\epsilon))} u_L^{\frac{2n}{n-2}}d\sigma_{\bar{g}} \geq 1 -\epsilon
\ \ \ \mbox{for any} \ \  \ L\geq L_0.
$$
Then the $C^2$-convergence $u_{L_j}\to u$ on $X(L(\epsilon))$ implies
$$
\int_{X(L(\epsilon))} u^{\frac{2n}{n-2}}d\sigma_{\bar{g}} \geq 1 -\epsilon.
$$
This implies that $u>0$ on $X$ and 
$$
\int_{X} u^{\frac{2n}{n-2}}d\sigma_{\bar{g}} = 1 .
$$
This completes the proof of Theorem \ref{yam:Th1}.
\end{lproof}

\noindent
{\bf \ref{s4}.3. Singularities of the Yamabe metric.}  Let
$\check{g}=u^{\frac{4}{n-2}}\bar{g}$ be the Yamabe metric obtained in
Theorem \ref{yam:Th1}.  We also prove the existence of a Yamabe metric
for the case $\lambda({\cal L}_h)=0$ (Theorem \ref{yam:Th3}).
Moreover, we study the singularities of the metric $\check{g}$ near
infinity of the tame end $Z\times[0,\infty)$. There are two very
different cases here: $\lambda({\cal L}_h)>0$ and $\lambda({\cal
L}_h)=0$.

We start with the case $\lambda({\cal L}_h)>0$. A canonical model here
is provided by the canonical open cone
$$
\Cone(Z)= (Z\times (0,1), r^2\cdot h + dr^2)
$$
over $(Z, h)$.  To obtain the cylindrical coordinates, we set $r =
e^{-t}$.  Then we identify the cone $\Cone(Z)$ with
$$
\Cone(Z)\cong (Z\times (0,\infty), e^{-2t}(h + dt^2)). 
$$
We use the cone $\Cone(Z)$ as a canonical model to introduce the
following definition.
\begin{figure}[ht]
\hspace*{20mm}
\PSbox{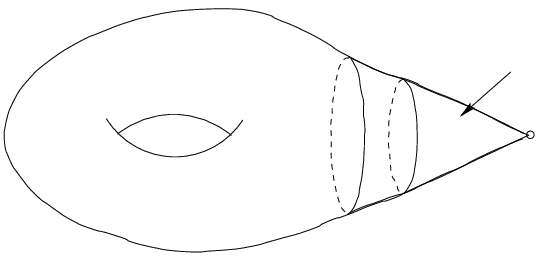}{10mm}{25mm}  
\begin{picture}(0,1)
\put(80,5){{\small $Z \times[1,\infty)$}}
\put(110,55){{\small $Z \times[L,\infty)$}}
\put(17,55){{\small $(N,g)$}}
\end{picture}
\caption{Almost conical metric.}\label{fig5-1}
\end{figure}
\begin{Definition}\label{def:cone}
{\rm Let $N$ be an open manifold with a connected tame end $Z \times
  [0, \infty)$. A metric $g\in \Riem(N)$ is called an \emph{almost
  conical metric} if there exist a coordinate system $(x,t)$ on
  $Z\times [0,\infty)$, a metric $h\in \Riem(Z)$, constants
  $0<\alpha\leq\beta$, $0<C_1\leq C_2$  and a positive function
  $\phi\in C_+^{\infty}(Z\times [1,\infty))$ such that}
\begin{enumerate}
\item[{\bf (i)}] $g(x,t)=\phi(x,t)(h(x)+dt^2)$ on $Z\times
[1,\infty)$,
\item[{\bf (ii)}] $C_1\cdot e^{-\beta t}\leq \phi(x,t) \leq C_2\cdot
e^{-\alpha t}$ on $Z\times [1,\infty)$.
\end{enumerate}
\end{Definition}
Below we state the result on singularities without using Convention
\ref{convention}.
\begin{Theorem}\label{yam:Th2}
Under the same assumptions as in Theorem \ref{yam:Th1}, let $u\in
C_+^{\infty}(X)\cap L^{1,2}_{\bar{g}}(X)$ be the Yamabe minimizer
obtained in Theorem \ref{yam:Th1}. Then for any constant $a>0$
satisfying $\displaystyle 0<a < \left(\min_Z \phi_h
\right)^{\frac{2}{n-2}}\cdot \sqrt{\frac{\lambda({\cal L}_h)}{\alpha_n}}$,
there exist constants $\ov{C}$, $\underline{C}> 0$ and $\ell>0$ such
that
$$
\underline{C} \cdot e^{-a_0t}\leq u(x,t)\leq \ov{C} \cdot e^{-at}
\ \ \mbox{for} \ \ (x,t)\in Z\times [\ell,\infty), \ \ \mbox{where} \ \ 
a_0= \sqrt{\frac{\lambda({\cal L}_h)}{\alpha_n}}.
$$
Here \ \ $\displaystyle\ov{C} = \ov{C}(a,n,\lambda({\cal L}_h),
\min_Z\phi_h)$, \ \ $\displaystyle \underline{C}=
\underline{C}(n,\lambda({\cal L}_h),\min_Z \phi_h,\min_{Z\times \{1\}}
u)$, $\displaystyle \ell=\ell(a,n,\lambda({\cal L}_h),\min_Z \phi_h)$,
and these constants also depend on $Y^{\cyl}_{[h+dt^2]}(Z\times \R)$,
$Y^{\cyl}_{\bar{C}}(X)$.  In particular, the Yamabe metric
$\check{g}=u^{\frac{4}{n-2}}\cdot \bar{g}$ on $X$ is an almost conical
metric.
\end{Theorem}
\begin{Proof}
In the proof of Proposition \ref{prop3}, we replace the cylindrical
metric $\bar{g}$ by another metric $\tilde{g}\in [\bar{g}]$
satisfying $\tilde{g}= \phi_h^{\frac{4}{n-2}}\cdot \bar{g}$ on
$Z\times [1,\infty)$, where $\phi_h$ is the positive smooth function
given in Convention \ref{convention}. 
\vspace{2mm}

\noindent
Let $\tilde{u}\in C_+^{\infty}(X)\cap L^{1,2}_{\bar{g}}(X)$ be the
Yamabe minimizer with respect to $\tilde{g}$ with
$\tilde{u}^{\frac{4}{n-2}}\cdot \tilde{g} = u^{\frac{4}{n-2}}\cdot
\bar{g}$. Recall that $R_{\tilde{g}} = \lambda({\cal L}_h)\cdot
\phi_h^{-\frac{4}{n-2}}$, $\displaystyle \max_Z\phi_h =1$ and
$$
\Delta_{\tilde{g}} f = \phi_h^{-\frac{4}{n-2}}\cdot f^{\prime\prime}
\ \ \ \ \ \mbox{for} \ \ f=f(t)\in C^{\infty}(Z\times [1,\infty)).
$$
By these properties, for any $\displaystyle 0<a < \left(\min_Z \phi_h
\right)^{\frac{2}{n-2}}\cdot \sqrt{\frac{\lambda({\cal
L}_h)}{\alpha_n}}$ there exist $\displaystyle
\bar{K}=\bar{K}(a,n,\lambda({\cal L}_h),\min_Z \phi_h)>0$ and
$\displaystyle\ell=\ell (a,n,\lambda({\cal L}_h),\min_Z\phi_h)>0$ such
that
$$
\tilde{u}\leq \bar{K}\cdot e^{-at} \ \ \ \mbox{on} \ \ \
Z\times[\ell,\infty) .
$$
A similar argument to the above and Proposition \ref{prop3}, combined
with the condition $Y^{\cyl}_{\bar{C}}(X)>0$, implies that there
exists a constant $\underline{K}=\underline{K}(n,\lambda({\cal
L}_h),\!\!\displaystyle\min_{Z\times \{1\}} u )>0$ such that
$$
u\geq \underline{K}\cdot(\min_Z\phi_h)\cdot e^{-a_0t} \ \ \ \mbox{on} \ \ \
Z\times[1,\infty).
$$
This completes the proof.
\end{Proof}
Next, we study the case when $\lambda({\cal L}_h)=0$. Here a canonical
model comes from hyperbolic geometry, namely, it is given by
\emph{cusp ends of hyperbolic manifolds}, see \cite[Chapter D]{BP}.
\vspace{2mm}

\noindent
A cusp end of a hyperbolic $n$-manifold of curvature $-1$ is given as
$$
\left(
(R^{n-1}/\Gamma)\times [1,\infty),\frac{1}{t^2}(h_{\fla} +dt^2)
\right).
$$
Here $(R^{n-1}/\Gamma, h_{\fla})$ is a closed Riemannian manifold
uniformized by a flat torus $T^{n-1}$. With this understood, we
introduce the following definition.
\begin{Definition}\label{yam:def2}
{\rm Let \ $N$ \ be an open manifold with a connected tame end \ $Z \times
 [0, \infty)$. \ A metric \ $g\in \Riem(N)$ \ is called an \ \emph{almost cusp
 metric} \ if there exist a coordinate system \ $(x,t)$ \ on the
 cylinder \ $ Z\times [0,\infty)$, \ a metric $h\in \Riem(Z)$,
 constants $0<C_1\leq C_2$ and a positive function $\phi\in
 C_+^{\infty}(Z\times [1,\infty))$ such that}
\begin{enumerate}
\item[{\bf (i)}] $g(x,t)=\phi(x,t)(h(x)+dt^2)$ on $Z\times [1,\infty)$,
\item[{\bf (ii)}] $C_1\cdot t^{-2}\leq \phi(x,t) \leq C_2\cdot t^{-2}$
on $Z\times [1,\infty)$.
\end{enumerate}
\end{Definition}
Recall that $\lambda({\cal L}_h)=0$ implies that
$Y^{\cyl}_{[h+dt^2]}(Z\times \R)=0$. Hence $Y^{\cyl}_{\bar{C}}(X)\leq
0$. 

\begin{figure}[h]
\hspace*{5mm}
\PSbox{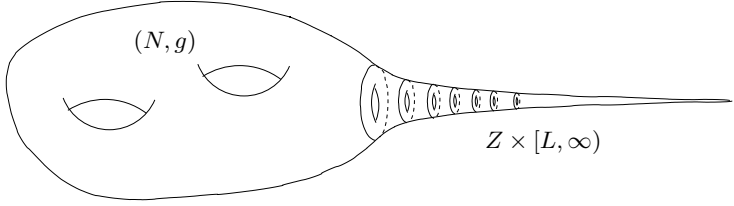}{10mm}{28mm}  
\begin{picture}(0,1)
\put(150,20){{\small $Z \times[L,\infty)$}}
\put(17,58){{\small $(N,g)$}}
\end{picture}
\vspace{4mm}
\caption{Almost cusp metric.}\label{fig5-2}
\end{figure}
\begin{Theorem}\label{yam:Th3}
Let \ $X$ \ be an open manifold with a \ connected tame end $Z\times
[0,\infty)$ and $h\in \Riem(Z)$ a metric with $\lambda({\cal
L}_h)=0$. Let $\bar{g}$ be a cylindrical metric on $X$ with
$\p_{\infty}\bar{g} = h$, and $\bar{C} = [\bar{g}]$.  Assume that
$Y^{\cyl}_{\bar{C}}(X)<0$.  Then there exists a Yamabe minimizer $u\in
C_+^{\infty}(X)\cap L^{1,2}_{\bar{g}}(X)$ with $\int_X
u^{\frac{2n}{n-2}} d\sigma_{\bar{g}}=1$ such that $Q_{(X,\bar{g})}(u)
= Y_{\bar{C}}^{\cyl}(X)$. In particular, the minimizer $u$ satisfies
the Yamabe equation:
$$
\L_{\bar{g}} u = -\frac{4(n-1)}{n-2}\Delta_{\bar{g}} u + R_{\bar{g}} u = 
Y_{\bar{C}}^{\cyl}(X)u^{\frac{n+2}{n-2}}.
$$
Moreover, there exist constants $0<\underline{C}\leq\ov{C}$ such that
$$
\underline{C}\cdot t^{-\frac{n-2}{2}} \leq u(x,t) \leq \ov{C}\cdot t^{-\frac{n-2}{2}}
\ \ \ \mbox{on} \ \ \ Z\times [1,\infty)\ \mbox{and}
$$
$
\displaystyle \ov{C}=\ov{C}(n, Y^{\cyl}_{\bar{C}}(X),\min_Z
\phi_h ,\max_{Z\times \{1\}} u), \ \displaystyle
\underline{C}=\underline{C}(n, Y^{\cyl}_{\bar{C}}(X), \min_Z\phi_h,
\min_{Z\times \{1\}} u). 
$
In particular, the Yamabe metric $\check{g}= u^{\frac{4}{n-2}}\cdot
\bar{g}$ on $X$ is an almost cusp metric.
\end{Theorem}
\begin{Proof}
As in Convention \ref{convention}, we may assume that $R_h\equiv 0$ on
$Z$. For any $L\geq 1$, there exists a nonnegative function $u_L\in
 C^0(X) \cap C_+^{\infty}(\intt(X(L)))\cap L_{\bar{g}}^{1,2}(X)$ satisfying the
same properties as in Lemma \ref{yam:L1}. Similarly to the proofs of
Lemmas \ref{yam:L2}, \ref{yam:L3}, there exist constants $L_0>> 1$
and $K_0>0$ such that
\begin{equation}\label{new:eq27}
u_L\leq K_0 \ \ \ \mbox{on} \ \ \ X
\end{equation}
for any $L\geq L_0$.  Here we may assume that
\begin{equation}\label{new:eq27-1}
Q_L \leq\frac{1}{2} Y_{\bar{C}}^{\cyl}(X)< 0 \ \ \ \mbox{for} \ \ \ L\geq L_0.
\end{equation}
We set the constant $\ov{C}$ by
$$
\ov{C}=\max\{\left(
\frac{n(n-2)\alpha_n}{2|Y_{\bar{C}}^{\cyl}(X)|}
\right)^{\frac{n-2}{4}}\!\!\!, K_0
\} >0 .
$$
Consider the function $\ov{w}_L= u_L - \ov{C}\cdot t^{-\frac{n-2}{2}}$
on $Z\times [1,L]$. It then follows from  
(\ref{yam:P1-eq2}) and (\ref{new:eq27-1}) that
$$
\begin{array}{rcl}
\Delta_{\bar{g}} \ov{w}_L & = &\displaystyle
\frac{|Q_L|}{\alpha_n} u_L^{\frac{n+2}{n-2}} -
\frac{n(n-2)}{4\ov{C}^{_{\frac{4}{n-2}}}}
\left(
\ov{C}\cdot t^{-\frac{n-2}{2}}
\right)^{\frac{n+2}{n-2}}
\\
\\
& \leq &\displaystyle
\frac{|Y_{\bar{C}}^{\cyl}(X)|}{2\alpha_n}
\{
u_L^{\frac{n+2}{n-2}} - 
\left(
\ov{C}\cdot t^{-\frac{n-2}{2}}
\right)^{\frac{n+2}{n-2}}
\} \ \ \ \mbox{on} \ \ Z\times [1,L],
\end{array}
$$
and $\ov{w}_L|_{Z\times\{1\}}\leq 0$, $\ov{w}_L|_{Z\times\{L\}}<
0$. By the maximum principle, we obtain that
\begin{equation}\label{new:eq27-2}
\ov{w}_L= u_L - \ov{C}\cdot t^{-\frac{n-2}{2}} \leq 0 \ \ \mbox{on} \
\ Z\times [1,L].
\end{equation}
From (\ref{new:eq27-2}), for any small $\epsilon>0$ there
exists $L(\epsilon)>0$ such that
\begin{equation}\label{new:eq29}
\int_{X(L(\epsilon))} u_L^{\frac{2n}{n-2}}d\sigma_{\bar{g}} \geq
1 - \ov{C}^{_{\frac{2n}{n-2}}}\cdot\Vol_h(Z)\cdot
\int_{L(\epsilon)}^{\infty} \frac{1}{t^n}\ dt \geq 1-\epsilon
\end{equation}
for any $L\geq L_0$. The estimates  (\ref{new:eq27}) and
(\ref{new:eq27-2}) combined with (\ref{new:eq29}) imply  that $u_L$
converges to a Yamabe minimizer $u\in C_+^{\infty}(X)\cap
L^{1,2}_{\bar{g}}(X)$ with $\int_X u^{\frac{2n}{n-2}}
d\sigma_{\bar{g}}=1$ in the $C^2$-topology on each $X(i)$ for $i=1,2,\ldots$.
Moreover, $u$ satisfies
\begin{equation}\label{new:eq29-1}
\begin{array}{l}
\displaystyle
u\leq \ov{C}\cdot t^{-\frac{n-2}{2}} \ \ \ \mbox{on} \ \ Z\times[1,\infty),
\\
\\
\displaystyle
-\alpha_n \Delta_{\bar{g}} u = Y_{\bar{C}}^{\cyl}(X) u^{\frac{n+2}{n-2}}  \ \ \ 
\mbox{on} \ \ X.
\end{array}
\end{equation}
We set $\displaystyle \ \underline{C} =\min\{ \left(
\frac{n(n-2)\alpha_n}{4|Y_{\bar{C}}^{\cyl}(X)|}
\right)^{\frac{n-2}{4}}, \min_{Z\times\{1\}} u \}>0$. Consider the
function 
$$
\underline{w}= \underline{C}\cdot t^{-\frac{n-2}{2}} -u
$$
on $Z\times [1,\infty)$. It then follows from (\ref{new:eq29-1}) that
$$
\begin{array}{rcl}
\displaystyle
\Delta_{\bar{g}} \underline{w} & =&
\displaystyle
\frac{n(n-2)}{4\underline{C}^{\frac{4}{n-2}}}
\left(
\underline{C}\cdot t^{-\frac{n-2}{2}}
\right)^{\frac{n+2}{n-2}}
-\frac{|Y_{\bar{C}}^{\cyl}(X)|}{\alpha_n}u^{\frac{n+2}{n-2}}
\\
\\
& \geq &\displaystyle
\frac{|Y_{\bar{C}}^{\cyl}(X)|}{\alpha_n} 
\{
\left(
\underline{C}\cdot t^{-\frac{n-2}{2}}
\right)^{\frac{n+2}{n-2}} - u^{\frac{n+2}{n-2}}
\} \ \ \ \mbox{on} \ \ Z\times [1,\infty),
\displaystyle
\end{array}
$$
and $\underline{w}|_{Z\times\{1\}} \leq 0$,
$\displaystyle\lim_{t\to\infty} \underline{w}(x,t)=0$. The maximum
principle now implies that
$$
 \underline{w} = \underline{C} \cdot t^{-\frac{n-2}{2}} - u \leq 0 \ \
 \ \mbox{on} \ \ \ Z\times [1,\infty).
$$
This completes the proof.
\end{Proof}
{\bf Remark.}  Let $(X, \bar{g})$ be a cylindrical manifold with tame
ends $Z \times [0, \infty)$ and $h = \p_{\infty}\bar{g} \in \Riem(Z)$
a metric with $\lambda({\mathcal L}_h) \geq 0$.  We emphasize that the
"strict" generalized Aubin's inequality $Y^{\cyl}_{[\bar{g}]}(X) <
Y^{\cyl}_{[h + dt^2]}(Z \times \R)$ is a crucial condition for the
above solutions of the Yamabe problem.  In the case when
\begin{equation}\label{remark_new}
Y^{\cyl}_{[\bar{g}]}(X) = Y^{\cyl}_{[h + dt^2]}(Z \times \R),
\end{equation}
Proposition \ref{yam:new1-1} implies that the Yamabe problem can not
be solved in general. Hence, it is a natural problem to characterize
conformally cylindrical manifolds $(X, \bar{g})$
satisfying (\ref{remark_new}).  In the next section, we solve this
problem in the supremum case, i.e.  when $Y^{\cyl}_{[\bar{g}]}(X) =
Y^{\cyl}_{[h + dt^2]}(Z \times \R) = Y(S^n)$.
\section{Canonical cylindrical manifolds}\label{s5}
Let $(Z, h)$ be a closed Riemannian manifold of $\dim Z = n-1 \geq 2$
with $\lambda({\mathcal L}_h) > 0$. Throughout in this section, we
also assume that $Z$ is connected (unless we specify otherwise).
First, we study the Yamabe problem on the \emph{canonical cylindrical
manifold} $(Z\times \R, \bar{h}=h+dt^2)$.  Clearly, canonical
cylindrical manifolds do not satisfy the strict generalized Aubin's
inequality in Theorem \ref{yam:Th1}.  Instead, these Riemannian
manifolds are invariant under parallel translations along the
$t$-coordinate. Using this symmetry for the renormalization technique,
we solve the Yamabe problem under the conditions $\lambda({\mathcal
L}_h) > 0$ and $Y^{\cyl}_{[\bar{h}]}(Z \times \R) < Y(S^n)$.  Second,
we characterize canonical cylindrical manifolds which satisfy the
supremum condition
$$
Y^{\cyl}_{[\bar{h}]}(Z \times \R) = Y(S^n)
$$
Furthermore, we characterize general cylindrical manifolds $(X,
\bar{g})$ which satisfy
$$
Y^{\cyl}_{[\bar{g}]}( X) = Y^{\cyl}_{[h + dt^2]}(Z \times \R) =
Y(S^n).
$$
Recall that there is no Yamabe minimizer for any canonical
cylindrical manifold $(Z\times \R, \bar{h}=h+dt^2)$ with $
\lambda({\cal L}_h)=0$. Hence we consider only the case $\lambda({\cal
L}_h)>0$.
\begin{Theorem}\label{can:Th1}
Let $(Z\times \R, \bar{h}=h+dt^2)$ be a canonical cylindrical manifold
of $\dim (Z \times \R) = n \geq 3$ with $\lambda({\cal
L}_h)>0$. Assume that
$$
Y^{\cyl}_{[\bar{h}]}(Z \times \R) < Y(S^n) 
   = Y^{\cyl}_{[h_+ + dt^2]}(S^{n-1} \times \R).
$$
Then there exists a Yamabe minimizer $u\in C_+^{\infty}(Z\times
\R)\cap L_{\bar{h}}^{1,2}(Z\times \R)$ with $\displaystyle
\int_{Z\times \R} u^{\frac{2n}{n-2}}d\sigma_{\bar{h}}= 1$ such that
$$
Q_{(Z\times \R, \bar{h})}(u) = Y^{\cyl}_{[\bar{h}]}( Z\times \R).
$$
Furthermore, for any constant $a>0$ satisfying $0<a <
\left(\min_Z \phi_h \right)^{\frac{2}{n-2}}\cdot
\sqrt{\frac{\lambda({\cal L}_h)}{\alpha_n}}$, there exist constants
$\ov{C}$, $\underline{C}> 0$ and $\ell>0$ such that
$$
\begin{array}{l}
\underline{C} \cdot e^{-a_0|t|}\leq u(x,t)\leq \ov{C} \cdot e^{-a|t|}
\ 
\end{array}
$$
for $(x,|t|)\in Z\times [\ell,\infty)$, where $ a_0\!=
\!\sqrt{\frac{\lambda({\cal L}_h)}{\alpha_n}}$.  Here
$\displaystyle\ov{C}\! = \!\ov{C}(a,n,\lambda({\cal L}_h),\!
\min_Z\phi_h)$, \ \ $\displaystyle \underline{C}=
\underline{C}(n,\lambda({\cal L}_h),\min_Z \phi_h,\min_{Z\times \{0\}}
u)$, $\displaystyle \ell=\ell(a,n,\lambda({\cal L}_h),\min_Z \phi_h,
u)$, and the constants $\ov{C}$, $\ell$ also depend on
$Y^{\cyl}_{[h+dt^2]}(Z\times \R)$.  In particular, the Yamabe metric
$\check{g}=u^{\frac{4}{n-2}}\cdot \bar{g}$ on $Z\times \R$ is an
almost conical metric.
\end{Theorem}
\begin{Proof}
As in Convention \ref{convention}, we may assume that $R_{\min}=\min_Z
R_h>0$. By the assumption that $Y:= Y^{\cyl}_{[\bar{h}]}( Z\times \R)<
Y(S^n)$, there exists $L_0>>1$ such that
\begin{equation}\label{can:eq1}
\{
\begin{array}{c}
\displaystyle \tilde{Q}_L \ \ := \inf_{\begin{array}{c} ^{f\in
L^{1,2}_{\bar{h}}(Z\times \R), \ \ f\not\equiv 0,} \\ ^{f\equiv 0 \
{\mathrm o}{\mathrm n} \ Z\times ((-\infty, -L]\sqcup [L,\infty))}
\end{array}}\!\!\!\!\!\!\!\!
Q_{(Z\times \R, \bar{h})}(f) < Y(S^n),
\\
\\
0<Y \leq \tilde{Q}_L \leq 2Y 
\end{array}
\right.
\end{equation} 
for any $L\geq L_0$. Then, 
\begin{equation}\label{can:eq2}
\tilde{Q}_L \to Y \ \ \ \mbox{as} \ \ \ L\to\infty.
\end{equation}
Note that the metric \ $\bar{h}=h+dt^2$ \ is invariant under parallel
translations along the $t$-coordinate. \ Using a similar argument to the
proof of \ \cite[Chapter 5, Theorem 2.1]{SY4} and this symmetry for the
renormalization technique, we obtain that for any $L\geq L_0$ there
exists a nonegative function $u_L\in C^0(Z\times \R)\cap
C^{\infty}(Z\times (-L,L))$ such that
\begin{equation}\label{can:eq3}
\!\!\!\!\!\!\!\!\!
\{
\begin{array}{l}
\displaystyle
Q_{(Z\times \R, \bar{h})}(u_L) = \tilde{Q}_L, \ \ \ \int_{Z\times \R}
u^{\frac{2n}{n-2}}_L d\sigma_{\bar{h}} = 1
\\
\displaystyle
u_L > 0 \ \ \mbox{on} \ \ Z\times(-L,L), \ \ u_L\equiv 0 \ \ \mbox{on} \ \ 
Z\times ((-\infty, -L]\sqcup [L,\infty)),
\\
\displaystyle
-\alpha_n\Delta_{\bar{h}} u_L + R_h u_L = \tilde{Q}_L u_L^{\frac{n+2}{n-2}} \ \ 
\mbox{on} \ \ Z\times (-L,L),
\\
u_L \leq K \ \ \mbox{on} \ \ Z\times \R,
\end{array}
\right.\!\!\!\!\!\!\!\!\!\!\!\!\!\!\!
\end{equation}
where $K>0$ is a constant independent of $L$.

By using a parallel translation along the $t$-coordinate for each $L
\geq L_0$, we may assume that
\begin{equation}\label{can:eq4}
u_L(x_L, 0) =\max_{Z\times \R} u_L \leq K,
\end{equation}
where $x_L\in Z$. In this case, the support $\supp(u_L)$ may be no
longer equal to $Z\times[-L,L]$ (as in (\ref{can:eq3})).  However,
clearly $\supp(u_L)\subset Z\times[-2L,2L]$, and $\diam (\supp(u))= 2L
\to \infty$ as $L\to \infty$.

Set $\supp(u_L)=Z\times
[t_L-2L, t_L]$, where $t_L >0$ for any $L\geq L_0$. By taking a
subsequence if necessary, there exist constants $T^-,T^+\in
[-\infty, \infty]$ with $T^-< T^+$ such that
$$
\begin{array}{c}
\displaystyle
T^- = \lim_{L\to\infty} (t_L -2L), \ \ \ T^+ = \lim_{L\to\infty} t_L,
\\
\\
-\infty \leq T^-\leq 0\leq T^+ \leq \infty.
\end{array}
$$
Clearly if $T^->-\infty$ (resp. $T^+<\infty$), then $T^+=\infty$
(resp. $T^-=-\infty$). 

We consider only the case $T^+=\infty$ since the argument in the case
$T^- = - \infty$ is similar. Then the properties
(\ref{can:eq1})--(\ref{can:eq4}) allow us to apply the $L^p$ and
Schauder interior estimates to conclude that
$$
\| u_L\|_{C^{2,\alpha}(Z\times [T^-(k), k])} \leq C(k)
$$
for any $L>> 1$ and $k=1,2,\ldots$. Here 
$$
T^-(k)=\{
\begin{array}{ll}
-k &\mbox{if} \ \ T^-=-\infty,
\\
T^- +\frac{1}{k} &\mbox{if} \ \ T^->-\infty,
\end{array}
\right.
$$
and each constant \ $C(k)$ \ is independent of \ $L$.  \ Then we use
the argument of diagonal subsequence to obtain a subsequence
$\{u_{L_j}\}$ and a nonnegative function $u \in C^0(Z \times \R) \cap
C^{\infty}(Z \times (T^-, \infty))$ such that $u_{L_j} \to u$
in the $C^2$-topology on each cylinder $Z \times [T^-(k), k]$.  Here
$u\equiv 0$ on $Z\times (-\infty, T^-]$ if $T^->-\infty$. Furthermore,
the nonnegative function $u$ satisfies the equation
\begin{equation}\label{can:eq5}
-\alpha_n \Delta_{\bar{h}}u + R_{h} u = Y u^{\frac{n+2}{n-2}} \ \ \ 
\mbox{on} \ \ Z\times (T^-,\infty).
\end{equation}
Combining (\ref{can:eq4}) and (\ref{can:eq5}) with \cite[Proposition
3.75]{Au3}, we have that $u>0$ or $u\equiv 0$ on $ Z\times
(T^-,\infty)$.

Now we recall that $Y>0$ by Proposition \ref{yam:P1-new} and $Y\leq
\tilde{Q}_L\leq 2Y$. Then at the maximum point $(x_L,0)$ for $u_L$,
$$
u_L(x_L,0) \geq \left(
\frac{R_h(x_L)}{\tilde{Q}_L}
\right)^{\frac{n-2}{4}} \geq \left(
\frac{R_{\min}}{2Y}
\right)^{\frac{n-2}{4}} > 0.
$$
This implies that $u>0$ on $Z\times (T^-,\infty)$, and hence $u\in
C^{\infty}_+(Z\times (T^-,\infty))$. We notice also that
$$
0< \int_{Z\times \R} u^{\frac{2n}{n-2}} d\sigma_{\bar{h}} \leq
\lim\!\!\!\!\!\!\!\!\!\inf_{L_j\to\infty \ \ \ \ \ \ }\!\!\!\!\!\! \int_{Z\times\R}  u^{\frac{2n}{n-2}}_{L_j} d\sigma_{\bar{h}} =1.
$$ 
If $\int_{Z\times \R} u^{\frac{2n}{n-2}} d\sigma_{\bar{h}}<1$, then $
 Q_{(Z\times \R,\bar{h})}(u)< Y $ from (\ref{can:eq5}).  This
 contradicts the definition of $Y$, and hence
\begin{equation}\label{can:eq6}
\int_{Z\times \R} u^{\frac{2n}{n-2}} d\sigma_{\bar{h}}=1.
\end{equation}
Then we have
$$
Y= Q_{(Z\times \R, \bar{h})}(u) = \inf_{\begin{array}{c}
^{f\in L^{1,2}_{\bar{h}}(Z\times \R), \ \ f\not\equiv 0,} 
\\
^{f\equiv 0 \ {\mathrm o}{\mathrm n} \ Z\times (-\infty, T^-]} 
\end{array}}\!\!\!\!\!\!\!\!
Q_{(Z\times \R, \bar{h})}(f) \ .
$$
Applying a similar argument in the proof of Proposition \ref{yam:P1}
 directly to the minimizer $u$ with (\ref{can:eq5}) and
 (\ref{can:eq6}), we obtain that
$$
u(x,t)=o(1) \ \ \mbox{as} \ \ t\to\infty.
$$
To complete the proof, we have to show that $T^-=-\infty$. Suppose that
$T^->-\infty$. By the boundary regularity for the Yamabe equation, we
notice that $u|_{Z\times [T^-,\infty)}\in C^{\infty}(Z\times
[T^-,\infty))$. If $\frac{\p u}{\p t}\equiv 0$ on $Z\times \{T^-\}$,
then $u$ satisfies 
$$
-\alpha_n \Delta_{\bar{h}} u + R_h u = Yu^{\frac{n+2}{n-2}} \ \ \mbox{on} \ \
Z\times \R
$$
in the distributional sense. Then the standard elliptic regularity
implies that $u\in C^{\infty}(Z\times \R)$,  and thus $u>0$ or $u\equiv
0$ on $Z\times \R$, and hence $T^-=-\infty$. 
\begin{figure}[h]
\hspace*{10mm}
\PSbox{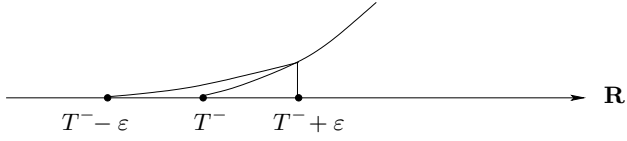}{10mm}{18mm}  
\begin{picture}(0,1)
\put(-10,-10){{\small $T^-\!\!-\epsilon$}}
\put(40,-10){{\small $T^-$}}
\put(70,-10){{\small $T^-\!+\epsilon$}}
\put(195,0){{\small $\R$}}
\end{picture}
\caption{The function $u_{\epsilon}(x,t)$.}\label{fig6-1}
\end{figure}
If $\frac{\p u}{\p t}\not\equiv 0$ on $Z\times \{T^-\}$, we define a
Lipschitz function $u_{\epsilon}$ by
$$
u_{\epsilon}(x,t) =
\{
\begin{array}{cl}
u(x,t) &\mbox{if} \ \ t\geq T^-+\epsilon,
\\
u(x,\frac{1}{2}(t+ T^-+\epsilon))  &\mbox{if} \ \ T^--\epsilon\leq t\leq T^-+\epsilon,
\\
0 &\mbox{if} \ \ t\geq T^--\epsilon.
\end{array}
\right.
$$
Then we have
$$
\begin{array}{c}
\displaystyle
E_{(Z\times \R,\bar{h})}(u_{\epsilon})= E_{(Z\times \R,\bar{h})}(u)
-\frac{1}{2}\alpha_n\left(
\int_{Z\times\{T^-\}}\left|\frac{\p u}{\p t}\right|^2 d\sigma_h 
\right)\epsilon + O(\epsilon^2),
\\
\\
\mbox{and}\ \ \ \displaystyle
\|u_{\epsilon}\|^2_{L_{\bar{h}}^{\frac{2n}{n-2}}(Z\times \R)} =
\|u\|^2_{L_{\bar{h}}^{\frac{2n}{n-2}}(Z\times \R)} + O(\epsilon^{1+\frac{2n}{n-2}}),
\end{array}
$$
and hence
$$
Q_{(Z\times \R,\bar{h})}(u_{\epsilon})= Q_{(Z\times \R,\bar{h})}(u) -
\frac{\frac{1}{2}\alpha_n\left( \int_{Z\times \{T^-\}} \left|\frac{\p
u}{\p t}\right|^2 d\sigma_h
\right)}{\|u\|^2_{L_{\bar{h}}^{\frac{2n}{n-2}}(Z\times \R)}}\cdot\epsilon +
O(\epsilon^2) < Y
$$
for $0<\epsilon<< 1$. This is a contradiction. Therefore, $T^-=-\infty$.

The above argument also implies that for each $i\geq 1$ there exist
constants $L_i^-$, $L_i^+$ and nonnegative functions 
$$
u_i\in C^0(Z\times \R)\cap C^{\infty}_+(Z\times (L_i^-, L_i^+))
$$
such that $u_i \rightarrow u$ in the $C^2$-topology 
on each cylinder $Z \times [-j, j]$ for $j = 1, 2, \ldots$, and 
$$
\{
\begin{array}{l}
\displaystyle
\lim_{i\to\infty} L_i^- = -\infty, \ \ \lim_{i\to\infty} L_i^+ = \infty,
\\
\displaystyle
Q_{(Z\times \R,\bar{h})}(u_i) = \tilde{Q}_i:=\!\!\!\!\!\!\!\!\!\!\!\!\!\!\!\!
\inf_{\begin{array}{c}
^{f\in L^{1,2}_{\bar{h}}(Z\times \R), f\not\equiv 0,}
\\
^{f\equiv 0 \ {\mathrm o}{\mathrm n} \ Z\times((-\infty,L_i^-]\sqcup[L_i^+,\infty))}
\end{array}}\!\!\!\!\!\!\!\!
Q_{(Z\times \R,\bar{h})}(f) ,
\\
\displaystyle -\alpha_n\Delta_{\bar{h}} u_i +R_h u_i = \tilde{Q}_i
u_i^{\frac{n+2}{n-2}} \ \ \mbox{on} \ \ Z\times (L_i^-, L_i^+), 
\\
\int_{Z\times \R} u_i^{\frac{2n}{n-2}} d\sigma_{\bar{h}} = 1 , \ \ u_i
\leq K \ \ \mbox{on} \ \ Z\times \R.
\end{array}
\right.
$$
Hence the decay estimate of $u$ follows from a similar argument to the
proof of Theorem \ref{yam:Th2}. This completes the proof of Theorem
\ref{can:Th1}.
\end{Proof}
Next, we show the following.
\begin{Theorem}\label{can:Th2}
Let $(Z\times \R, \bar{h}=h+dt^2)$ be a canonical cylindrical manifold
of $\dim (Z\times \R) =n \geq 3$. Assume that
$$
Y^{\cyl}_{[\bar{h}]}(Z \times \R) = Y(S^n) 
 = Y^{\cyl}_{[h_+ + dt^2]}(S^{n-1} \times \R)
$$
Then $(Z,h)$ is homothetic to the standard sphere $S^{n-1}(1)=(S^{n-1},h_+)$.
\end{Theorem}
As a corollary, we also obtain the following.
\begin{Corollary}\label{can:Th4}
Let $(X,\bar{g})$ be a cylindrical manifold of $\dim X =n \geq 3$ with
tame ends $Z\times [0,\infty)$ and $\p_{\infty} \bar{g}=h\in
\Riem(Z)$. Let $Z =\bigsqcup_{i=1}^k Z_j$, where each $Z_j$ is
connected. Assume that
$$
Y^{\cyl}_{[\bar{g}]}(X)= Y(S^n).
$$
Then there exist $k$ points $\{p_1,\ldots, p_k\}$ in $S^n$ such that
\begin{enumerate}
\item[{\bf (i)}] each manifold $(Z_j,h_j)$ is homothetic to
$(S^{n-1},h_+)$,
\item[{\bf (ii)}] the manifold $(X,[\bar{g}])$ is conformally
equivalent to the punctured sphere $(S^n\setminus \{p_1,\ldots,p_k\},
C_{\can})$.
\end{enumerate}
Here $h_j=h|_{Z_j}$ and $C_{\can}$ denotes the canonical conformal class on
$S^n$.
\end{Corollary}
\begin{lproof}{Proof of Corollary \ref{can:Th4}} 
We first notice that Proposition \ref{cyl:univ-bound} implies 
$$
Y^{\cyl}_{[\bar{g}]}(X)\leq Y^{\cyl}_{[h + dt^2]}(Z \times \R) =
\min_{1\leq j\leq k}
Y^{\cyl}_{[h_j+dt^2]}(Z_j\times \R),
$$
and hence
$$
Y^{\cyl}_{[h_j+dt^2]}(Z_j\times \R) = Y(S^n)
$$
for $j=1,\ldots,k$. Then Theorem \ref{can:Th2} implies that
$(Z_j,h_j)$ is homothetic to $(S^{n-1},h_+)$ for $j=1,\ldots,k$. By
the definition of cylindrical manifold modeled by $(S^{n-1},h_+)$,
there exists a smooth conformal compactification $(\hat{X}, \hat{C})$,
$\hat{X}= X\sqcup \{p_1,\ldots, p_k\}$, of the conformal manifold
$(X,[\bar{g}])$ such that $Y_{\hat{C}}(\hat{X})=
Y^{\cyl}_{[\bar{g}]}(X)$.  Here $Y_{\hat{C}}(\hat{X})$ stands for the
Yamabe constant of the closed conformal manifold $(\hat{X},
\hat{C})$.

By Schoen's Theorem \cite[Theorem 2]{Sc1} (cf. \cite{SY3} and
\cite{Au3}, \cite{LP}, \cite{SY4}), the equality
$Y_{\hat{C}}(\hat{X})=Y(S^n)$ implies that $(\hat{X}, \hat{C})$ is
conformally equivalent to $(S^n, C_{\can})$. Hence $(X,[\bar{g}])$ is
conformally equivalent to the punctured sphere $(S^n\setminus
\{p_1,\ldots,p_k\}, C_{\can})$.
\end{lproof}
For the proof of Theorem \ref{can:Th2}, we first prepare several
lemmas and propositions. Let $\Gamma$ be a spherical space form group
acting freely on $S^{n-1}(1)$ by isometries.  Let denote by
$S^{n-1}(1)/\Gamma$ the spherical space form, that is, the smooth
quotient of $S^{n-1}(1)$ by $\Gamma$.
\begin{Lemma}\label{can:L4}
Let $(Z\times \R, \bar{h}=h+dt^2)$ be a canonical cylindrical manifold
of $\dim (Z\times \R) =n \geq 3$ with $\lambda({\cal L}_h)>0$. If the
manifold $(Z\times \R, [\bar{h}])$ is locally conformally flat, then
$(Z,h)$ is homothetic to a smooth quotient $S^{n-1}(1)/\Gamma$ of
$S^{n-1}(1)$.
\end{Lemma}
\begin{Proof}
First we assume that $n\geq 4$. Let $W_{\bar{h}} =
(\bar{W}_{\alpha\beta\gamma\delta})$ denote the Weyl curvature tensor
of the metric $\bar{h}$, where $\alpha,\beta, \gamma, \delta =
0,1,\ldots,n-1$. Here the index $0$ corresponds to the $t$-coordinate, and
the indices $1,\ldots,n-1$ correspond to local coordinates
$x=(x^1,\ldots,x^{n-1})$ on $Z$ respectively. Then we  obtain that
\begin{equation}\label{can:eq12}
\begin{array}{l}
0 = \bar{W}_{0i0j} = 
-\frac{1}{n-2}\left(
\bar{R}_{ij} -\frac{R_{\bar{h}}}{n-1} \bar{h}_{ij}
\right)
= 
-\frac{1}{n-2}\left(
R_{ij} -\frac{R_h}{n-1} h_{ij}
\right)
\end{array}
\end{equation}
for all $1\leq i,j\leq n-1$. Here $\left(\bar{R}_{ij}\right)$ and
$\left(R_{ij}\right)$ denote respectively the Ricci curvature tensors
of the metrics $\bar{h}$ and $h$. Then (\ref{can:eq12}) combined with
$\lambda({\cal L}_h)>0$ implies that $h$ is an Einstein metric of
positive scalar curvature on $Z$. We use this to obtain that
$$
0 = \bar{W}_{ijk\ell} = R_{ijk\ell} -
\frac{R_h}{(n-1)(n-2)}(h_{ik}h_{j\ell} - h_{i\ell}h_{jk}) 
$$
for all $1\leq i,j,k,\ell\leq n-1$. The right hand side of the above
equation is nothing but the concircular curvature tensor of $h$. Thus
$h$ is a metric of positive constant curvature, and hence $(Z,h)$ is a
homothetic to a smooth quotient $S^{n-1}(1)/\Gamma$ of $S^{n-1}(1)$.

Next we consider the case $n=3$. Let
$$
B_{\bar{h}} = \left(
\bar{B}_{\alpha\beta\gamma}
\right)=
\left(
(\bar{\nabla}_{\gamma} \bar{R}_{\alpha\beta}-
\bar{\nabla}_{\beta} \bar{R}_{\alpha\gamma}) -\frac{1}{4}
(\p_{\gamma} R_{\bar{h}}\cdot \bar{h}_{\alpha\beta} - 
\p_{\beta} R_{\bar{h}}\cdot \bar{h}_{\alpha\gamma})
\right)
$$
denote the Bak tensor of $\bar{h}$. Then
$$
0 = \bar{B}_{00i} = -\frac{1}{4}\p_i R_{\bar{h}} = -\frac{1}{4}\p_i R_{h} 
$$
for $i=1,2$. This identity combined with the assumption $\lambda({\cal
L}_h)>0$ and $\dim Z =2$ implies that $h$ is a metric of positive
constant Gaussian curvature on $Z$. Therefore, $(Z,h)$ is homothetic
to either $S^2(1)$ or the projective space $\R\P^2= S^2(1)/\Z_2$.
\end{Proof}
\begin{Proposition}\label{can:P5}
Under the same assumptions as in Theorem \ref{can:Th2}, we also assume
that $(Z\times \R, \bar{h}=h+dt^2)$ is locally conformally flat. Then 
$(Z,h)$ is homothetic to the standard sphere $S^{n-1}(1)$.
\end{Proposition}
\begin{Proof}
The condition $Y^{\cyl}_{[\bar{h}]}( Z\times \R) = Y(S^n) >0 $ implies
that $\lambda({\cal L}_h)>0$, and hence from Lemma \ref{can:L4}
$(Z,h)$ is homothetic to a smooth quotient $S^{n-1}(1)/\Gamma$ of
$S^{n-1}(1)$.

Let $h_{\Gamma}$ denote the metric on $S^{n-1}(1)/\Gamma$ induced by
the metric $h_+$, i.e. $(S^{n-1}/\Gamma,
h_{\Gamma})=S^{n-1}(1)/\Gamma$. Then we notice that the first
eigenvalues $\lambda({\cal L}_{h_+})$ and $\lambda({\cal
L}_{h_{\Gamma}})$ are equal.

Now we use the above remark and Proposition \ref{cyl:univ-bound-2} to
get the estimate
$$
Y^{\cyl}_{[\bar{g}]}(Z\times \R) \leq Y(S^n)\cdot
\left(
\frac{\Vol(S^{n-1}(1)/\Gamma)}{\Vol(S^{n-1}(1))}
\right)^{\frac{2}{n}} =  \frac{Y(S^n)}{|\Gamma|^{\frac{2}{n}}}\ ,
$$
where $|\Gamma|$ denotes the order of $\Gamma$. If $|\Gamma| \geq 2$,
then $Y^{\cyl}_{[\bar{g}]}(Z\times \R)< Y(S^n)$, which contradicts the
assumption. Therefore $|\Gamma|=1$, and hence $(Z,h)$ is homothetic to
$S^{n-1}(1)$.
\end{Proof}
{\bf Remark.} The above argument implies that 
$Y^{\cyl}_{[h_{\Gamma} + dt^2]}((S^{n-1}/\Gamma) \times \R) 
\leq Y(S^n)/|\Gamma|^{\frac{2}{n}}$ for any spherical space form 
group $\Gamma$. 
Moreover, combining this inequality with 
$Y^{\cyl}_{[h_+ + dt^2]}(S^{n-1} \times \R) = Y(S^n)$, 
we then obtain that 
$$
Y^{\cyl}_{[h_{\Gamma} + dt^2]}((S^{n-1}/\Gamma) \times \R) 
= \frac{Y(S^n)}{|\Gamma|^{\frac{2}{n}}}. 
$$

\noindent
When $\dim (Z\times \R)=n\geq 6$, the following result combined with
Proposition \ref{can:P5} completes the proof of Theorem \ref{can:Th2}.
\begin{Proposition}\label{can:P6}
Let $(Z\times \R, \bar{h}=h+dt^2)$ be a canonical cylindrical manifold
of $\dim (Z\times \R) =n \geq 6$ with $\lambda({\cal L}_h)>0$. Assume
that the manifold $(Z\times \R, [\bar{h}])$ is not locally conformally
flat. Then
$$
Y^{\cyl}_{[\bar{h}]}(Z\times \R)< Y(S^n)=Y_{[h_++dt^2]}^{\cyl}(S^{n-1}\times
\R).
$$
\end{Proposition}
For the proof of Proposition \ref{can:P6}, we first recall the
following useful result (see \cite[Theorem 5.1]{LP} and \cite{Ca},
\cite{Gu}).
\begin{Lemma}\label{can:L7}{\rm (Conformal normal coordinates)}
Let $(M,\bar{g})$ be a Riemannian manifold and $p$ a point in
$M$. Then there exists a conformal metric $g\in [\bar{g}]$ on $M$ such
that
$$
d\sigma_{g} = dx := dx^1\wedge\cdots\wedge dx^n , \ \ \ \mbox{i.e.} \ \ 
\det(g_{ij}(x))=1
$$
in a $g$-normal coordinate neighborhood. Here $x=(x^1,\ldots,x^n)$
denotes $g$-normal coordinates at $p$. 
\end{Lemma}
\begin{lproof}{Proof of Proposition \ref{can:P6}} 
Since $(Z\times \R, \bar{h}=h+dt^2)$ is not locally conformally flat,
there exists a point $p=(x_0,0)\in Z\times \R$ such that
$|W_{\bar{h}}(p)|\neq 0$.  Let $g\in [\bar{h}]$ be a conformal metric
on $Z\times \R$ such that
$$
d\sigma_g = dx \ \ \ \mbox{on} \ \ \
B_{\rho_0}(0)=\{|x|<\rho_0\},
$$
where $x=(x^1,\ldots,x^n)$ are $g$-normal coordinates at
$p$. Similarly to Aubin's argument in \cite{Au2} (cf. \cite[Theorem
B]{LP}), we construct a family of test functions as follows. Let
$\{u_{\epsilon}\}_{\epsilon>0}$ be the family of positive functions on
$\R^n$ given by
$$
u_{\epsilon}(x) =\left(
\frac{\epsilon}{\epsilon^2+|x|^2}
\right)^{\frac{n-2}{2}} \ \ \ \mbox{for} \ \ \ x\in \R^n,
$$
i.e. $\{u_{\epsilon}\}_{\epsilon>0}$ are the \emph{instantons} on $(
\R^n\cong S^n\setminus \{{\mathrm p}{\mathrm o}{\mathrm i}{\mathrm
n}{\mathrm t}\}, C_{\can})$ centered at $x=0$. Let $\rho>0$ be a small
constant with $2\rho<\rho_0$, and $\eta$ a smooth cut-off function of
$r=|x|$ which satisfies
$$
\{
\begin{array}{ll}
\eta(x)=1 &\mbox{for} \ \ |x|\leq \rho,
\\
\eta(x)=0 &\mbox{for} \ \ |x|\geq 2\rho,
\\
|\nabla\eta|\leq \frac{2}{\rho} &\mbox{for} \ \ \rho\leq |x|\leq 2\rho.
\end{array}
\right.
$$
Set $\psi_{\epsilon}=\eta\cdot u_{\epsilon}$ on $Z\times \R$ for
$\epsilon>0$. Then we have the following estimate (see \cite[Proof of
Theorem B]{LP} for details):
$$
\begin{array}{l}
\displaystyle
Q_{(Z\times \R,g)}(\psi_{\epsilon})\leq
\\
\\
\{
\begin{array}{ll}
\displaystyle
Y(S^6) + \frac{1}{\|\psi_{\epsilon}\|^2_{L_{g}^3(Z\times \R)}}
\left[\begin{picture}(0,20)\end{picture}
-c_6\cdot |W_{g}(p)|^2 \epsilon^4\log(1/\epsilon) + O(\epsilon^4)
\right]
& \!\!\!\mbox{if} \ \ n=6,
\\
\\
\displaystyle
Y(S^n) + \frac{1}{\|\psi_{\epsilon}\|^2_{L_{g}^{\frac{2n}{n-2}}(Z\times \R)}}
\left[\begin{picture}(0,20)\end{picture}
- c_n\cdot |W_{g}(p)|^2\epsilon^4 + o(\epsilon^4)
\right]
& \!\!\!\mbox{if} \ \ n\geq 7,
\end{array}
\right.
\end{array}
$$
and $0< \|\psi_{\epsilon}\|^2_{L_g^{\frac{2n}{n-2}}(Z\times \R)}\leq
K$ for any small $\varepsilon > 0$. Here $c_6$, $c_n$ and $K$ denote
positive constants depending only on $n$. Recall that
$|W_{g}(p)|\neq 0$. Choosing $\epsilon>0$ sufficiently small, we
obtain that
\begin{equation}\label{can:eq12a}
Q_{(Z\times \R,g)}(\psi_{\epsilon}) < Y(S^n). 
\end{equation}
Then it follows from (\ref{can:eq12a}) and Fact \ref{cyl:Th1} that
$$
\begin{array}{rcl}
Y^{\cyl}_{[\bar{h}]}(Z\times \R) & = &\displaystyle
\!\!\!\!\! \inf_{\begin{array}{c}
^{f\in C_c^{\infty}(Z\times \R)}
\\
^{f\not\equiv 0}
\end{array}} \!\!\!\!\!\! Q_{(Z\times \R, \bar{h})}(f)
\\
\\
& = & \displaystyle 
\!\!\!\!\! \inf_{\begin{array}{c}
^{f\in C_c^{\infty}(Z\times \R)}
\\
^{f\not\equiv 0}
\end{array}} \!\!\!\!\!\! Q_{(Z\times \R, g)}(f)
\leq Q_{(Z\times \R, g)}(\psi_{\epsilon}) < Y(S^n).
\end{array}
$$
This completes the proof.
\end{lproof}
Now we consider the remaining case when $n=3,4,5$ in Theorem
\ref{can:Th2}. We start with the existence and uniqueness properties of minimal Green's
functions of the conformal Laplacians in dimensions $n\geq 3$.
\begin{Lemma}\label{can:L8} 
Let $(Z\times \R, \bar{h}=h+dt^2)$ be a canonical cylindrical manifold
of $\dim (Z\times \R) = n \geq 3$ with $\lambda({\cal L}_h)>0$. Let
$p\in Z\times \R$ be an arbitrary point. Then there exists a unique
normalized minimal positive Green's function $\bar{G}_p$ (with pole at
$p$) of the conformal Laplacian $\L_{\bar{h}}$ on $Z\times
\R$. Namely, for each point $p\in Z\times \R$ the Green's function
$\bar{G}_p\in C^{\infty}((Z\times \R)\setminus \{p\})$ satisfies the
following properties:
\begin{enumerate}
\item[{\bf (i)}] $\L_{\bar{h}} \bar{G}_p = \beta_n\cdot\delta_p$,
where $\beta_n=4(n-1)\Vol(S^{n-2}(1))>0$, and $\delta_p$ is the Dirac
$\delta$-function at $p$.
\item[{\bf (ii)}] $\bar{G}_p > 0$ on $(Z\times \R)\setminus \{p\}$.
\item[{\bf (iii)}] If $\bar{G}^{\prime}_p\in
C^{\infty}((Z\times \R)\setminus \{p\})$ is a normalized positive
Green's function (i.e. $\bar{G}^{\prime}_p$ satisfies the conditions
{\bf (i)} and {\bf (ii)}), then $\bar{G}^{\prime}_p \geq \bar{G}_p $.
\end{enumerate}
\end{Lemma}
\begin{Proof} Let $G_p$ be any normalized positive Green's function 
with pole at $p$ of $\L_{\bar{h}}$. We notice that for $u\in
C^{\infty}_+(Z\times \R)$ the function $u(p)^{\frac{n+2}{n-2}}\cdot
u^{-1} \cdot G_p$ is also such a normalized positive Green's function
of the conformal Laplacian $\L_{\tilde{h}}$, with respect to the
metric $\tilde{h}:=u^{\frac{4}{n-2}}\cdot \bar{h}$. Using this
observation combined with $\lambda({\cal L}_h)>0$, we assume that
$R_h>0$ on $Z$ (as in Convention \ref{convention}). We may also assume
that $p=(x_0,0)\in Z\times \R$. Then, since $R_h>0$, for each $i\geq
1$ there exists a unique positive Green's function 
$$ 
\begin{array}{c}
G_p^{(i)}\in
C_+^{\infty}((Z\times (-i,i))\setminus \{p\})\cap C^0((Z\times [-i,
i])\setminus \{p\}) \ \ \mbox{such that}
\\
\\
\{
\begin{array}{ll}
\L_{\bar{h}} G_p^{(i)} = \beta_n\cdot \delta_p & \mbox{on} \ \ Z\times (-i,i),
\\
G_p^{(i)} = 0  & \mbox{on} \ \ Z\times \{\pm i\}.
\end{array}
\right.
\end{array}
$$
For each $ G_p^{(i)}$ we use the normalization condition $\L_{\bar{h}}
G_p^{(i)} = \beta_n\cdot \delta_p $ to obtain the following expansion
in a fixed normal coordinate system $x$ centered at the point $p$:
$$
G_p^{(i)}(x) = |x|^{2-n}\left(1+ o(1)\right) \ \ \ \mbox{as} \ \
|x|\to 0.
$$
Then by the maximum principle, for any $i\geq 1$,
$$
G_p^{(i)} \leq G_p^{(i+1)} \leq G_p^{(i+2)}\leq \cdots \ \ \ \mbox{on}
\ \ \ Z\times [-i,i].
$$
We set $S^1=[-1,1]/_{\{-1\}\sim \{1\}}$. Let $N$ denote the closed
manifold $Z\times S^1$, and let $g$ be the metric on $N$ induced by
$\bar{h}=h+dt^2$ via the Riemannian submersion $\Phi: (Z\times \R,
\bar{h}) \to (Z\times S^1,g)$. Then the condition $\lambda({\cal
L}_h)>0$ implies that the first eigenvalue $\lambda(\L_{g})>0$. Hence
there exists a unique normalized minimal positive Green's function
$G_{\check{p}}$ of $\L_g$ on $N$ with pole at $\check{p}$, where
$\Phi(p) =\check{p}$. The maximum principle also implies that
$$
G_p^{(i)}\leq \Phi^*G_{\check{p}} \ \ \ \ \mbox{on} \ \ Z\times [-i,i] 
$$
for any $i\geq 1$. Then by Harnack's convergence theorem, there exists
a normalized positive Green's function $\bar{G}_p$ with pole at $p$ of
$\L_{\bar{h}}$ on $Z\times \R$ such that the sequence $\{
\bar{G}_p^{(i)}\}$ converges uniformly to $\bar{G}_p$ on each cylinder
$$
\begin{array}{c}
\left(Z\times [-j, j]\right)\setminus B_{1/j}(p;\bar{h}),  \ \mbox{where} \ \ 
\\
\\
B_{1/j}(p;\bar{h}) = \{ x\in Z\times \R \ | \
\dist_{\bar{h}}(p,x)<1/j\}.
\end{array}
$$
Now, by the construction of the Green's function $\bar{G}_p$ and the
maximum principle, $\bar{G}_p$ is a normalized minimal
positive Green's function of $\L_{\bar{h}}$ with pole at $p$.
The uniqueness for such $\bar{G}_p$ follows directly from 
the minimality condition (iii). 
\end{Proof}
{\bf Remark.} Let $\iota: Z\times \R \to Z\times \R$ denote the
involution defined as $\iota(x,t)=(x,-t)$. Let $o=(x_0,0)\in Z\times
\R $ be a point. Since $\iota^*\bar{h}= \bar{h}$, the uniqueness of
a normalized minimal positive Green's function implies that
$\bar{G}_o$ is $\iota$-invariant.
\vspace{2mm}

\noindent
Similarly to Theorem \ref{yam:Th2}, we obtain the following.
\begin{Lemma}\label{can:L10}
There exist positive constants $a\leq b$, $C_1\leq C_2$ such that
\begin{equation}\label{can:eq13}
C_1\cdot e^{-b|t|} \leq \bar{G}_o(x,t) \leq C_2 \cdot e^{-a|t|} 
\ \ \  \mbox{on} \ \ Z\times((-\infty, -1]\sqcup [1,\infty)). \ \ 
\end{equation}
\end{Lemma}
To proceed further, we set $S^1_{\ell} =
[-2^{\ell},2^{\ell}]/_{\{-2^{\ell}\}\sim \{2^{\ell}\}}$ for
$\ell=1,2,\ldots$. Let $N_{\ell}$ denote the closed manifold $Z\times
S^1_{\ell}$, and let $g_{\ell}$ be the metric on $N_{\ell}$ defined
via the Riemannian submersion $\Phi_{\ell}: (Z\times \R ,\bar{h}) \to
(N_{\ell},g_{\ell})$. Let $G_{\ell}$ be the normalized minimal
positive Green's function of $\L_{g_{\ell}}$ with pole at
$\check{o}=\Phi_{\ell}(o)\in N_{\ell}$. Then we have the following
fact (see \cite[Chapter 6, Proposition 2.4]{SY4} for the proof):
\begin{Lemma}\label{can:L11}
Let $\Gamma_{\ell}\subset \pi_1(N_{\ell})$ denote  the deck transformation
group of the covering $\Phi_{\ell} : Z\times \R \to N_{\ell}$. Then
\begin{equation}\label{can:eq14}
\Phi_{\ell}^*G_{\ell} = \sum_{\gamma\in \Gamma_{\ell}}
 \gamma_*\bar{G}_o \ \ \ \mbox{on} \ \ \ (Z\times \R)\setminus \{o\}.
\end{equation}
\end{Lemma}
From now on, we assume that $ \dim(Z \times \R) = n = 3, 4, 5$. Let $o
= (x_0,0)\in Z\times \R$ be a fixed point and $\bar{G}=\bar{G}_o$ the
normalized minimal positive Green's function given in Lemma \ref{can:L8}.
Below we suppress the dependence of $\bar{G} = \bar{G}_{o}$ on the
point $o$.

\begin{figure}[h]

\hspace*{10mm}
\PSbox{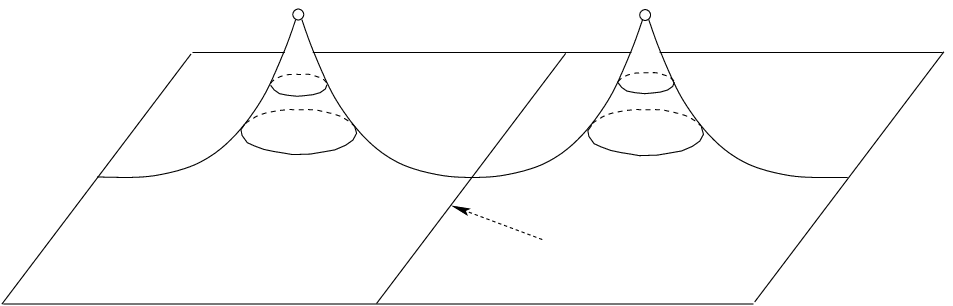}{10mm}{32mm}  
\begin{picture}(0,1)
\put(100,10){{\small $(Z\setminus \{x_0\})\times \{0\}$}}
\put(135,33){{\small $Z \times \{L\}$}}
\put(28,33){{\small $Z \times \{-L\}$}}
\end{picture}
\caption{The manifold $(\hat{X},\hat{h})$.}\label{fig6-2}
\end{figure}

We set $\hat{X}= (Z\times \R)\setminus \{o\}$ with the metric
$\hat{h}= \bar{G}^{\frac{4}{n-2}}\cdot h$. Then the condition
$\L_{\bar{h}}\bar{G} =\beta_n \delta_o$ implies that
\begin{equation}\label{can:eq15}
R_{\hat{h}}\equiv 0 \ \ \ \mbox{on} \ \ \hat{X}.
\end{equation}
Let $x=(x^1,\ldots,x^n)$ be conformal normal coordinates at
$o=(x_0,0)\in Z\times \R$ satisfying 
$$
d\sigma_{\bar{h}}(x) = dx .
$$
Then we use Lemma 6.4 and Theorem 6.5 in \cite{LP} to obtain the
following result. First, we need some notations. We denote by
$\bar{\nabla}$ for the covariant derivative with respect to the metric
$\bar{h}$. Then we write $f=O^{\prime\prime}(r^k)$ to mean that
$f=O(r^k)$, $\bar{\nabla}f=O(r^{k-1})$ and
$\bar{\nabla}^2f=O(r^{k-2})$ (following the notations from \cite{LP}).
\begin{Lemma}\label{can:L12}
In the conformal normal coordinates $x=(x^i)$, the Green's function
$\bar{G}$ has the following expansion:
\begin{equation}\label{can:eq17}
\bar{G}(x) = r^{2-n} + A + O^{\prime\prime}(r), \ \ \ A\equiv \const.
\end{equation}
Furthermore, in the inverted conformal normal coordinates $y=(y^i =
r^{-2}x^i)$, the metric $\hat{h}$ has the following expansion:
\begin{equation}\label{can:eq18}
\{
\begin{array}{rcl}
\hat{h}_{ij}(y) &=& \displaystyle
\gamma(y)^{\frac{4}{n-2}}\left(\delta_{ij} + O^{\prime\prime}(\rho^{-2})\right),
\\
\gamma(y)&=& \displaystyle
1 + A\rho^{2-n} + O^{\prime\prime}(\rho^{1-n}),
\end{array}
\right.
\end{equation}
where $\rho=|y|=r^{-1}$.
\end{Lemma}
A similar argument to the proofs of \cite[Theorem 1]{Sc1} and
\cite[Chapter 5, Theorem 4.1]{SY4} implies the following.
\begin{Proposition}\label{can:P13}
If $A>0$ in {\rm (\ref{can:eq17})}, then $Y^{\cyl}_{[\bar{h}]}(Z\times \R) <
Y(S^n)$. 
\end{Proposition}
Now we notice that (\ref{can:eq15}) and (\ref{can:eq18}) imply that
$(\hat{X},\hat{h})$ is a scalar-flat, asymptotically flat manifold of
order $1$ (if $n=3$), and of order 2 (if $n=4,5$). We emphasize that
$(\hat{X},\hat{h})$ has two almost conical singularities and $\hat{h}$
is $\iota$-invariant, see Fig. \ref{fig6-2}.
We notice that in this case the \emph{mass} ${\mathfrak
m}={\mathfrak m}(\hat{h},(y^i))$ is well-defined in the standard
way:
$$
{\mathfrak m}(\hat{h},(y^i)) = \lim_{R\to\infty}
\frac{1}{\Vol(S^{n-1})}\int _{\{|y|=R\}} \ \ \sum_{i,j=1}^n
\left( \p_i \hat{h}_{ij} - \p_j\hat{h}_{ii}\right)\p_j \lr \ dy.
$$
We emphasize that the mass ${\mathfrak m}(\hat{h})={\mathfrak
m}(\hat{h},(y^i))$ depends only on the metric $\hat{h}$, see
\cite[Theorem 4.2]{Bar}.
However, for any large constant $L>0$, the region $(X\times
 [-L,L])\cap \hat{X}$ is not convex with respect to the metric
 $\hat{h}$. Hence we cannot apply directly the technique (developed in
\cite{SY1}) for asymptotically flat manifolds with many ends to proving
that the mass ${\mathfrak m}\geq 0$.

We also remark that if $n=3$ or $\bar{h}$ is locally conformally
flat near $o$, then $\hat{h}$ has the following expansion:
$$
\hat{h}_{ij}(y) =\left(
1+\frac{4}{n-2} A\rho^{2-n}
\right)\delta_{ij} + O^{\prime\prime}(\rho^{1-n}).
$$
Then the mass ${\mathfrak m}(\hat{h})={\mathfrak m}(\hat{h},(y^i))$ is
given as
\begin{equation}\label{can:eq19}
{\mathfrak m}(\hat{h}) = 4(n-1) A .
\end{equation}
Even if $\bar{h}$ is not locally conformally flat on any open set
inside of $Z \times \R$, the equality (\ref{can:eq19}) still holds
(if $3\leq n\leq 5$), see \cite[Lemma 9.7]{LP}.

The following result is a version of the positive mass theorem for the
asymptotically flat manifold $(\hat{X},\hat{h})$ with almost conical
singularities.
\begin{Theorem}\label{can:Th14}{\rm (First part of the positive mass theorem)}
Let $(\hat{X},\hat{h})$ be the asymptotically flat manifold with two
almost conical singularities and $\dim\hat{X}=n=3,4,5$, as above.
Then
\begin{equation}\label{can:eq20}
A\geq 0.
\end{equation}
\end{Theorem}
\begin{Proof}
Let $\hat{N}_{\ell}= N_{\ell}\setminus \{\check{o}\}$, where
$\check{o}=\Phi_{\ell}(o)$, with the metric $\hat{g}_{\ell} =
G_{\ell}^{\frac{4}{n-2}}\cdot g_{\ell}$ for $\ell=1,2,\ldots$. Since
$\Phi_{\ell}$ is a local isometry near $o$ with respect to $\bar{h}$
and $g_{\ell}$, we can use the same conformal normal coordinates
$x=(x^i)$ near $\check{o}\in N_{\ell}$. Then each Green's function
$G_{\ell}$ has the following expansion:
\begin{equation}\label{can:eq21}
G_{\ell}(x) = r^{2-n} + A_{\ell} + O^{\prime\prime}(r), \ \ \
A_{\ell}\equiv\const.
\end{equation}
The mass of the asymptotically flat manifold
$(\hat{N}_{\ell},\hat{g}_{\ell})$ is also given by
$$
{\mathfrak m}(\hat{g}_{\ell})= 4(n-1) A_{\ell}
$$
It then follows from the positive mass theorem \cite{SY4,SY1,SY2}
(cf. \cite{Bar} and \cite{LP}) that
\begin{equation}\label{can:eq22}
A_{\ell} \geq 0 \ \ \mbox{for any $\ell \geq 1$}.
\end{equation}
Under the identification $\Phi_{\ell}|_{Z\times [-2^{\ell},2^{\ell}]}$
and using that $\iota^* G_{\ell}= G_{\ell}$, we may regard $G_{\ell}$
as a function on $Z\times [-2^{\ell},2^{\ell}]$. Then the maximum
principle combined with (\ref{can:eq14}) implies that
\begin{equation}\label{can:eq23}
G_{\ell}\geq G_{\ell+1}\geq \cdots\geq \bar{G} >0 \ \ \ \mbox{on}  \ \ 
(Z\times [-2^{\ell},2^{\ell}])\setminus \{o\}.
\end{equation}
Hence $\{G_{\ell}\}$ converges uniformly to $\bar{G}$ on each cylinder
$ (Z\times [-i,i])\setminus B_{o}(1/i;\bar{h})$. The estimate (\ref{can:eq13})
and the properties (\ref{can:eq14}), (\ref{can:eq23}) give that
$$
\bar{G}(x) \leq G_{\ell}(x) \leq \bar{G}(x) + 
 \frac{2C_2\cdot e^{-a\cdot 2^{\ell}}}{1 - e^{-a\cdot 2^{\ell}}}
$$
for $\ell>>1$. We then obtain that
\begin{equation}\label{can:eq24}
|\bar{G}(x) -G_{\ell}(x)| \leq C\cdot e^{-a \cdot 2^{\ell}} \ \ \
 \mbox{on} \ \ U_{o}\setminus \{o\}
\end{equation}
for $\ell>>1$. Here $U_{o}$ denotes a coordinate neighborhood cenered
at $o$ and $C>0$ is a constant independent of $\ell$ and $x$.  From
(\ref{can:eq17}) and (\ref{can:eq21}), the function
$\bar{G}(x)-G_{\ell}(x)$ is smooth on $U_o$. Thus by letting $|x|\to 0$ in
(\ref{can:eq24}), we obtain
$$
|A-A_{\ell}| \leq C\cdot e^{-a\cdot 2^{\ell}} \ \ \ \mbox{for} \ \ \ell>>1,
$$
and hence by letting $\ell\to\infty$, we prove that
\begin{equation}\label{can:eq25}
A=\lim_{\ell\to\infty} A_{\ell}.
\end{equation}
Therefore, it follows from (\ref{can:eq22}) and (\ref{can:eq25}) that
$A\geq 0$.
\end{Proof}
\begin{Proposition}\label{can:P15}
If the metric $\hat{h}= \bar{G}^{\frac{4}{n-2}}\cdot \bar{h}$ is not
Ricci-flat on $\hat{X}$, then $A>0$, where $A$ is the same constant as
in {\rm (\ref{can:eq17})}.
\end{Proposition}
\begin{Proof}
For any symmetric 2-tensor $S=(S_{ij})$ with compact support in
$\hat{X}=(Z\times \R)\setminus \{o\}$, we define the family of smooth
metrics
$$
\bar{h}^{s} = \bar{h}+ s\cdot \bar{G}^{-\frac{4}{n-2}}\cdot S \ \ \ \mbox{on} \ \ \
Z\times \R
$$
for small $s\geq 0$. Notice that $\bar{h}^{s}$ is no longer a product
metric on $Z\times \R$ for $s>0$. However, the arguments in Lemmas
\ref{can:L8}, \ref{can:L10}, \ref{can:L11}, \ref{can:L12} and Theorem
\ref{can:Th14} are still valid since $\supp(S)$ is compact in
$\hat{X}$. Hence there exists a small constant $\delta_0>0$ such that
the following holds:
\vspace{2mm}

\noindent
{\sl For any $s$ with $0\leq s \leq \delta_0$ there exists a unique
normalized minimal positive Green's function $\bar{G}^s$ with pole at
$o$ of $\L_{\bar{h}^s}$ on $\hat{X}$ such that $\bar{G}^s$ has the
expansion (as in {\rm (\ref{can:eq17})}):}
$$
\bar{G}^s(x) = |x|^{2-n} +A_s + O^{\prime\prime}(|x|), \ \ \ A_s\equiv\const\geq 0.
$$
Then we can apply the argument \cite[Lemma 3]{Sc1} to proving that $A>0$. 
\end{Proof}
From Propositions \ref{can:P13}, \ref{can:P15} and Theorem
\ref{can:Th14}, in order to complete the proof of Theorem
\ref{can:Th2}, it is enough to show the following.
\begin{Proposition}\label{can:P16}
Under the same assumptions as in Theorem \ref{can:Th14}, we also
assume that $\hat{h}$ is Ricci-flat
\begin{equation}\label{can:eq82}
\Ric_{\hat{h}}\equiv 0 \ \ \ \mbox{on} \ \ \hat{X} = (Z\times \R)\setminus \{o\}.
\end{equation}
Moreover, assume that $A=0$ in {\rm (\ref{can:eq20})} when $n=3$. Then
$(\hat{X},\hat{h})$ is isometric to $\R^n\setminus\{\mbox{2 points}\}$
with the Euclidean metric. In particular, $(Z,h)$ is homothetic to the
sphere $S^{n-1}(1)$.
\end{Proposition}
As a complementary assertion to Theorem
\ref{can:Th14}, we also obtain the following.
\begin{Theorem}\label{can:Cor17}{\rm (Second part of the 
positive mass theorem)} Let $(\hat{X},\hat{h})$ be an asymptotically
flat manifold with two almost conical singularities and 
$\dim\hat{X}=n=3,4,5$, as above.  Then if $A=0$, then
$(\hat{X},\hat{h})$ is isometric to $\R^n\setminus\{2 \ {\mathrm
p}{\mathrm o}{\mathrm i}{\mathrm n}{\mathrm t}{\mathrm s}\}$ with the
Euclidian metric.
\end{Theorem}
\begin{lproof}{Proof of Theorem \ref{can:Cor17}}
The condition $A=0$ combined with Theorem \ref{can:Th14} and
Proposition \ref{can:P15} implies that $\hat{h}$ is Ricci-flat on
$\hat{X}$. It then follows from Proposition \ref{can:P16} that
$(\hat{X},\hat{h})$ is isometric to $\R^n\setminus\{\mbox{2
points}\}$. 
\end{lproof}
\begin{lproof}{Proof of Proposition \ref{can:P16}}
Recall that, in the inverted conformal normal coordinates
$y=(y^1,\ldots,y^n)$, the metric $\hat{h}$ has the following expansion
($n=3$ with $A=0$ or $n=4,5$):
$$
\hat{h}_{ij}(y)= \delta_{ij} + O^{\prime\prime}(\rho^{-2}) \ \ \ \mbox{as} \ \ 
\ \rho=|y| \to \infty.
$$
We choose a sufficiently large $L_0> 0$ and fix it. Set
$\hat{X}_0=\hat{X}\cap (Z\times [-L_0,L_0])$. 

Let $L_{\delta}^{k,p}(\hat{X}_0)$ denote the weighted Sobolev space
with weight $\delta\in \R$ on $(\hat{X}_0,\hat{h}; (y^i))$, $k=0,1,2$,
$1\leq p\leq \infty$, defined in \cite{Bar}. We need the
following result to complete the proof of Proposition \ref{can:P16}.
\begin{Lemma}\label{can:L18}
{\rm (Harmonic coordinates near infinity)} Let $\varepsilon$ {\rm ($0 <
\varepsilon < 1$)}, $q > n$ and $\tau$ {\rm ($\frac{3}{2} \leq \tau < 2$)} be
positive constants satisfying $q>\frac{n}{1-\epsilon}$. Then there
exist smooth functions $z^i\in C^{\infty}(\hat{X}_0)$ for $i=1,\ldots,n$
such that
\begin{equation}\label{can:eq83}
\{
\begin{array}{l}
\Delta_{\hat{h}} z^i = 0 \ \ \mbox{on} \ \ \hat{X}_0,
\\
\frac{\p z^i}{\p \nu} = 0 \ \ \ \mbox{on} \ \ \p\hat{X}_0,
\\
y^i-z^i \in L_{1-\tau}^{2,q}(\hat{X}_0),
\\
y^i - z^i = O^{\prime\prime}(\rho^{-(1-\epsilon)}) \ \ \mbox{as} \ \
\rho=|y|\to\infty.
\end{array}
\right.
\end{equation}
Here $\nu$ is the outward unit vector field normal to the boundary
$\p\hat{X}_0$ with respect to the metric $\hat{h}$. In particular,
$z=(z^i)$ are harmonic coordinates near infinity.
\end{Lemma}
\begin{lproof}{Proof of Lemma \ref{can:L18}}
We extend each function $y^i$ to a smooth function on $ \hat{X}_0$
satisfying $\frac{\p y^i}{\p\nu}=0$ on $\p\hat{X}_0$. We notice that
$$
\Delta_{\hat{h}} y^i = \hat{h}^{jk} \hat{\Gamma}^i_{jk} = O^{\prime}(\rho^{-3})
\ \ \ \mbox{on} \ \ \ \hat{X}_0.
$$
Here we write $f=O^{\prime}(\rho^{k})$ if $f=O(\rho^{k})$ and
$\hat{\nabla} f=O(\rho^{k-1})$.

Modifying the $L^{\frac{2n}{n-2}}$-theory in \cite[Lemmas 3.1,
3.2]{SY1} and using the $L^q$-estimates in the linear elliptic theory,
we can show the following:
\begin{Claim}
There exist unique smooth functions $u^i\in L^q(\hat{X}_0)$ such that
$$
\{
\begin{array}{l}
\Delta_{\hat{h}} u^i = \Delta_{\hat{h}} y^i = O^{\prime}(\rho^{-3}) \ \ \mbox{on} \ \ 
\hat{X}_0,
\\
\frac{\p u^i}{\p\nu}=0 \ \ \ \mbox{on} \ \ \p\hat{X}_0,
\\
u^i = O^{\prime\prime}(\rho^{-(1-\epsilon)}) \ \ \ \mbox{as} \ \ \rho\to\infty.
\end{array}
\right.
$$
\end{Claim}
Then the weighted Sobolev estimate \cite[Proposition 1.6]{Bar} and
the scale-broken estimate \cite[Theorem 1.10]{Bar} imply that
$u^i\in L^{2,q}_{1-\tau}(\hat{X}_0)$. Set $z^i= y^i-u^i\in
C^{\infty}(\hat{X}_0)$. Then the functions $z=(z^1,\ldots,z^n)$
satisfy the properties (\ref{can:eq83}). This completes the proof of
Lemma \ref{can:L18}.
\end{lproof}
We return to the proof of Proposition \ref{can:P16}. Let
$z=(z^1,\ldots,z^n)$ be the harmonic coordinates near infinity of
$\hat{X}_0$ obtained in Lemma \ref{can:L18}. Then the Ricci-flatness of
$\hat{h}$ combined with \cite[Proposition 3.3]{Bar} implies that,
in the harmonic coordinates $z=(z^i)$, the metric
$\hat{h}(z)=(\hat{h}_{ij}(z))$ satisfies
\begin{equation}\label{can:eq84}
\hat{h}_{ij}(z) -\delta_{ij} \in L^{2,q}_{-\eta}(\hat{X}_0), \ \ \
i,j=1,\ldots,n
\end{equation}
for any $\eta>n-2$. It then follows from (\ref{can:eq82}), (\ref{can:eq84}) and
\cite[Theorem 4.3]{Bar} that the mass ${\mathfrak m}(\hat{h})
={\mathfrak m}(\hat{h}, (z^i))=0$. 

Now we can apply the Bochner technique to complete the proof. The
harmonicity of $(z^i)$ implies that $\{d z^i\}$ are harmonic $1$-forms
on $(\hat{X}_0,\hat{h})$. From the Bochner formula for $1$-forms
$(z^i)$ and the condition that $\frac{\p z^i}{\p \nu}=0$ on
$\p\hat{X}_0$, we obtain that
\begin{equation}\label{can:eq85}
\begin{array}{l}
\displaystyle
\sum_{i=1}^n\int_{\hat{X}_0\cap\{|z|<R\}}\left[
|\hat{\nabla}dz^i|^2 + \Ric_{\hat{h}}(dz^i,dz^i)
\right]d\sigma_{\hat{h}} =
\\
\\
\displaystyle
\sum_{i,j=1}^n\int_{\{|z|=R\}} \hat{h}^{jk}\< dz^i, \hat{\nabla}_j
dz^i\> \p_k \lr \ d\sigma_{\hat{h}}
\end{array}
\end{equation}
for any $R>>1$. Letting $R\to\infty$ in (\ref{can:eq85}) and using the
Ricci-flatness $\Ric_{\hat{h}}\equiv 0$ on $\hat{X}$, then we have
\begin{equation}\label{can:eq86}
{\mathfrak m}(\hat{h})=\frac{1}{\Vol(S^{n-1}(1))}
\sum_{i=1}^n\int_{\hat{X}_0} |\hat{\nabla} dz^i|^2 d\sigma_{\hat{h}},
\end{equation}
see \cite[Theorem 4.4]{Bar} or \cite[Proposition 10.2]{LP}.  Now
we use that ${\mathfrak m}(\hat{h})=0$ in (\ref{can:eq86}). Then we
obtain that the $1$-forms $\{d z^i\}$ are parallel on $\hat{X}_0$ with
respect to $\hat{h}$. Since the coframe $\{d z^i\}$ is orthonormal at
infinity, then $\{d z^i\}$ is a parallel orthonormal coframe
everywhere on $\hat{X}_0$. This implies that the map
$z=(z^1,\ldots,z^n): (\hat{X}_0,\hat{h}) \to \R^n$ is a local
isometry, and hence $\bar{h}$ is locally conformally flat on
$\hat{X}$. From Lemma \ref{can:L4}, $(Z,h)$ is homothetic to a smooth
quotient $S^{n-1}(1)/\Gamma$ of $S^{n-1}(1)$.
\begin{Claim}\label{can:C19}
The manifold $(Z,h)$ is homothetic to $S^{n-1}(1)$.
\end{Claim}
\begin{lproof}{Proof of Claim \ref{can:C19}}
Let $\iota : Z\times \R \to Z\times \R $ be the involution given by
$\iota(x,t)=(x,-t)$.  Recall that $\bar{G}$ is $\iota$-invariant, and
hence the metric $\hat{h}= \bar{G}^{\frac{4}{n-2}}\cdot \bar{h}$ is an
$\iota$-invariant metric.  This implies that $(Z\setminus
\{x_0\})\times \{0\}$ is a totally geodesic submanifold of $(\hat{X},
\hat{h})$. Using this fact and that $z$ is a local isometry, we
conclude that the restriction $z$ to $(Z\setminus \{x_0\})\times
\{0\}$ is a global isometry onto a hyperplane in $\R^n$. Hence $(Z,h)$
is homothetic to $S^{n-1}(1)$.
\end{lproof}
From Claim \ref{can:C19}, the manifold $(\hat{X},\hat{h})$, $\hat{h}=
\bar{G}^{\frac{4}{n-2}}\cdot \bar{h}$, is conformally equivalent to
$\R^n\setminus \{\mbox{2 points}\}$ with the Euclidean metric $g_0$.
Under the identification $\hat{X}\cong \R^n\setminus \{\mbox{2
points}\}$, the metric $g_0$ is represented as $ g_0 =
G^{\frac{4}{n-2}}\cdot \bar{h}$ on $\hat{X}$. Here $G$ is a normalized
positive Green's function with pole at $o$ of $\L_{\bar{h}}$, which
satisfies $G(x,t)=o(1)$ as $|t|\to\infty$ for $(x,t)\in Z\times \R$.
Using the maximum principle combined with the minimality of $\bar{G}$
and the normalization for $\bar{G}$ and $G$, we obtain that
$\bar{G}=G$ on $\hat{X}$. Hence $\hat{h}= g_0$ on $\hat{X}$.
This completes the proof of Proposition \ref{can:P16}.
\end{lproof}
\section{Appendix}
Let $(X, \bar{g})$ be a cylindrical manifold of dim$X = n \geq 3$
 with tame ends $Z \times [0, \infty)$ and
 $\p_{\infty}\bar{g} = h \in \Riem(Z)$.  Here we study some properties
 of the conformal Laplacian $\L_{\bar{g}}$, and its relationship to the
 constant $Y^{\cyl}_{[\bar{g}]}(X)$ and the sign of the scalar curvature
 $R_{\hat{g}}$ of a conformal metric $\hat{g} \in [\bar{g}]$.  In the case
 of a closed conformal manifold $(M, C)$, the signs of the first eigenvalues
 of $\L_g$ and $\L_{\tilde{g}}$ for $g, \tilde{g} \in C$ are identical.
Since $(X, \bar{g})$ is not compact,
 $\L_{\bar{g}}$ has no first eigenvalue in general.
Instead, we consider the \emph{bottom of the spectrum of
$\L_{\bar{g}}$}, defined as
$$
\begin{array}{c}
\lambda(\L_{\bar{g}}) := \inf_{f\in L^{1,2}_{\bar{g}}(X), \
f\not\equiv 0} \frac{E_{(X,\bar{g})}(f)}{\|f\|^2_{L^{2}_{\bar{g}}}}, \
\ \mbox{where} 
\\
\\
E_{(X,\bar{g})}(f) = \int_X\[\alpha_n|df|^2 +
R_{\bar{g}}f^2\] d\sigma_{\bar{g}}.
\end{array}
$$
It is easy to show that $\lambda(\L_{\bar{g}}) > -\infty$ since
$\bar{g}= h+dt^2$ on $Z\times [1,\infty)$. Indeed, one has the
following estimate
$$
\lambda(\L_{\bar{g}}) \geq \min_{X} R_{\bar{g}} = \min \{\min_{X(1)}
R_{\bar{g}}, \min_{Z} R_h \} >-\infty .
$$
Similarly to the case of closed manifolds, the sign of
$\lambda(\L_{\bar{g}})$ is uniquely determined by the conformal
class $[\bar{g}]$ provided that $\bar{g}$ is a cylindrical metric.
\begin{Proposition}\label{app:P1}
Let $\bar{C}$ be a cylindrical conformal class on $X$ and $\bar{g},
\tilde{g} \in \bar{C}\cap \Riem^{\cyl}(X)$ two cylindrical metrics.
Then $\lambda(\L_{\bar{g}})$ and $\lambda(\L_{\tilde{g}})$ have the
same sign or are both zero. 
\end{Proposition}
\begin{Proof} The proof of the following statement is straightforward. 
\begin{Fact}\label{app:L1} 
Let $\tilde{g}=\phi^{\frac{4}{n-2}}\cdot \bar{g}$ (here $\tilde{g}$
and $\bar{g}$ are not necessarily cylindrical), where $\phi\in
C_+^{\infty}(X)$. Then $\L_{\tilde{g}} = \phi^{-\frac{n+2}{n-2}}\cdot
L_{\bar{g}}\circ \phi$, and $ E_{(X,\bar{g})}(\phi f) =
E_{(X,\tilde{g})}(f)$ for any function $f\in C_c^{\infty}(X)$.
\end{Fact}
Clearly we have $d\sigma_{\tilde{g}}= \phi^{\frac{2n}{n-2}}
d\sigma_{\bar{g}}$. This implies
$$
\left(
\inf_X \phi
\right)^{\frac{4}{n-2}} \|\varphi  f\|^2_{L_{\bar{g}}^2} \leq  \| f\|^2_{L_{\tilde{g}}^2}
\leq \left(
\sup_X \phi
\right)^{\frac{4}{n-2}} \|\varphi  f\|^2_{L_{\bar{g}}^2} 
$$
for any $f\in C^{\infty}_c(X)$. We proved in
Proposition \ref{conf:new} that for cylindrical metrics $\bar{g}$
and $\tilde{g}$ with $\tilde{g} = \varphi^{\frac{4}{n-2}}\cdot
\bar{g}$ there exists a constant $K\geq 1$ such that $0< K^{-1}\leq
\phi \leq K <\infty$ on $X$. This implies the estimate
\begin{equation}\label{app:eq2}
K^{-\frac{4}{n-2}} \cdot \frac{E_{(X,\bar{g})}(\phi f)}{\|\phi
f\|^2_{L^{2}_{\bar{g}}}} \leq
\frac{E_{(X,\tilde{g})}(f)}{\|f\|^2_{L^{2}_{\tilde{g}}}}\leq
K^{\frac{4}{n-2}} \cdot \frac{E_{(X,\bar{g})}(\phi f)}{\|\phi
f\|^2_{L^{2}_{\bar{g}}}}
\end{equation}
for any $f \in C^{\infty}_c(X)$.  It then follows from (\ref{app:eq2})
that
$$
K^{-\frac{4}{n-2}} \lambda(\L_{\bar{g}}) \leq \lambda(\L_{\tilde{g}}) \leq
K^{\frac{4}{n-2}} \lambda(\L_{\bar{g}}). 
$$
This completes the proof.
\end{Proof}
Recall that $Y^{\cyl}_{[\bar{g}]} (X) = -\infty$ when $\lambda({\cal
L}_h)<0$. Hence the finitness of $\lambda(\L_{\bar{g}})$ does not
imply that of the cylindrical Yamabe constant
$Y^{\cyl}_{[\bar{g}]}(X)$. However, the signs of
$\lambda(\L_{\bar{g}})$ and $Y^{\cyl}_{[\bar{g}]} (X) $ are still
related as follows.
\begin{Proposition}\label{app:P2}
Let $\bar{g}$ be a cylindrical metric on $X$. Then
\begin{enumerate}
\item[{\bf (i)}] $\lambda(\L_{\bar{g}}) \geq 0$ if and only if
$Y^{\cyl}_{[\bar{g}]} (X)\geq 0$,
\item[{\bf (ii)}] $\lambda(\L_{\bar{g}}) < 0$ if and only if
$Y^{\cyl}_{[\bar{g}]} (X)<0$.
\end{enumerate}
\end{Proposition}
\begin{Proof}
It is enough to prove only {\bf (ii)}.  We first note that
\begin{equation}\label{new78}
Q_{(X, \bar{g})}(f) = \frac{E_{(X,
\bar{g})}(f)}{||f||^2_{L^{\frac{2n}{n-2}}_{\bar{g}}}} = 
\frac{E_{(X,
\bar{g})}(f)}{||f||^2_{L^2_{\bar{g}}}} \cdot
\frac{||f||^2_{L^2_{\bar{g}}}}{||f||^2_{L^{\frac{2n}{n-2}}_{\bar{g}}}}
\end{equation}
for any $f \in C^{\infty}_c(X)$ with $f \not\equiv 0$.  This implies that
$\lambda(\L_{\bar{g}}) < 0$ if and only if $Y^{\cyl}_{[\bar{g}]}(X) <
0$.
\end{Proof}
It follows from (i) that $\lambda(\L_{\bar{g}}) \geq 0$ implies
$\lambda({\cal L}_h)\geq 0$ for any cylindrical metric $\bar{g}$ with
$\p_{\infty}\bar{g}=h$. 
\vspace{2mm}

For a closed Riemannian manifold $(M, g)$, the following result for
the first eigenvalue $\lambda(\L_g)$ of $\L_g$ is well-known
(cf. \cite{Au1}, \cite{KW}): 
\vspace{2mm}

\noindent
\emph{There exists a conformal metric
$\hat{g} \in [g]$ such that the scalar curvature $R_{\hat{g}}$
satisfies $R_{\bar{g}} > 0, R_{\bar{g}} < 0$, or $R_{\bar{g}} \equiv
0$ on $X$.  Moreover, the sign of $R_{\hat{g}}$ is identical with the
sign of $\lambda(\L_g)$.}
\vspace{2mm}

Under the assumption $\lambda({\mathcal L}_h) > 0$ for $(X, \bar{g})$,
a similar relationship holds.
\begin{Proposition}\label{app:P3}
Let $(X, \bar{g})$ be a cylindrical manifold of $\dim X = n \geq 3$
with tame ends $Z \times [0, \infty)$ and $h = \p_{\infty}\bar{g} \in
\Riem(Z)$ a metric with $\lambda({\mathcal L}_h) > 0$.  Then there
exists a conformal metric $\hat{g} = v^{\frac{4}{n-2}}\cdot \bar{g}$
with $v \in C^{\infty}_+(X) \cap L^{1,2}_{\bar{g}}(X)$ such that the
scalar curvature $R_{\hat{g}}$ satisfies $R_{\bar{g}} > 0, R_{\bar{g}}
< 0$, or $R_{\bar{g}} \equiv 0$ on $X$.  Moreover, the sign of
$R_{\hat{g}}$ is identical with the sign of $\lambda(\L_{\bar{g}})$.
\end{Proposition}
\begin{Proof}
By the conformal-rescaling argument below, we may assume that
$\lambda(\L_{\bar{g}}) < \lambda({\mathcal L}_h) $.  Indeed, if
$\lambda(\L_{\bar{g}}) \geq \lambda({\mathcal L}_h)$, it is enough to
change $\bar{g}$ to an appropriate conformal metric $\bar{g}_{\varphi}
= e^{2\varphi}\cdot \bar{g}$ with
$$
\varphi =
\{\begin{array}{cl}
k \equiv \textrm{const.} > 0 &
\textrm{on} \quad  W = X \setminus (Z \times (0, \infty)), 
\\
0 & \textrm{on}\quad Z \times [1, \infty).
\end{array}
\right.
$$
In particular, $\L_{\bar{g}_{\varphi}} = e^{-2k}\cdot \L_{\bar{g}}$ on
$W$ and $\bar{g}_{\varphi} = \bar{g}$ on $Z \times [1, \infty)$.  By
using this combined with $\lambda({\mathcal L}_h) > 0$, the Dirichlet
first eigenvalue $\lambda(\L_{\bar{g}_{\varphi}} ; W)$ on $W$
satisfies $\lambda(\L_{\bar{g}_{\varphi}} ; W) < \lambda({\mathcal
L}_h)$ for sufficiently large $k > 0$.  Hence by the domain
monotonicity of the Dirichlet eigenvalues,
$\lambda(\L_{\bar{g}_{\varphi}}) \leq \lambda(\L_{\bar{g}_{\varphi}} ;
W) < \lambda({\mathcal L}_h)$.  Using this and Proposition
\ref{app:P1}, we may assume that $\bar{g}$ itself satisfies
$\lambda(\L_{\bar{g}}) < \lambda({\mathcal
L}_h)$.  

Then the standard argument implies the existence of a positive
minimizer $v \in C^{\infty}_+(X) \cap L^{1,2}_{\bar{g}}(X)$ with
$\lambda(\L_{\bar{g}}) = E_{(X, \bar{g})}(v)/||v||^2_{L^2_{\bar{g}}}$.
Hence $v$ satisfies
$$
- \alpha_n\Delta_{\bar{g}}v + R_{\bar{g}}v
= \lambda(\L_{\bar{g}})v \quad \textrm{on} \ \ \ X,
$$
and thus $R_{\hat{g}} = \lambda(\L_{\bar{g}})v^{-\frac{4}{n-2}}$
for $\hat{g} = v^{\frac{4}{n-2}}\cdot \bar{g}$.
This completes the proof.
\end{Proof}
From Theorem \ref{yam:Th1} and Propositions \ref{app:P2}, \ref{app:P3}, we
obtain the following.
\begin{Corollary}\label{app:P4}
Under the same assumptions as in Proposition \ref{app:P3}, then
\begin{enumerate}
\item[{\bf (i)}] $\lambda(\L_{\bar{g}}) > 0$ if and only if
$Y^{\cyl}_{[\bar{g}]} (X) > 0$, 
\item[{\bf (ii)}] $\lambda(\L_{\bar{g}}) = 0$ if and
only if $Y^{\cyl}_{[\bar{g}]} (X) = 0$, 
\item[{\bf (iii)}] $\lambda(\L_{\bar{g}}) < 0$ if and only if
$Y^{\cyl}_{[\bar{g}]} (X) < 0$.
\end{enumerate}
\end{Corollary}
\begin{Proof}
From Proposition \ref{app:P2}, it is enough to prove that
$\lambda(\L_{\bar{g}}) = 0$ if and only if $Y^{\cyl}_{[\bar{g}]} (X) = 0$.

First we assume that $\lambda(\L_{\bar{g}}) = 0$. From the argument in
the proof of Proposition \ref{app:P3}, there exists a non-zero
function $v \in C^{\infty}_+(X) \cap L^{1,2}_{\bar{g}}(X)$ such that
$E_{(X, \bar{g})}(v)/||v||^2_{L^2_{\bar{g}}} = \lambda(\L_{\bar{g}}) =
0$.  Using this in (\ref{new78}), we then obtain
$$
Y^{\cyl}_{[\bar{g}]}(X) \leq Q_{(X, \bar{g})}(v) = 0.
$$
If $Y^{\cyl}_{[\bar{g}]}(X) < 0$, then $\lambda(\L_{\bar{g}}) < 0$
by Proposition \ref{app:P2}, (ii).  This is a contradiction.
Hence, $Y^{\cyl}_{[\bar{g}]}(X) = 0$.

Next we assume that $Y^{\cyl}_{[\bar{g}]}(X) = 0$.
Note that $Y^{\cyl}_{[\bar{g}]}(X) = 0 < \lambda({\mathcal L}_h)$.
From Theorem \ref{yam:Th1}, there exists a Yamabe minimizer
$u \in C^{\infty}_+(X) \cap L^{1,2}_{\bar{g}}(X)$ such that
$Q_{(X, \bar{g})}(u) = Y^{\cyl}_{[\bar{g}]}(X) = 0$.
Using this also in (\ref{new78}), we then obtain
$$
\lambda(\L_{\bar{g}}) \leq
\frac{E_{(X, \bar{g})}(u)}{||u||^2_{L^2_{\bar{g}}}} = 0.
$$
If $\lambda(\L_{\bar{g}}) < 0$, then $Y^{\cyl}_{[\bar{g}]}(X) < 0$
by Proposition \ref{app:P2}, (ii).  This is also a contradiction.
Hence, $\lambda(\L_{\bar{g}}) = 0$.
This completes the proof.
\end{Proof}

\vspace{10mm}

\begin{small}
{\sf
\noindent
Kazuo Akutagawa, Shizuoka University, Shizuoka, Japan
\\
e-mail: smkacta@ipc.shizuoka.ac.jp
\\
\\
Boris Botvinnik, University of Oregon, Eugene, USA
\\
e-mail: 
botvinn@math.uoregon.edu
}
\end{small}
\end{document}